\documentclass[10pt]{amsart}
\pdfoutput=1
\usepackage{amsmath, amsthm, amssymb,slashed}
\usepackage{ifpdf}
\usepackage[pdftex]{graphicx}
\usepackage{tikz}
\usetikzlibrary{matrix,arrows,calc}

\usepackage[pdftex,plainpages=false,hypertexnames=false,pdfpagelabels]{hyperref}
 \theoremstyle{plain}
 \newtheorem{thm}{Theorem}[section]
    
 \newtheorem{cor}[thm]{Corollary}
 \newtheorem{lem}[thm]{Lemma}
 \newtheorem{prop}[thm]{Proposition}
 \theoremstyle{definition}
 \newtheorem{defn}[thm]{Definition}
 \newtheorem{notation}[thm]{Notation}
 \newtheorem{ex}[thm]{Example}
  
 \newtheorem*{thm*}{Theorem 2.4}
 \theoremstyle{remark}
 \newtheorem{rmk}[thm]{Remark}

\def\beq{\begin{eqnarray}}
\def\eeq{\end{eqnarray}}
 \newcommand{\bp}{\begin{proof}[Proof]}
 \newcommand{\ep}{\end{proof}}

\DeclareSymbolFont{bbold}{U}{bbold}{m}{n}
\DeclareSymbolFontAlphabet{\mathbbold}{bbold}

\def\Fun{{\sf Fun}}
\def\dil{{u}}
\def\Eu{{\rm Eu}}
\def\twist{{\alpha}}
\def\Mst{{\mathcal{M}}}
\def\Ind{{\rm {Ind}}}

\def\Mell{\mathcal{M}_{\rm ell}}

\def\H{{\rm H}}
\def\HH{{\mathbb H}}
\def\Spin{{\rm Spin}}
\def\SO{{\rm SO}}
\def\U{{\rm U}}

\def\Map{{\sf Map}}

\def\TMF{{\rm TMF}}
\def\Ell{{\rm Ell}}
\def\EE{{\mathcal{E}}}

\def\pt{{\rm pt}}
\def\cl{{\rm cl}}
\def\F{\mathcal{F}}

\def\vol{{\rm vol}}

\def\ev{{\rm ev}}
\def\odd{{\rm odd}}

\def\R{{\mathbb{R}}}

\def\M{{\mathbb{M}}}
\def\E{{\mathbb{E}}}

\def\CP{{\mathbb{CP}}}
\def\N{{\mathbb{N}}}
\def\id{{{\rm id}}}
\def\K{{\rm {K}}}
\def\C{{\mathbb{C}}}

\def\Z{{\mathbb{Z}}}

\def\End{{\sf End}}
\def\Aut{{\sf Aut}}
\def\Bun{{\sf Bun}}

\def\SL{{\rm SL}}
\def\GL{{\rm GL}}

\def\Rep{{\rm Rep}}

\def\Iso{{\sf Iso}}

\def\Lat{{\sf Lat}}
\def\F{{\mathcal{F}}}

\newcommand{\op}{{\sf{op}}}   
\newcommand{\sq}{\mathord{/\!\!/}}

\def\twocommute{\ensuremath{\rotatebox[origin=c]{30}{$\Rightarrow$}}}

\newcommand\nc{\newcommand}
\nc\mf\mathfrak
\nc\mc\mathcal
\nc\mb\mathbb

\begin{document}

\title[Equivariant elliptic cohomology and discrete torsion]{Equivariant elliptic cohomology, gauged sigma models, and discrete torsion}

\author{Daniel Berwick-Evans}

\address{Department of Mathematics, University of Illinois at Urbana--Champaign}

\email{danbe@illinois.edu}


\begin{abstract}
For $G$ a finite group, we show that functions on fields for the 2-dimensional supersymmetric sigma model with background $G$-symmetry determine cocycles for complex analytic $G$-equivariant elliptic cohomology. Similar structures in supersymmetric mechanics determine cocycles for equivariant K-theory with complex coefficients. The path integral for gauge theory with a finite group constructs wrong-way maps associated to group homomorphisms. When applied to an inclusion of groups, we obtain the induced character formula of Hopkins, Kuhn, and Ravenel. For the homomorphism $G\to *$ we obtain Vafa's formula for gauging with discrete torsion. The image of equivariant Euler classes under gauging constructs modular form-valued invariants of representations that depend on a choice of string structure.  We illustrate nontrivial dependence on the string structure for a 16-dimensional representation of the Klein 4-group. 
\end{abstract}

\maketitle

\maketitle 
\setcounter{tocdepth}{1}
\tableofcontents

\section{Introduction}

The goal of this paper is to explain how super spaces of classical fields with background $G$-symmetry in dimensions 1 and 2 naturally lead to cocycle models for equivariant K-theory and equivariant elliptic cohomology with complex coefficients, respectively. In an effort to make this connection as accessible and explicit as possible, we focus on the easiest case of a finite group acting on a smooth manifold. On the topological side, the finite group assumption simplifies the construction and calculation of complexified equivariant K-theory and complex analytic equivariant elliptic cohomology. On the physics side, this assumption simplifies the space of classical fields. Both simplifications arise because the relevant objects are parameterized by moduli spaces of $G$-bundles, and when~$G$ is finite this parameterizing space is discrete. 

The context of finite groups also also permits rigorous construction of certain path integrals as (finite) sums over principal bundles~\cite{DW,FreedQuinn}. This determines wrong-way maps in cohomology. In particular, we realize Vafa's \emph{discrete torsion}~\cite{Vafatorsion} as a structure internal to elliptic cohomology. We also show that the path integral constructs the Hopkins--Kuhn--Ravenel  formula for (higher) induced characters~\cite{HKR}. These wrong-way maps are highly computable, and we provide several illustrative examples. 
In particular we construct modular form-valued invariants of a 16-dimensional representation of $\Z/2\times \Z/2$ that distinguish the two possible choices of string structure on $\R^{16}\sq (\Z/2\times \Z/2)$.

\subsection{Statement of results}
Let $M$ be a compact manifold with an action by a finite group~$G$. For $d=1,2$ we construct a groupoid of fields $\mathcal{L}^{d|1}(M\sq G)$ with objects
$$
{\rm Obj}(\mathcal{L}^{d|1}(M\sq G))=\{(\ell,{\bf g},\phi)\mid T^{d|1}_\ell\leftarrow P_{\bf g}\stackrel{\phi}{\to} M\}
$$
where $\ell$ is a conformal structure on $T^{d|1}=\R^{d|1}/\Z^d$, ${\bf g}\colon \Z^d\to G$ determines a principal $G$-bundle $P_{\bf g}\to T^{d|1}_\ell$, and $\phi$ is a $G$-equivariant map. For $d=1,2$, the groupoids $\mathcal{L}^{d|1}(M\sq G)$ support line bundles $\omega^{k/2+\twist}$ for $k\in \Z$ and $[\twist]\in \H^{d+1}(BG;U(1))$. There is a finite-dimensional smooth subgroupoid $\mathcal{L}_0^{d|1}(M\sq G)\subset \mathcal{L}^{d|1}(M\sq G)$, and we identify sections of~$\omega^{k/2+\twist}$ on this subgroupoid with cocycles for a cohomology theory as follows.

\begin{thm}\label{mainthm}
There are cocycle maps 
\beq
\Gamma(\mathcal{L}^{1|1}_0(M\sq G);\omega^{\bullet/2+\twist})&\twoheadrightarrow& \K_G^{\bullet+\twist}(M)\otimes \C \label{eq:Kthy}\\
\Gamma(\mathcal{L}^{2|1}_0(M\sq G);\omega^{\bullet/2+\twist})&\twoheadrightarrow& \TMF_G^{\bullet+\twist}(M)\otimes \C \label{eq:TMF}
\eeq
realizing~\eqref{eq:Kthy} as a model for the $\twist$-twisted $G$-equivariant K-theory of $M$ with complex coefficients and~\eqref{eq:TMF} as a model for (the global sections of) the complex analytic $\twist$-twisted $G$-equivariant elliptic cohomology of $M$. 
\end{thm}

We also refine~\eqref{eq:Kthy} and~\eqref{eq:TMF} to morphisms of sheaves of commutative differential graded algebras, recovering Theorem~\ref{mainthm} as the map on (derived) global sections; see Theorem~\ref{mainthmsheaf}. 


The subgroupoid $\mathcal{L}^{1|1}_0(M\sq G)\subset \mathcal{L}^{1|1}(M\sq G)$ is a superspace generalization of the \emph{inertia groupoid}, $\Map(\pt\sq \Z,M\sq G)\subset \Map(\R\sq \Z,M\sq G)$. 
It plays the role of the constant loops in the action groupoid $M\sq G$ with extra structure coming from nontrivial automorphisms. Similarly, the super Lie groupoid $\mathcal{L}^{2|1}_0(M\sq G)$ is a superspace generalization of the \emph{double inertia groupoid}, $\Map(\pt\sq \Z^2,M\sq G)\subset \Map(\R^2\sq \Z^2,M\sq G)$.
There is an $\SL_2(\Z)$-action on the double inertia groupoid from the precomposition action on $\pt\sq \Z^2$. Geometrically, this comes from the action of the mapping class group on the 2-torus $T^2\simeq \R^2\sq \Z^2$. The modularity properties of elliptic cocycles under~\eqref{eq:TMF} are inherited from this $\SL_2(\Z)$-action.

We explore two applications of Theorem~\ref{mainthm}. First, physical theories construct cocycle representatives of cohomology classes via their partition functions; we focus on the case of free fermion theories, whose cohomological counterparts are (twisted, equivariant) Euler classes as explained after Corollary~\ref{cor:main1} below. Second, operations on field theories determine operations on the associated cohomology theories. The operation we consider is a version of quantization of a finite gauge theory, which leads to a wrong-way map in cohomology. 
Path integral quantization of a gauge theory is (formally) an integral over the space of connections modulo gauge. Although this is generally fraught with technical difficulties, when~$G$ is finite there is a well-known measure~\cite{DW,FreedQuinn} that permits rigorous construction of the path integral as a finite sum over isomorphism classes of principal bundles. 

\begin{thm}[Theorem~\ref{prop:pushforward}]\label{thm:pushforward}
Given a homomorphism $\zeta\colon H\to G$, a $G$-manifold $M$, and a $d$-cocycle representing a class $[B]\in \H^d(BG;U(1))$, there exists a wrong-way map
$$
\zeta_!^B\colon \Gamma(\mathcal{L}^{d|1}_0(M\sq H);\omega^{\bullet/2+\zeta^*\alpha})\to \Gamma(\mathcal{L}^{d|1}_0(M\sq G);\omega^{\bullet/2+\alpha})
$$
determined by the sum over principal bundles
\beq
(\zeta_!f)(\ell,{\bf g},\phi)=\sum_{[P_{\bf g}\simeq P_{\bf h}\times_\zeta G]} \frac{(\epsilon_B\cdot f)(\ell,{\bf h},\phi')}{|{\rm Aut}(P_{\bf h})|/|{\rm Aut}(P_{\bf g})|}\label{eq:induction}
\eeq
indexed by the finite set of isomorphism classes of $H$-principal bundles~$P_{\bf h}$ with an isomorphism of $G$-bundles $P_{\bf g}\simeq P_{\bf h}\times_\zeta G$ that is compatible with the map to $M$. Here, ${\rm Aut}(P_{\bf h})$ is the automorphism group of the $H$-bundle $P_{\bf h}$, $P_{\bf h}\times_\zeta G$ is the associated $G$-bundle for the left $H$-action on $G$ via $\zeta$, and the function $\epsilon_B$ is determined by $B$ via transgression. 
\end{thm}

Theorem~\ref{mainthm} together with a little computation gives the following. 

\begin{cor}[Corollaries~\ref{cor:inductionexists},~\ref{cor:HKR} and~\ref{cor:Vafa}]\label{cor:main1}
The formula~\eqref{eq:induction} determines wrong-way maps
$$
\zeta_!^B\colon \K_H^{\bullet+\zeta^*\alpha}(M)\to \K_G^{\bullet+\alpha}(M),\qquad \zeta_!^B\colon \TMF_H^{\bullet+\zeta^*\alpha}(M)\otimes \C\to \TMF_G^{\bullet+\alpha}(M)\otimes \C. 
$$
When $M=\pt$, $\zeta\colon H\hookrightarrow G$ is an inclusion, and $B\equiv 1$, this recovers the height~$d=1,2$ induction formulas of Hopkins--Kuhn--Ravenel~\cite[Theorem~D]{HKR}. When $\zeta\colon G\to \{e\}$ and $d=2$, it recovers the formula for Vafa's gauging with discrete torsion $\epsilon_B$. 
\end{cor}

The pushforward determined by~\eqref{eq:induction} is highly computable, and we provide some physics-inspired applications with a focus on~$d=2$. Namely, pushing forward equivariant Euler classes of $G$-representations with string structure along~$\zeta\colon G\to \{e\}$ constructs invariants of such representations values in weak modular forms. We compute this explicitly for examples built from the sign representation~$\sigma:=\R\sq \Z/2$, i.e., $\Z/2$ acting on $\R$ by $\{\pm 1\}$. Specifically, we consider $\Eu(\sigma^{\oplus k})$ for the $\Z/2$-equivariant class associated with $k$~copies of the sign representation and $\Eu(\sigma^{\oplus k})^{\boxtimes j}$ for the $(\Z/2)^{\times j}$-equivariant class associated with the action on $(\R^{k})^{\times j}$ where each $\Z/2$ factor acts by the sign representation on $\R^k$. Then we compute the values of the indicated wrong-way maps (see Examples~\ref{ex:Jacobi}-\ref{ex:discretetorsion}), 
\beq
&&\resizebox{.95\textwidth}{!}{$
\begin{array}{rll}
\TMF_{\Z/2}^\bullet(\pt)\otimes \C\to \TMF^\bullet(\pt)\otimes \C:& \zeta_!(\Eu(\sigma^{\oplus 8}))=0, & \zeta_!( \Eu(\sigma^{\oplus 16}))=\frac{\vartheta_{E_8}}{\Delta},\\ 
\quad \TMF_{\Z/2\times \Z/2}^\bullet(\pt)\otimes \C\to \TMF^\bullet(\pt)\otimes \C:&  \zeta_!(\Eu(\sigma^{\oplus 8})^{\boxtimes 2})=0,& \zeta_!^B(\Eu(\sigma^{\oplus 8})^{\boxtimes 2})=\frac{\vartheta_{E_8}}{\Delta},\\
\TMF_{(\Z/2)^{\times 3}}^\bullet(\pt)\otimes \C\to \TMF^\bullet(\pt)\otimes \C:& \zeta_!(\Eu(\sigma^{\oplus 16})^{\boxtimes 3})=\frac{j}{\Delta^2}.
\end{array}$}\label{eq:examples}
\eeq
Above, $\Delta$ denotes the modular discriminant, $\vartheta_{E_8}$ is the theta function of the $E_8$ lattice (which also equals the 4th Eisenstein series), and $j$ is Klein's $j$-function. The computations~\eqref{eq:examples} boil down to identities amongst Jacobi theta functions. They depend on the fact that the representation~$\sigma^{\oplus 8}$ has a canonical choice of string structure, which endows $(\sigma^{\oplus 8})^{\boxtimes j}$ and $\sigma^{\oplus 8k}$ with string structures. The wrong-way maps denoted $\zeta_!$ above use this canonical string structure. There is a second choice of string structure on $(\sigma^{\oplus 8})^{\boxtimes 2}$ from the fact that $\H^2(B(\Z/2\times \Z/2);\Z/2)=\Z/2$, and $\zeta_!^B$ is the pushforward with respect to this other choice. Since in the middle line $\frac{\vartheta_{E_8}}{\Delta}\ne 0$, the modular form-valued invariant of the representation $(\sigma^{\oplus 8})^{\boxtimes 2}=\R^{16}\sq (\Z/2\times \Z/2)$ depends on the choice of string structure. 

The equivariant Euler classes in~\eqref{eq:examples} have an interpretation as partition functions of free fermion theories. In terms of theta functions, this is a well-known fact in the physics literature: the Klein form that defines the equivariant Euler class literally \emph{is} the formula for the partition function of the associated physical theory, e.g., see~\cite{Vafatorsion}. The precise connection with equivariant elliptic cohomology---as well as a rigorous construction of these partition functions---is given in~\cite[\S7]{BET1}. Hence, from the point of view of physics the equalities~\eqref{eq:examples} compute partition functions of specific gauged free fermion theories, where the choice of string structure is a choice of \emph{anomaly cancelation}. 

In the mathematically well-understood examples relating elliptic cohomology and 2-dimensional field theories, the choice of anomaly cancelation does not affect the resulting elliptic invariant. For example, the Witten genus of an oriented manifold~$M$ is the $q$-expansion of a modular form when $M$ has a rational string structure~\cite{Witten_Dirac}. We emphasize that modularity in this case only depends on the existence of a rational string structure, i.e., the vanishing of $\frac{p_1}{2}(TM)$ in rational cohomology. It does not require an integral vanishing of $\frac{p_1}{2}(TM)$ nor does the resulting invariant depend on a choice of coboundary $\frac{p_1}{2}(TM)=dH$. Recent work has constructed (at a physics level of rigor) torsion invariants of string manifolds where the invariant depends on a \emph{choice} of integral string structure, manifesting in the physics as a choice of anomaly cancelation~\cite{GJFW,GJF2}. Furthermore, these torsion invariants have an interpretation in terms of topological modular forms (TMF). Hence, a deeper mathematical understanding of anomaly cancelation and string structures seems to be a key step toward understanding how the rich torsion in TMF might be captured by invariants of 2-dimensional supersymmetric field theories. The discrete torsion invariants constructed by Theorem~\ref{thm:pushforward} are one piece of this puzzle: the second line in~\eqref{eq:examples} gives an example where the invariant in elliptic cohomology depends on the choice of integral string structure.

\subsection{Relationship with other work}

The proposed relationship between gauged sigma models and equivariant elliptic cohomology is an old one, making it difficult to comprehensively describe previous work on the problem. In our view there are two important strands of thought. On the one hand, physicists observed connections between equivariant modular objects (e.g., theta functions) and 1-loop partition functions of orbifolds. The most relevant reference for the work below is Vafa's construction of discrete torsion invariants~\cite{Vafatorsion}, which draws on~\cite{DHVW1,Theta}. Around the same time, mathematicians began to recognize analogies between $d$-dimensional gauge theories and height $d$ chromatic phenomena. For example Segal noticed~\cite[\S5]{Segal_Elliptic} that formulas in 2-dimensional gauge theory resembled those in height~2 Hopkins--Kuhn--Ravenel character theory~\cite{HKR}. Many authors elaborated on this analogy (albeit decades later); we were particularly inspired by Morava~\cite{MoravaHKR,Morava}, Ando--French~\cite{AndoFrench}, Ganter~\cite{GanterHecke}, Stolz--Teichner~\cite{ST04,ST11}, and Barthel--Stapleton~\cite{Stapleton,BarthelStapleton}, who explored beautiful and striking similarities between structures in physics and ones in homotopy theory. 

The main contribution of this work is to provide a \emph{direct} link between these two strands of thought from physics and mathematics. Indeed, starting from the superspace formalism of supersymmetric sigma models, Theorem~\ref{mainthm} constructs cocycle models in terms of functions on a certain subspace of fields. Theorem~\ref{thm:pushforward} further shows that operations on field theories construct non-trivial operations on the corresponding cohomology theories. This parallels the nonequivariant cocycle model and construction of pushforwards in~\cite{DBE_MQ}, further cementing the connection between height~2 homotopy theory and 2-dimensional physics. Just as in the nonequivariant cocycle models from~\cite{DBE_MQ}, the constructions below are designed to communicate with the more ambitious goals of the Stolz--Teichner program to realize supersymmetric field theories as cocycles for cohomology theories. Indeed, the fields $\mathcal{L}^{d|1}_0(M\sq G)$ determine objects in the bordism categories appearing in the proposed model for equivariant K-theory and equivariant topological modular forms~\cite{ST11,BET1}. In particular, partition functions in this framework determine sections of line bundles on the super Lie groupoids $\mathcal{L}^{d|1}_0(M\sq G)$. By Theorem~\ref{mainthm}, this determines cocycles in complexifications of the expected (equivariant) cohomology theories.\footnote{In the $2|1$-dimensional case there is a subtlety regarding the holomorphic dependence of these sections on the moduli of elliptic curves. This is the subject of~\cite{DBEChern}; see Remarks~\ref{rmk:supercircle} and~\ref{rmk:holophys} below.}

The super Lie groupoids $\mathcal{L}^{d|1}_0(M\sq G)$ in Theorem~\ref{mainthm} can also be studied in the language of super stacks. This is the approach taken in~\cite{Stoffel} when $d=1$ and~\cite{BET1} and~\cite{Powerops} when $d=2$. Indeed, the geometric power operation in~\cite{Powerops} is only well-defined on the stack underlying the super Lie groupoid~$\mathcal{L}^{d|1}_0(M\sq G)$, and various manipulations in~\cite{Stoffel} and~\cite{BET1} are similarly only well-defined on the stack. In this paper, we have deliberately avoided super stacks both because they are not widely-known and also because they are unnecessary for our intended applications. Indeed, precisely when~$G$ is finite the stacks associated to $\mathcal{L}^{d|1}_0(M\sq G)$ are \emph{geometric}, meaning they are presented by super Lie groupoids; see Remark~\ref{rmk:fail}. Hence, although some of the results below can be deduced (with effort) from results in~\cite{BET0,BET1}, the methodology below is likely to be more accessible.

Not only does the finite group assumption provide a less strenuous path to Theorem~\ref{mainthm}, but it also permits constructions and computations unique to this setting. For example, the derived global sections of the complex analytic equivariant elliptic cohomology sheaf is still unknown in general (see~\cite[Remark~3.7]{BET0}); for finite groups we give a complete answer in Proposition~\ref{prop:dgs}. Gauging and discrete torsion (Theorem~\ref{thm:pushforward}) are also only possible in the finite group case, arising as a (finite) sum over principal $G$-bundles. The equivariant Euler classes for finite groups are furthermore closely related to Klein forms and Jacobi's theta functions, permitting computations that become significantly more technical (if not ill-posed) in the case of a general compact Lie group. For example, a representation of a connected Lie group admits a string structure if and only if the representation is trivial. Hence invariants of representation depending on a choice of string structure as in~\eqref{eq:examples} are only interesting in the case that $G$ is finite. 

\subsection{Outline} Section~\ref{sec:groupoids} provides background for two standard operations in the 2-category of (finite) groupoids: transgression and groupoid integration. Section~\ref{sec:cohomo} defines complexified equivariant K-theory and complex analytic equivariant elliptic cohomology as derived global sections of sheaves on certain moduli spaces of $G$-bundles. These definitions seem to be known to the experts, but have not appeared in print. Section~\ref{sec:super} gives a brief introduction to supermanifolds and super Lie groupoids. Section~\ref{sec:circlesandtori} constructs groupoids of $G$-bundles on super conformal circles and super conformal tori. Section~\ref{sec:cocycles} defines inertia fields $\mathcal{L}^{d|1}_0(M\sq G)$ and proves they define a super Lie groupoid. Section~\ref{sec:coccyles} constructs the cocycle models and proves Theorem~\ref{mainthm} as well as the sheaf-theoretic refinement Theorem~\ref{mainthmsheaf}. Section~\ref{sec:discretetorsion} constructs the finite path integral, proves Theorem~\ref{thm:pushforward}, and computes the examples~\eqref{eq:examples}. 


\subsection{Acknowledgements}
It is a pleasure to thank Tobi Barthel, Nora Ganter, Theo Johnson-Freyd, Charles Rezk, Nat Stapleton, and Arnav Tripathy for the generously sharing their numerous insights into this subject. 

\section{Inertia groupoids, transgression, and groupoid integration}\label{sec:groupoids}

This section reviews constructions of inertia groupoids, transgression maps, and integration on finite groupoids. This material is all standard, though is a bit spread out in the literature. In~\S\ref{sec:HKR} we show that integration on inertia groupoids recovers the induction formula of Hopkins--Kuhn--Ravenel~\cite[Theorem~D]{HKR}. This connection is surely known to the experts, but we do not know of a reference. 

\subsection{Groupoids}

A \emph{groupoid} is a (small) category whose morphisms are all invertible. The data of a groupoid is hence a set of objects $\mathcal{G}_0$, a set of morphisms $\mathcal{G}_1$ and maps
\beq
s,t\colon \mathcal{G}_1\to \mathcal{G}_0,\quad u\colon \mathcal{G}_0\to \mathcal{G}_1, \quad c\colon \mathcal{G}_1\times_{\mathcal{G}_0}\mathcal{G}_1\to \mathcal{G}_1,\label{eq:grpd}
\eeq
called source, target, unit, and composition, respectively. We often abbreviate this, writing a groupoid as $\mathcal{G}=\{\mathcal{G}_1\rightrightarrows \mathcal{G}_0\}$. There is a strict 2-category whose objects are groupoids, 1-morphisms are functors, and 2-morphisms are natural isomorphisms of functors. A groupoid is \emph{finite} if $\mathcal{G}_0$ and $\mathcal{G}_1$ are finite sets. 


\begin{ex} \label{ex:action} Suppose a group $G$ acts on a set~$X$ on the left. The \emph{action groupoid} $X\sq G$ has objects $\mathcal{G}_0=X$ and morphisms $\mathcal{G}_1=G\times X$. The source map is the projection $s=p_2$, and the target map is the action. The unit map is inclusion along the identity element of~$G$, and composition is determined by multiplication in $G$. If $G$ is a finite group and $X$ is a finite set, $X\sq G$ is a finite groupoid. 
\end{ex}

\begin{ex} \label{ex:discrete} Given a set $X$, the \emph{discrete groupoid on $X$} has objects $X$ and only identity morphisms. Conversely, a groupoid $\mathcal{G}$ is \emph{discrete} if its only morphisms are identities.
\end{ex}

\begin{ex}
For groupoids $\mathcal{G}$ and $\mathcal{H}$, $\Fun(\mathcal{G},\mathcal{H})$ is the groupoid of functors and natural isomorphisms of functors. If $\mathcal{G}$ and $\mathcal{H}$ are finite groupoids, then so is $\Fun(\mathcal{G},\mathcal{H})$. 
\end{ex} 

\begin{defn}\label{defn:inertia} Given a groupoid $\mathcal{G}$, the \emph{$d$-fold inertia groupoid} of $\mathcal{G}$ is $\Fun(\pt\sq \Z^d,\mathcal{G})$. \end{defn}

The $d$-fold inertia groupoid carries an action by $\SL_d(\Z)$, acting by precomposing with automorphisms of $\pt\sq \Z^d$ coming from the standard action of $\SL_d(\Z)$ on $\Z^d$. 


\begin{ex}
The $d$-fold inertia groupoid of $\pt\sq G$ has the explicit description 
$$
\Fun(\pt\sq \Z^d,\pt\sq G)\simeq G^{(d)}\sq G
$$
where $G^{(d)}\subset G^{\times d}$ denotes $d$-tuples of pairwise commuting elements in $G$ and $G^{(d)}\sq G$ is the action groupoid for the left $G$-action by simultaneous conjugation on $d$-tuples in $G^{(d)}$. If we include the action of $\SL_d(\Z)$, we obtain the groupoid $G^{(d)}\sq (G\times \SL_d(\Z))$, for the following left $\SL_d(\Z)$-action on $G^{(d)}$. View a $d$-tuple of commuting elements as defining a homomorphism $\Z^d\to G$, and precompose with the action by $\SL_d(\Z)$ on $\Z^d$. Note the precomposition actions are naturally a right action, and so to make this into a left action an element $\gamma\in \SL_d(\Z)$ acts through~$\gamma^{-1}$; see~\eqref{eq:G2act} for a formula when $d=2$. The $G$ and $\SL_d(\Z)$-actions on $G^{(d)}$ commute, yielding the action groupoid $G^{(d)}\sq (G\times \SL_d(\Z))$.
\end{ex}

\begin{ex}[$G$-bundles on tori]\label{ex:tori} Let $T^d= \R^d/\Z^d$ denote the $d$-dimensional torus and let $G$ be a finite group. Define a groupoid $\Bun_G(T^d)$ whose objects are principal $G$-bundles $P\to T^d$ and whose morphisms are isomorphisms of principal $G$-bundles
\beq
\begin{tikzpicture}[baseline=(basepoint)];
\node (A) at (0,0) {$P$};
\node (B) at (3,0) {$P'$};
\node (C) at (0,-1) {$T^d$};
\node (D) at (3,-1) {$T^d$};
\draw[->] (A) to node [above] {$\sim$}  (B);
\draw[->] (A) to (C);
\draw[->] (B) to (D);
\draw[->] (C) to node [above] {$\gamma$} (D);
\path (0,-.75) coordinate (basepoint);
\end{tikzpicture}\nonumber
\eeq
covering an automorphism of the torus determined by $\gamma\in \SL_d(\Z)$ acting on~$\R^d$ preserving the image $\Z^d\subset \R^d$. We caution that $\Bun_G(T^d)$ does not have a set of objects and morphisms, so is not small; however it is equivalent to a small groupoid, so is \emph{essentially small}.
%
Indeed, define a functor $G^{(d)}\sq (G\times \SL_d(\Z))\to {\Bun}_G(T^d)$ by the assignments
\beq&&\resizebox{.93\textwidth}{!}{
${\bf g}\mapsto P_{\bf g}=(G\times \R^d)/\Z^d\qquad ( (h,\gamma) \colon {\bf g}\to {\bf g}')\mapsto P_{\bf g}=(G\times \R^d)/\Z^d\stackrel{h\cdot}{\to} (G\times \R^d)/\Z^d=P_{\bf g'}$}. \label{eq:functorT}
\eeq
In more detail, the object ${\bf g}=(g_1,\dots, g_d)\in G^{(d)}$ is sent to the principal $G$-bundle $P_{\bf g}$ gotten by the quotient of the left $\Z^d$-action on~$G\times\R^d$ from the standard action of $\Z^d$ on $\R^d$ and the $\Z^d$-action on $G$ generated by $\{g_1,\dots, g_d\}$. To a morphism $(h,\gamma)\in G\times \SL_2(\Z)$, left multiplication by $h\in G$ and the action of $\SL_d(\Z)$ on ${\bf g}$ gives a morphism $P_{\bf g}\to P_{\bf g'}$ covering the automorphism of the torus. Identifying a $G$-bundle with a representation of $\pi_1(T^d)$ and choosing an isomorphism $\pi_1(T^d)\simeq \Z^d$ (i.e., generators of the fundamental group) one finds that this functor determined by~\eqref{eq:functorT} is essentially surjective. Since isomorphic $G$-bundles correspond to conjugate homomorphisms and the mapping class group $\SL_d(\Z)$ acts by changing the choice of generators under the isomorphism $\pi_1(T^d)\simeq \Z^d$, the functor is also fully faithful. Hence $\Bun_G(T^d)$ is equivalent to $G^{(d)}\sq (G\times \SL_d(\Z))$. 
\end{ex}

\begin{defn} 
A \emph{function on a groupoid} is an object in the groupoid $\Fun(\mathcal{G},\C)$, where $\C$ is viewed as a discrete groupoid. 
\end{defn}

Note the only natural isomorphisms between functors in $\Fun(\mathcal{G},\C)$ are identity isomorphisms, so that $\Fun(\mathcal{G},\C)$ is a discrete groupoid determined by its underlying set. The algebra structure on $\C$ gives $\Fun(\mathcal{G},\C)$ the structure of a $\C$-algebra. 

\begin{notation}
To fit with the (super) Lie groupoid notation later in the paper, let $C^\infty(\mathcal{G})=\Fun(\mathcal{G},\C)$ denote the algebra of functions on a finite groupoid. This is consistent because functions on a finite set are the same as smooth functions on the associated 0-manifold. 
\end{notation}

\begin{defn} Given a groupoid $\mathcal{G}$, let $\pi_0(\mathcal{G})$ denote the set of isomorphism classes of objects of $\mathcal{G}$. \end{defn}

\begin{ex}[Generalized class functions]
A function $f\in C^\infty(\mathcal{G})$ assigns a complex number $f(g)$ to each object $g\in \mathcal{G}$ with the property that isomorphic objects $g$ and $g'$ are assigned the same number, $f(g)=f(g')$. This is equivalent to a function on $\pi_0(\mathcal{G})$. Hence, functions on the $d$-fold inertia groupoid $\Fun(\pt\sq \Z^d,\pt\sq G)$ are functions on conjugacy classes of $d$-tuples of commuting elements in $G$. We identify this with the ring of class functions when $d=1$. For $d>1$, these are \emph{generalized class functions}~\cite[\S1]{Stapleton}. 
\end{ex}

\begin{defn} A \emph{line bundle} on a finite groupoid $\mathcal{G}$ is a complex line bundle on $\mathcal{G}_0$ (viewed as a topological space with the discrete topology) with an equivariant structure for the action by $\mathcal{G}_1$. \end{defn}

To unwind this definition, note first that a line bundle $L\to \mathcal{G}_0$ on a discrete space is always trivializable,~$L\simeq \mathcal{G}_0\times \C$. Then for each $x\in \mathcal{G}_1$ we require the data of an isomorphism $\C\simeq L_{s(x)}\to L_{t(x)}\simeq \C$ between the fibers of the line at the source and target $s(x),t(x)\in \mathcal{G}_0$ of $x$, and these isomorphisms need to be compatible with composition in $\mathcal{G}$. An isomorphism of lines is an element of $\C^\times$. Hence, a line bundle is precisely the data of a functor $\mathcal{G}\to\pt\sq \C^\times$, or equivalently a 1-cocycle on $\mathcal{G}$ with values in $\C^\times$. 


\subsection{Transgression for finite groupoids}
Let $A$ be an abelian group. Groupoid homology is defined as the total homology of the double complex $(C_*(\mathcal{G}_\bullet;A),\partial,d)$ where $\mathcal{G}_\bullet$ denotes the nerve of $\mathcal{G}$ and $C_*$ is the singular chains functor where the set $\mathcal{G}_k$ is viewed as a space with the discrete topology\footnote{This definition of (co)homology can be applied to any topological groupoid $\mathcal{G}$.}. The differentials $\partial$ and $d$ are the boundary operator and bar differential, respectively. Groupoid cohomology is defined similarly, applying cochains rather than chains. Groupoid (co)homology coincides with the (co)homology of the classifying space $B\mathcal{G}=|\mathcal{G}_\bullet|$ defined as the geometric realization of the nerve. We refer to~\cite[pg.~280-288]{BehrendCoh} for details. 

For a finite groupoid $\mathcal{G}$, consider the span
\beq
\Fun(\pt\sq \Z^d,\mathcal{G})\stackrel{p}{\leftarrow}  \pt\sq \Z^d \times\Fun(\pt\sq \Z^d,\mathcal{G})\stackrel{\ev}{\to} \mathcal{G}\label{eq:span}
\eeq
where $\Fun(\pt\sq \Z^d,\mathcal{G})$ is the finite groupoid whose objects are functors and morphisms are natural transformations of functors. In~\eqref{eq:span}, the map $\ev$ denotes the evaluation functor (adjoint to the identity functor) and $p$ is the projection. Fix an abelian group $A$. After choosing a representative of the fundamental class of the $d$-dimensional torus, $[T^d]\in \H_d(T^d;A)\simeq \H_d(B\Z^d;A)\simeq \H_d(\pt\sq \Z^d;A)$, define \emph{transgression} as the map on cocycles given by the composition 
$$
Z^k(\mathcal{G};A)\stackrel{\ev^*}{\to} Z^k(\pt\sq \Z^d \times\Fun(\pt\sq \Z^d,\mathcal{G});A)\stackrel{\int_{T^d}}{\to} Z^{k-d}(\Fun(\pt\sq \Z^d,\mathcal{G});A)
$$
where $\int_{T^d}$ denotes the pairing with the chosen cycle whose underlying class is the image of~$[T^d]$. The transgression map sends a $k$-cocycle on $\mathcal{G}$ to a $(k-d)$-cocycle on $\Fun(\pt\sq \Z^d,\mathcal{G})$. 

\begin{rmk} With $A=U(1)$, transgression gives a map from $(k-1)$-gerbes on $\mathcal{G}$ to $(k-d-1)$-gerbes on $\Fun(\pt\sq \Z^d,\mathcal{G})$, where a $0$-gerbe is a line bundle and a $(-1)$-gerbe is a $U(1)$-valued function. We refer to~\cite{Willerton} for a detailed description of this perspective. 
\end{rmk}

We will require the following specific instances of the transgression map. 

\begin{ex}[Transgression to loop space] \label{1transgr}
A 1-cocycle $B\colon G\to U(1)$ is a homomorphism, i.e., a 1-dimensional unitary representation. Transgression of this 1-cocycle determines an element of $C^\infty(G\sq G)=C^\infty(G)^G$, which is simply the 1-dimensional representation viewed as a character. Given a 2-cocycle $\twist\colon G\times G\to U(1)$, trangression returns the 1-cococyle on $\Fun(\pt\sq \Z,\pt \sq G)\simeq G\sq G$. A 1-cocycle is the same as a functor $G\sq G\to \pt\sq U(1)$, which we identify with a $G$-equivariant line bundle on $G$ for the conjugation action. A formula for this 1-cocycle is (e.g., see~\cite[\S2]{Willerton})
\beq
\mathcal{T}^\alpha\colon G\times G\to U(1),\qquad (g,h)\mapsto \frac{\twist(ghg^{-1},g)}{\twist(g,h)}\in U(1). \label{eq:11twist}
\eeq
Explicitly, we identify the cocycle $\mathcal{T}^\alpha$ with the equivariant line bundle on~$G\sq G$ whose line bundle on $G$ is trivial, and to a morphism $h\stackrel{g}{\to} ghg^{-1}$ in $G\sq G$ associates the linear map $\C\to \C$ between fibers of the line bundle specified by $\mathcal{T}^\alpha(g,h)$. 
\end{ex}

\begin{ex}[Transgression to double loops] \label{2transgr}
Given a 2-cocycle $B\colon G\times G\to U(1)$, double transgression gives a $G$-invariant $U(1)$-valued function on $\Fun(\pt\sq \Z^2,\pt\sq G)\simeq G^{(2)}\sq G$. A formula for this function is
$$
\epsilon_B(g_1,g_2)=\frac{B(g_1,g_2)}{B(g_2,g_1)}\in C^\infty(G^{(2)})^G\simeq C^\infty(\Fun(\pt\sq \Z^2,\pt\sq G)). 
$$
We note that the span~\eqref{eq:span} when $d=2$ is in fact $\SL_2(\Z)$-equivariant for the $\SL_2(\Z)$-action on the groupoid $\pt\sq \Z^2$. As a consequence, the function $\epsilon_B$ is $\SL_2(\Z)$-invariant, so descends to a function on the groupoid $G^{(2)}\sq G\times\SL_2(\Z)$ for the action
\beq
&&(\gamma,h)\cdot (g_1,g_2)= (hg_1^{d}g_2^{-b}h^{-1},hg_1^{-c}g_2^ah^{-1})\quad h\in G, \gamma=\left[\begin{array}{cc} a & b \\ c & d\end{array}\right]\in \SL_2(\Z). \label{eq:G2act}
\eeq
Given a 3-cocycle $\alpha\colon G^{\times 3}\to U(1)$, transgression determines a 1-cocycle on $G^{(2)}\sq G$. Moreover, this 1-cocycle is $\SL_2(\Z)$-equivariant, meaning it extends uniquely to a 1-cocycle
\beq
\mathcal{T}^\twist\colon G^{(2)}\times G\times \SL_2(\Z)\to U(1).\label{eq:21twist}
\eeq
Explicit formulas for this cocycle can be computed, e.g., see~\cite[\S2.3-2.4]{GanterHecke} and~\cite[\S2]{Willerton}. We identify $\mathcal{T}^\alpha$ with a line bundle on $G^{(2)}\sq (G\times \SL_2(\Z))$. 
\end{ex}


\begin{rmk}\label{rmk:projchar}
For a 2-cocycle $\alpha$, global sections of $\mathcal{T}^\twist$ over $G\sq G$ are spanned by characters of irreducible $\twist$-projective representations, i.e., projective representations determined by a map $\rho\colon G\to {\rm End}(V)$ with the properties $\rho(gh)=\twist(g,h)\rho(g)\rho(h)$ and $\rho(e)={\rm id}$, e.g., see~\cite[Lemma~5]{Willerton}. Similarly, for a 3-cocycle~$\alpha$, sections of $\mathcal{T}^\twist$ over $G^{(2)}\sq (G\times \SL_2(\Z))$ are spanned by characters of representations of the \emph{2-group} or \emph{categorical group}~$G^\twist$~\cite{GanterUsher}. This 2-group is an extension $1\to \pt\sq U(1) \to G^\twist\to G\to 1$ in the category of group objects internal to groupoids, where the extension is classified by~$\twist$.
\end{rmk}

\begin{ex} \label{ex:Conway}
Sporadic groups have interesting $U(1)$-valued cocycles that have emerged from recent studies of generalized moonshine. For example, let ${\rm Co}_0$ denote Conway's group of automorphisms of the Leech lattice~$\Lambda\subset \R^{24}$. In~\cite{TheoDavid} a degree~3 class of order~$24$ on ${\rm Co}_0$ is constructed, and in~\cite{TheoAnomaly} a related degree~$24$ class was constructed on the Monster group~$\M$. Hence, by applying double transgression from Example~\ref{2transgr} and using the description of $G^{(2)}\sq (G\times \SL_2(\Z))$ afforded by Example~\ref{ex:tori}, there are line bundles of order~24 on the groupoids of ${\rm Co}_0$- and $\M$-principal bundles on $2$-dimensional tori. 
\end{ex}
\begin{rmk}
The line bundles in the previous example are expected to encode anomalies of 2-dimensional conformal field theories, namely the monster conformal field theory for~$\M$, and the theory of 24 free fermions gauged by the action of~${\rm Co}_0$ on~$\R^{24}\simeq \Lambda\otimes \R$ (see Example~\ref{ex:ConwayEuler} below). 
\end{rmk}

\subsection{Integration on finite groupoids}\label{sec:integration}

Next we review the concept of groupoid cardinality, following \cite{FreedQuinn}, \cite{BaezDolan} and \cite{weinstein}. For a finite group $G$, let $d\mu$ denote the volume form on $\pt\sq G$ with total volume $1/|G|$. For $\mathcal{G}$ a finite groupoid there is an equivalence,
\beq
\coprod_{[g]\in \pi_0(\mathcal{G})} \pt\sq \Aut(g)\hookrightarrow \mathcal{G}\label{eq:finiteequiv}
\eeq
where the coproduct is indexed by choices of representatives of isomorphisms classes in $\mathcal{G}$ and $\Aut(g)$ is the automorphism group of an object $g\in \mathcal{G}$. 

\begin{defn}[\cite{FreedQuinn} \S2.1, \cite{BaezDolan} \S3]\label{defn:vol}
Define the volume form $d\mu$ on the groupoid $\mathcal{G}$ as the pullback of the volume form $\coprod d\mu_g$  along the inverse to the equivalence~\eqref{eq:finiteequiv}, where $d\mu_g$ is the previously defined volume form on $\pt\sq \Aut(g)$ with total volume $1/|\Aut(g)|$. 
\end{defn}

More explicitly, $d\mu$ is characterized by the values of integrals of functions $f\in C^\infty(\mathcal{G})$,
\beq
\int_{\mathcal{G}} fd\mu=\sum_{[g]\in \pi_0(\mathcal{G})} \frac{f(g)}{|\Aut(g)|}.\label{eq:grpdintegral}
\eeq

\begin{lem}\label{lem:measureinvariant}
The integral~\eqref{eq:grpdintegral} is independent of the choices of representatives $[g]\in \pi_0(\mathcal{G})$ and is invariant under equivalences of groupoids $\mathcal{G}\to \mathcal{G}'$. 
\end{lem}
\bp
By definition of $f\in C^\infty(\mathcal{G})$, $f(g)=f(g')$ for isomorphic objects $g$ and $g'$. Furthermore, $\Aut(g)\simeq \Aut(g')$ for isomorphic objects $g$ and $g'$, and the first claim follows. The second claim follows from the existence of a span of equivalences
$$
\mathcal{G} \hookleftarrow \coprod_{[g]\in \pi_0(\mathcal{G})} \pt\sq \Aut(g)\hookrightarrow \mathcal{G}',
$$
so that there are choices for which the formula for the integral~\eqref{eq:grpdintegral} is the same for integration on $\mathcal{G}$ and $\mathcal{G}'$. 
\ep

\begin{ex} Let $G$ be a finite group acting on a finite set $X$. Then an application of the orbit stabilizer theorem shows that $\vol(X\sq G):=\int_{X\sq G} 1d\mu=|X|/|G|$. More generally, Definition~\ref{defn:vol} is compatible with Euler characteristics of classifying spaces: $\chi(B\mathcal{G})= \int_\mathcal{G} 1 d\mu$, see~\cite[Example~2.7]{Leinster}. \end{ex}

Next we define integration along the fibers of a map $\zeta\colon \mathcal{H}\to \mathcal{G}$ of finite groupoids. We recall that the \emph{homotopy fiber} of an object $g\in \mathcal{G}$ is defined as the 2-pullback
\beq
\begin{tikzpicture}[baseline=(basepoint)];
\node (A) at (0,0) {${\rm Fib}(g)$};
\node (B) at (4,0) {$\mathcal{H}$};
\node (C) at (0,-1.5) {$*$};
\node (D) at (4,-1.5) {$\mathcal{G}$}; 
\draw[->] (A) to node [above] {$i$} (B);
\draw[->] (A) to (C);
\draw[->] (C) to node [above] {$\iota_g$} (D);
\draw[->] (B) to node [right] {$\zeta$} (D);
\path (0,-.75) coordinate (basepoint);
\end{tikzpicture}\label{eq:homotopyfiber}
\eeq
where $*$ is the (terminal) category with a single object and identity morphism, and $\iota_g$ picks out the object $g\in \mathcal{G}_0$ and identity isomorphism $\id_g\in \mathcal{G}_1$. Objects of ${\rm Fib}(g)$ therefore consist of pairs $(h,x)\in \mathcal{H}_0\times\mathcal{G}_1$ with the property that $x$ is an isomorphism from $g$ to $\zeta(h)$. Isomorphisms $(h,x)\to (h',x')$ in ${\rm Fib}(g)$ are the data of an isomorphism $h\stackrel{y}{\to} h'$ in $\mathcal{H}$ such that we have the equality $\zeta(y)\circ x=x'$ of isomorphisms $g\stackrel{\sim}{\to}\zeta(h')$ in $\mathcal{G}$.

\begin{defn} \label{defn:intfib}
For $f\in C^\infty(\mathcal{H})$, define the \emph{fiberwise integral} along $\zeta\colon \mathcal{H}\to \mathcal{G}$ by
\beq
(\int_{\mathcal{H}/\mathcal{G}} f d\mu)(g)=\int_{{\rm Fib}(g)} i^*f d\mu\label{eq:fibintegral}
\eeq
where $d\mu$ is the volume form on the finite groupoid ${\rm Fib}(g)$ from Definition~\ref{defn:vol}, and $i\colon {\rm Fib}(g)\to\mathcal{H}$ is the functor from~\eqref{eq:homotopyfiber}. If the homotopy fiber is empty, define the value of the integral to be zero. 
\end{defn} 

\begin{lem}
The values~\eqref{eq:fibintegral} determine a function $\int_{\mathcal{H}/\mathcal{G}} f d\mu\in C^\infty(\mathcal{G})$.
\end{lem}
\bp
The claim is equivalent to showing that the values~\eqref{eq:fibintegral} are the same for isomorphic objects $g$ and $g'$. Using Lemma~\ref{lem:measureinvariant}, this follows from the fact that an isomorphism between objects $g\stackrel{\sim}{\to} g'$ has a canonically associated equivalence of categories ${\rm Fib}(g)\stackrel{\sim}{\to} {\rm Fib}(g')$ over~$\mathcal{H}$. 
\ep

\begin{lem} 
The fiberwise integral can be computed as
$$
(\int_{\mathcal{H}/\mathcal{G}} f d\mu)(g)=\left(|\Aut(g)| \sum_{\zeta(h)=g} \frac{f(h)}{|\Aut(h)|}\right).
$$
where the sum is indexed by the subset of $[h]\in \pi_0(\mathcal{H})$ such that $\zeta(h)=g$. If this set is empty, the value of the function is zero. 
\end{lem}
\bp 
The equivalences~\eqref{eq:finiteequiv} give a 2-commuting square of groupoids
\beq
\begin{tikzpicture}[baseline=(basepoint)];
\node (A) at (0,0) {$\coprod \pt\sq \Aut(h)$};
\node (B) at (4,0) {$\mathcal{H}$};
\node (C) at (0,-1.5) {$\coprod \pt\sq \Aut(g)$};
\node (D) at (4,-1.5) {$\mathcal{G}$}; 
\draw[->] (A) to node [above] {$\sim$} (B);
\draw[->] (A) to node [left] {$\coprod \zeta_{h,g}$} (C);
\draw[->] (C) to node [above] {$\sim$} (D);
\draw[->] (B) to node [right] {$\zeta$} (D);
\path (0,-.75) coordinate (basepoint);
\end{tikzpicture}\label{eq:grpodsquare}
\eeq
where each $\zeta_{h,g}\colon \pt\sq \Aut(h)\to \pt\sq \Aut(g)$ is determined by a homomorphism $\Aut(h)\to \Aut(g)$. Since homotopy fibers are natural under equivalences of groupoids, we may compute ${\rm Fib}(g)$ in terms of the left vertical arrow in this diagram. Objects of ${\rm Fib}(g)$ in this description then consist of pairs $([h],x)\in \pi_0(\mathcal{H})\times \Aut(g)$ where $\zeta(h)=g$. Morphisms $([h],x)\to ([h],x')$ are determined by an element $y\in \Aut(h)$. This gives
$$
{\rm Fib}(g)\simeq \coprod_{\zeta(h)=g} \Aut(g)\sq \Aut(h),
$$
and the result follows. 
\ep

\begin{lem} 
The fiberwise integral of Definition~\ref{defn:intfib} satisfies the Fubini formula for compositions, $\mathcal{K}\to \mathcal{H}\to \mathcal{G}$,
\beq
\int_{\mathcal{K}/\mathcal{G}} f d\mu=\int_{\mathcal{H}/\mathcal{G}} \left(\int_{\mathcal{K}/\mathcal{H}} f d\mu\right)d\mu.\label{eq:Fubini}
\eeq
\end{lem}

\bp
By choosing equivalences~\eqref{eq:finiteequiv} for $\mathcal{K},\mathcal{H}$ and $\mathcal{G}$, this can be reduced to the case that $\mathcal{K}=\pt\sq K$, $\mathcal{H}=\pt\sq H$ and $\mathcal{G}=\pt\sq G$. Then~\eqref{eq:Fubini} becomes
$$
\frac{|K|}{|G|}f(\pt)=\frac{|K|}{|H|}\frac{|H|}{|G|}f(\pt)
$$
and the result follows. 
\ep

\subsection{Integration on inertia groupoids and character theory}\label{sec:HKR}

Next we explain how integration along maps of finite groupoids
$$
\Fun(\pt\sq\Z^d,\pt\sq H)\to \Fun(\pt\sq\Z^d,\pt\sq G)
$$ 
recovers the induction formulas from Hopkins--Kuhn--Ravenel character theory~\cite{HKR}. 
We start with the classical case for the character of an induced representation. Given a representation $\rho\colon H\to \End(V)$ of a finite group~$H$ and an inclusion of groups $H\hookrightarrow G$, we recall that the \emph{induced representation} ${\rm Ind}_H^G(\rho)$ is the $G$-representation whose underlying vector space consists of the $H$-invariant $V$-valued functions on $G$. The character of the induced representation is given by
\beq
\chi({\rm Ind}_H^G(\rho))(g)=\frac{1}{|H|}\sum_{g=xhx^{-1}} \chi(\rho)(h)\label{eq:inducedcharacter}
\eeq
where the sum is indexed by pairs $(h,x)\in H\times G$ such that $g=xhx^{-1}$. Since the character of an $H$-representation is a conjugation-invariant function on $H$, it defines an element of $C^\infty(H\sq H)$, and similarly characters of $G$-representations define elements of $C^\infty(G\sq G)$. A homomorphism $\zeta\colon H\to G$ then determines a map of finite groupoids,
\beq
H\sq H\simeq \Fun(\pt\sq \Z,\pt\sq H)\to \Fun(\pt\sq \Z,\pt\sq G)\simeq G\sq G, \label{eq:1charactermap}
\eeq
and integration along this gives a map $\int(-)d\mu\colon C^\infty(H\sq H)\to C^\infty(G\sq G)$ on characters. 

\begin{lem}\label{lem:Frob}
For a map $\zeta\colon H\to G$ of finite groups, the fiberwise integral along~\eqref{eq:1charactermap} recovers the formula~\eqref{eq:inducedcharacter}.
\end{lem}

\bp
By Definition~\ref{defn:intfib}, for $f\in C^\infty(H\sq H)$ and $g\in G$ we have
\beq
(\int f d\mu)(g)=\sum_{[g]=[\zeta(h)]} |C(g)|\frac{f(h)}{|C(h)|}\label{eq:Frobenious0}
\eeq
where the sum is indexed by conjugacy classes $[h]$ in $H$ such that $[\zeta(h)]=[g]$, and $C(h)<H$, $C(g)<G$ denote the centralizer of $h\in H$ and $g\in G$, respectively (i.e., the stabilizer for the conjugation actions). Applying the orbit stabilizer theorem for the conjugation action of $H$ on itself, $|H\cdot h|=|H|/|C(h)|$, we find
\beq
(\int f d\mu)(g)=\sum_{g=x\zeta(h)x^{-1}} \frac{|C(g)|}{|C(g)|}\frac{f(h)}{|C(h)||H\cdot h|}=\frac{1}{|H|}\sum_{g=x\zeta(h)x^{-1}} f(h)\label{eq:Frobenious}
\eeq
where the sum is indexed by $h\in H$ and $x\in G$ such that $g=x\zeta(h)x^{-1}$. Note that there exists an $x\in G$ with $g=x\zeta(h)x^{-1}$, then any other $x'$ with $g=x'\zeta(h)x'^{-1}$ can be written as $x'=yx$ for $y\in C(g)$, and hence the number of solutions $x$ to $g=x\zeta(h)x^{-1}$ is either zero or~$|C(g)|$. This count explains the factor of $1/|C(g)|$ in~\eqref{eq:Frobenious}, and the statement is proved. 
\ep

One can generalize integration along~\eqref{eq:1charactermap} to integration along the map
\beq
H^{(d)}\sq H=\Fun(\pt\sq \Z^d,\pt\sq H)\to \Fun(\pt\sq \Z^d,\pt\sq G)=G^{(d)}\sq G,\label{eq:inducedcharacterd}
\eeq
induced by a homomorphism $\zeta\colon H\to G$. 

\begin{prop}\label{prop:HKR}
Integration along~\eqref{eq:inducedcharacterd} reproduces the height~$d$ Hopkins--Kuhn--Ravenel induction formula~\cite[Theorem~D]{HKR}, 
\beq
(\int f d\mu)({\bf g})=\frac{1}{|H|}\sum_{{\bf g}=x\zeta({\bf h})x^{-1}} f({\bf h})\label{eq:HKR}
\eeq
where the sum is indexed by $({\bf h},x)\in H^{(d)}\times G$ such that ${\bf g}=x\zeta({\bf h})x^{-1}\in G^{(d)}$. Furthermore, if $f\in C^\infty(H^{(d)})^{H\times \SL_d(\Z)}$ is invariant under the $\SL_d(\Z)$-action on $H^{(d)}$, then $\int f d\mu\in C^\infty(G^{(d)})^{G\times \SL_d(\Z)}$ is also. 
\end{prop}

\bp
The first statement follows from the same orbit-stabilizer argument as in the equalities~\eqref{eq:Frobenious} in the proof of Lemma~\ref{lem:Frob}. The second statement follows by observing that $\gamma\in \SL_d(\Z)$ gives a bijection between solutions to ${\bf g}=x\zeta({\bf h})x^{-1}\in G^{(d)}$ and solutions to $\gamma\cdot {\bf g}=x\zeta(\gamma\cdot {\bf h})x^{-1}\in G^{(d)}$, and $f(\gamma\cdot{\bf h})=f({\bf h})$ by assumption. 
\ep

\begin{rmk}
A version of~\eqref{eq:HKR} for height 2 was stated in~\cite[\S5.2]{GanterHecke}. 
\end{rmk}

Lemma~\ref{lem:Frob} and Proposition~\ref{prop:HKR} admit a geometric rephrasing by way of Example~\ref{ex:tori}. Namely, a homomorphism $\zeta\colon H\to G$ allows one to form an associated $G$-bundle,
\beq
\zeta_*\colon \Bun_H(T^d)\to \Bun_G(T^d),\qquad P_{\bf h}\mapsto P_{\bf h}\times_\zeta G\simeq P_{\zeta({\bf h})}\label{eq:associated}
\eeq
where $H$ acts on $G$ via $\zeta$. The sums~\eqref{eq:Frobenious} and~\eqref{eq:HKR} can then be re-interpreted as being indexed by pairs given by (1) an $H$-bundle $P_{\bf h}$ and (2) an isomorphism of $G$-bundles $x\colon P_{\bf g}\stackrel{\sim}{\to} P_{\bf h}\times_\zeta G$. 
Similar sums over principal bundles appear when quantizing finite gauge theories, e.g., Dijkgraaf--Witten theory~\cite{DW} or Chern--Simons theory with finite gauge group~\cite{FreedQuinn}. We will return to this observation in~\S\ref{sec:discretetorsion}. 


\subsection{Internal groupoids}
In many applications, the objects and morphisms of a groupoid have additional structure. This leads to the notion of a groupoid internal to a category. An important instance of this later in the paper is the category of super Lie groupoids reviewed in~\S\ref{sec:superLie}. For now, we consider the following flavors of internal groupoids. 

\begin{defn} 
A \emph{Lie groupoid} is a groupoid whose objects and morphisms are manifolds, the maps~\eqref{eq:grpd} are smooth, and whose source map is a submersion. A \emph{complex analytic groupoid} is defined similarly, but where objects and morphisms are complex manifolds and the maps are holomorphic. 
\end{defn}

\begin{rmk} The submersion condition in the above definition guarantees that the fibered product~$\mathcal{G}_1\times_{\mathcal{G}_0}\mathcal{G}_1$ is a smooth manifold. \end{rmk}

We often regard finite groupoids as complex analytic or Lie groupoids, using that a finite set is a complex or smooth 0-manifold. Another class of examples comes from actions. 

\begin{ex}
For a Lie group $G$ acting on a manifold $M$, define the action Lie groupoid $M\sq G$ analogously to the action groupoid of Example~\ref{ex:action}. Similarly, for a complex analytic Lie group $G$, acting on a complex manifold $M$, let $M\sq G$ denote the action complex analytic groupoid. 
\end{ex}

\begin{defn}\label{defn:grpdfun} 
A \emph{(smooth) function} on a Lie groupoid $\mathcal{G}=\{\mathcal{G}_1\rightrightarrows \mathcal{G}_0\}$ is $f\in C^\infty(\mathcal{G}_0)$ such that $s^*f=t^*f$ on $\mathcal{G}_1$. If $\mathcal{G}$ is given a complex analytic structure, a \emph{holomorphic function} on $\mathcal{G}$ is $f\in \mathcal{O}(\mathcal{G}_0)$ such that $s^*f=t^*f$. 
\end{defn} 

\begin{ex} A function on $M\sq G$ is a $G$-invariant function on $M$. \end{ex}

\begin{defn}\label{defn:grpdvb}
A \emph{vector bundle} on a Lie groupoid is a vector bundle $V\to \mathcal{G}_0$ and an isomorphism of vector bundles $s^*V\simeq t^*V$ on $\mathcal{G}_1$ satisfying a cocycle condition on $\mathcal{G}_1\times_{\mathcal{G}_0}\mathcal{G}_1$. If $\mathcal{G}$ is given a complex analytic structure, a \emph{holomorphic vector bundle} on $\mathcal{G}$ is a holomorphic vector bundle $V\to \mathcal{G}_0$ and an isomorphism $s^*V\simeq t^*V$ of holomorphic vector bundles on~$\mathcal{G}_1$ satisfying a cocycle condition on $\mathcal{G}_1\times_{\mathcal{G}_0}\mathcal{G}_1$. 
\end{defn}

\begin{ex} A vector bundle on $M\sq G$ is a $G$-equivariant vector bundle on $M$. \end{ex}

\begin{defn} \label{defn:grpdsheaf} A \emph{sheaf} on a complex analytic or Lie groupoid $\mathcal{G}=\{\mathcal{G}_1\rightrightarrows\mathcal{G}_0\}$ is a sheaf $\F$ on $\mathcal{G}_0$ and an isomorphism of sheaves $s^*\F\simeq t^*\F$ on $\mathcal{G}_1$ satisfying a cocycle condition on $\mathcal{G}_1\times_{\mathcal{G}_0}\mathcal{G}_1$. A \emph{morphism of sheaves} $\F\to \F'$ on $\mathcal{G}$ is a morphism of sheaves on $\mathcal{G}_0$ compatible with the isomorphisms $s^*\F\simeq t^*\F$ and $s^*\F'\simeq t^*\F'$.\end{defn} 

Explicitly, a sheaf on $\mathcal{G}$ assigns to each subset $U\subset \mathcal{G}_0$ a set $\F(U)$ and for each $x\in \mathcal{G}_1$ determining an isomorphism $x\colon U\to U'$ it assigns bijection $\F(U')\to \F(U)$. As per the usual terminology, a sheaf of abelian groups, rings, modules, etc. will endow the values $\F(U)$ with a structure, require the restriction maps to respect the structure, and require the maps $\F(U')\to \F(U)$ associated with isomorphisms in $\mathcal{G}$ to be isomorphisms for the structure. A \emph{complex of sheaves} on $\mathcal{G}$ is a collection of sheaves $\{\F^k\}_{k\in \Z}$ of abelian groups on $\mathcal{G}$ with morphisms $d\colon \F^k\to \mathcal{F}^{k+1}$ such that $d\circ d=0$. A sheaf of \emph{commutative differential graded algebras} (cdgas) will be taken to mean a complex of sheaves with an associative, graded commutative multiplication $\F^k\times \F^j\to \F^{k+j}$ compatible with the differential. 

\section{Twisted equivariant K-theory and elliptic cohomology over $\C$}\label{sec:cohomo}

This section reviews the definition of complexified twisted equivariant K-theory and complex analytic twisted equivariant elliptic cohomology for a finite group~$G$ acting on a compact manifold~$M$. In both cases, we view these theories as sheaves of commutative differential graded algebras on certain groupoids. For K-theory, the relevant groupoid is $G\sq G$, the inertia groupoid of $\pt\sq G$. In the elliptic case, the groupoid is the complex analytic groupoid $\Bun_G(\EE)$ of $G$-bundles over complex elliptic curves. This is an enhancement of the double inertia groupoid of~$\pt\sq G$. Twistings of these cohomology theories arise from line bundles on the respective groupoids. We construct examples of such line bundles by transgression of classes in~$\H^d(BG;U(1))$ for $d=2$ in the case of K-theory and $d=3$ for elliptic cohomology. We also explain the relationship between $\Z/n$-equivariant K-theory (respectively, elliptic cohomology) with $n$-torsion points on the multiplicative group (respectively, the universal dual elliptic curve), and provide examples of $(\Z/n)^{\times k}$-equivariant Euler classes in this description. Most of the material in this section seems to be previously-known to experts, but we do not know of a reference. Throughout, differential forms will be taken with complex coefficients, $\Omega^\bullet(M)=\Omega^\bullet(M;\C)$.

\subsection{Complexified equivariant K-theory as global sections of a sheaf}

For $M$ a compact manifold with the action of a finite group $G$ and $\twist\colon G\times G\to U(1)$ a 2-cocycle, a description of the $\twist$-twisted complexified $G$-equivariant K-theory group ${\rm K}_G^\twist(M)\otimes \C$ was given by Adem and Ruan (also see~\cite[Proposition~3.11]{FHTcomplex}) generalizing the untwisted calculation by Atiyah and Segal \cite[Theorem~2]{AtiyahSegalEuler}.

\begin{thm}[\cite{AdemRuan} Theorem~7.4] Let $\twist\colon G\times G\to U(1)$ be a 
2-cocycle. For a manifold~$M$ with $G$-action, the $\twist$-twisted $G$-equivariant ${\rm K}$-theory of~$M$ with complex coefficients is isomorphic to the product indexed by conjugacy classes of $g\in G$
\beq
&&{\rm K}_G^{\bullet+\twist}(M)\otimes \C \cong\prod_{[g]} ({\rm H}_{\rm dR}^{\bullet}(M^g;\C[\beta,\beta^{-1}])\otimes \chi^\twist_g)^{C(g)}\qquad |\beta|=-2\label{eq:AdemRuan}
\eeq
where $C(g)<G$ denotes the centralizer of $g$, and $\chi^\twist_g\colon C(g)\to U(1)$ is the 1-dimensional representation given by $h\mapsto \twist(h,g)\twist(g,h)^{-1}$. \label{thm:FHTcomplex}
\end{thm}

Next we repackage ${\rm K}_G^{\bullet+\twist}(M)\otimes \C$ as the derived global sections of a complex of sheaves on~$G\sq G$. When $\alpha$ is trivial, this is essentially a simplification of the definitions from~\cite[\S1]{BlockGetzler} and~\cite[Definition~22]{Vergne} given for an arbitrary compact Lie group~$G$.

\begin{defn}\label{defn:Ksheaf}
Define the sheaf $\mathcal{K}^\bullet_G(M)$ of cdgas on the set $G$ by the assignment 
\beq
&&U\mapsto \Big(\prod_{g\in U} \Omega^\bullet(M^g;\C[\beta,\beta^{-1}]),d\Big) ,\qquad U\subset G.\label{eq:Ksheaf}
\eeq
For maps of sets $h\colon U\to U'$ determined by the conjugation action by $h\in G$, take the isomorphisms of sheaves gotten by pulling back differential forms along the diffeomorphisms $h\colon M^g\to M^{ghg^{-1}}$. Finally, for a line bundle $\mathcal{T}$ on $G\sq G$, define the sheaf 
$$
\mathcal{K}_G^{\bullet+\mathcal{T}}(M):= \mathcal{K}_G^\bullet(M)\otimes \mathcal{T}
$$ 
as the tensor product of sheaves on $G\sq G$. We use the notation $\mathcal{K}_G^{\bullet+\alpha}(M)$ when $\mathcal{T}=\mathcal{T}^\alpha$ for $\alpha\colon G\times G\to U(1)$ a 2-cocycle as in Example~\ref{1transgr}. 
\end{defn}

\begin{prop}\label{prop:dgs2}
The derived global sections of $\mathcal{K}^{\bullet+\twist}_G(M)$ compute~${\rm K}_G^{\bullet+\twist}(M)\otimes \C$. 
\end{prop}

\bp
Derived global sections can be identified with the groupoid cohomology for $G\sq G$ with values in the sheaf~\eqref{eq:Ksheaf}; see~\cite[{\S}3.3]{BehrendXu}. Groupoid cohomology is invariant under equivalences, so we restrict along the fully faithful and essentially surjective functor
\beq
\coprod_{[g]} \pt\sq C(g)\hookrightarrow G\sq G.\label{eq:coproductofshit}
\eeq
On each component $\pt\sq C(g)$, the groupoid cohomology is the group hypercohomology with values in the complex of $C(g)$-modules $(\Omega^\bullet(M^g;\C[\beta,\beta^{-1}])\otimes\mathcal{T}^\twist_g,d)$. 
This is concentrated in degree zero: the order of $C(g)$ is invertible in the module. Hence, we identify the groupoid cohomology on a component of~\eqref{eq:coproductofshit} with the $C(g)$-invariant subcomplex $((\Omega^\bullet(M^g;\C[\beta,\beta^{-1}])\otimes \mathcal{T}^\twist_g)^{C(g)},d)$. The action of $C(g)$ on $\mathcal{T}^\twist_g$ agrees with the one-dimensional representation $\chi_g^\alpha\colon \C(g)\to \C^\times$ in the statement of Theorem~\ref{thm:FHTcomplex}. So in total we find that derived global sections are 
$$
\R\Gamma(G\sq G,\mathcal{K}_G^{\bullet+\twist}(M))\simeq \prod_{[g]}\H( (\Omega^\bullet(M^g;\C[\beta,\beta^{-1}])\otimes \chi_g^\alpha)^{C(g)},d),
$$
which agrees with the right hand side of~\eqref{eq:AdemRuan}.
\ep

We provide a few examples of structures in complexified equivariant K-theory inherited from the description from Definition~\ref{eq:Ksheaf}. These can be contrasted with their elliptic counterparts, Examples~\ref{ex:ellAS},~\ref{ex:Ellinduction}~\ref{ex:torsionpts}, and~\ref{ex:Klein}. 

\begin{ex}[Atiyah--Segal completion] Borel equivariant complexified K-theory is $G$-equivariant de~Rham cohomology with coefficients in~$\C[\beta,\beta^{-1}]$, i.e., it assigns to a $G$-manifold $M$ the ring $\H^\bullet_G(M;\C[\beta,\beta^{-1}])\simeq \H^\bullet(M;\C[\beta,\beta^{-1}])^G$. The sections of $\mathcal{K}_G^\bullet(M)$ over the subgroupoid $\pt\sq G\simeq e\sq C(e)\hookrightarrow G\sq G$ associated to the identity element $e\in G$ is therefore a cocycle model for Borel equivariant complexified K-theory. By Proposition~\ref{prop:dgs2}, the restriction of a derived global section to this fiber determines a map
\beq
\K_G^\bullet(M)\otimes \C\to \H^\bullet_G(M;\C[\beta,\beta^{-1}])\label{eq:AScompletion}
\eeq
which we recognize as the complexified Atiyah--Segal completion map~\cite{Atiyah_Segal_complete}. Using  Example~\ref{ex:tori}, the map~\eqref{eq:AScompletion} can be described in terms of the groupoid of $G$-bundles on $S^1$, namely it restricts the sheaf $\mathcal{K}_G^\bullet(M)$ to sections over the trivial $G$-bundle. We emphasize that this restriction forgets information. For example, the image of a character of a representation $\chi(\rho)\in C^\infty(G\sq G)\simeq \K_G^0(\pt)\otimes \C$ under~\eqref{eq:AScompletion} reads off the dimension of the representation. 
\end{ex}
\begin{ex}[Induction] From Lemma~\ref{lem:Frob}, applying groupoid integration for a map $\zeta\colon H\sq H\to G\sq G$ induced by an inclusion of groups determines a map of sheaves
$$
\Ind_H^G\colon \mathcal{K}_H(\pt)\to \mathcal{K}_G(\pt)
$$
that on global sections recovers the formula for an induced representation. We return to this structure in greater detail in~\S\ref{sec:discretetorsion}. 
\end{ex}

\begin{ex}[Torsion points on the multiplicative group] \label{eq:Kex1}
Identify the function $1-e^{2\pi i \theta}\in C^\infty(S^1)$ with the character of the virtual representation ${\bf 1}\ominus {\bf R}$ where ${\bf 1}$ is the trivial representation of $S^1$ and ${\bf R}$ is the 1-dimensional complex representation gotten from the identification $S^1 \simeq U(1)$. Consider the restriction of $1-e^{2\pi i \theta}$ along the inclusion $\mu_n\colon \Z/n\hookrightarrow S^1\subset \C^\times$ as the $n$th roots of unity. This determines a global section $c(\mu_n):=1-e^{2\pi i \theta}\in C^\infty(\Z/n\sq \Z/n)\simeq \Gamma(\Z/n\sq \Z/n;\mathcal{K}^0_{\Z/n}(\pt))$. Multiplication and projection provide three homomorphisms $\Z/n\times \Z/n\to \Z/n$ that determine three functors between finite groupoids
$$
m,p_1,p_2\colon (\Z/n\sq \Z/n)\times (\Z/n\sq \Z/n)\simeq (\Z/n\times \Z/n)\sq (\Z/n\times \Z/n)\to \Z/n\sq \Z/n.
$$
We may consider the pullback sheaves $\mathcal{K}^0_{\Z/n}(\pt)$ along these maps. This yields the equality of sections of $\mathcal{K}^0_{\Z/n\times \Z/n}(\pt)$
\beq
&&m^*(c(\mu_n))=c_1+c_2-c_1c_2\qquad c_i=p_i^*(c(\mu_n)),\label{eq:multFGL}
\eeq
which is the formula for the multiplicative group law. Hence, we may view $\Z/n\sq \Z/n$ with its sheaf $\mathcal{K}^0_{\Z/n}(\pt)$ as the $n$-torsion points on the multiplicative group, and the class $c(\mu_n)$ is the restriction of the standard coordinate to the torsion points. 
\end{ex}

\begin{rmk} The virtual representation $[{\bf 1}\ominus {\bf R}]\in \Rep(U(1))\simeq \K_{U(1)}(\pt)$ is a $U(1)$-equivariant refinement (relative to the Atiyah--Segal completion map) of the Chern class class $c(\mathcal{O}(-1))=[1-\mathcal{O}(-1)]\in \K^0(\CP^\infty)$ of the tautological line bundle associated with the usual complex orientation of K-theory. From Quillen's work relating complex orientations and formal group laws~\cite{Quillen}, this both explains the notation $c(\mu_n)$ above and the appearance of the multiplicative group law~\eqref{eq:multFGL}. \end{rmk}

\begin{ex}[Euler classes of Spin representations]\label{eq:Kex2}
Following Remark~\ref{rmk:projchar}, we can obtain twisted classes from characters of projective representations. So consider the projective representation of $S^1\simeq \SO(2)$ given by $\Delta^+\ominus \Delta^-$ where $\Delta^\pm$ are the irreducible representations of $\Spin(2)$. The character of this representation is a section of a line bundle over $S^1=\R/\Z$ given by the function 
\beq
e^{\pi i \theta}-e^{-\pi i \theta}=2i \sin(\pi \theta),\in C^\infty(\R)\qquad 2i\sin(\pi (\theta+1))=-2i\sin(\pi \theta),\label{eq:mobius}
\eeq
where the failure to be invariant under $\theta\mapsto \theta+1$ shows that this is a section of the (M\"obius) line bundle on $S^1=\R/\Z$. Restricting the line bundle along the inclusion $\Z/n\hookrightarrow S^1\subset \C^\times$ as the $n$th roots of unity defines twisted Euler classes, 
$$
\Eu(\mu_n):=\beta^{-1} 2i \sin(\pi \theta)\in \mathcal{K}_{\Z/n}^{2+\alpha}(\pt),\qquad |\beta|=-2
$$
where $\beta^{-1}$ puts the Euler class in the correct degree, and the twisting reflects the sign ambiguity and obstructs the existence of a \emph{spin structure}. Indeed, a spin structure on a representation is a lift $B\tilde{\rho}$ (depending on if we are viewing the representation as real or complex)
\beq
\begin{tikzpicture}[baseline=(basepoint)];
\node (A) at (1,0) {$BG$};
\node (C) at (5,0) {$BO(k)$};
\node (D) at (5,1.2) {$B\Spin(k)$};
\draw[->] (A) to node [below] {$B\rho$} (C);
\draw[->] (D) to (C);
\draw[->,dashed] (A) to node [above] {$B\tilde\rho$} (D);
\path (0,.5) coordinate (basepoint);
\end{tikzpicture}\qquad 
\begin{tikzpicture}[baseline=(basepoint)];
\node (A) at (1,0) {$BG$};
\node (B) at (3,0) {$B\U(k)$};
\node (C) at (5,0) {$B\SO(2k).$};
\node (D) at (5,1.2) {$B\Spin(2k)$};
\draw[->] (A) to node [below] {$B\rho$} (B);
\draw[->] (B) to (C);
\draw[->] (D) to (C);
\draw[->,dashed] (A) to node [above] {$B\tilde\rho$} (D);
\path (0,.5) coordinate (basepoint);
\end{tikzpicture}\label{diagam:spinlift}
\eeq
Assuming $\rho$ is orientation-preserving, the existence of a spin structure is obstructed by $[w_2(\rho)]\in \H^2(BG;\Z/2)$, the pullback of the universal Stiefel--Whitney class $w_2\in \H^2(B\SO(k);\Z/2)$. If $[w_2(\rho)]=0$, then the homotopy classes of lifts $B\tilde\rho$ are parameterized by~$H^1(BG;\Z/2)$. 
\end{ex}
\begin{ex}[Untwisted $\Z/2$-equivariant Euler classes]\label{ex:Kthyz2}
We observe that the square of the Euler class $\Eu(\mu_2)^2=\Eu(\mu_2\oplus \mu_2)\in \mathcal{K}_{\Z/2}(\pt)$ determines a canonically untwisted class, and hence determines a choice of spin structure on~$\mu_2^{\oplus 2}\simeq \C^2\sq \Z/2$. The existence of a spin structure can be seen by identifying $\mu_2^{\oplus 2}\simeq \sigma^{\oplus 4}$ where $\sigma\simeq \R\sq \Z/2$ is the sign representation of $\Z/2$ on $\R$. The representation $\sigma$ has total Stiefel--Whitney class $1+w_1\in \H(B\Z/2;\Z/2)\simeq \H(\R \mathbb{P}^\infty;\Z/2)\simeq (\Z/2)[w_1]$, where $w_1$ is the polynomial generator. Using multiplicativity of the total  Stiefel--Whitney class, the smallest $k\in \N$ for which $w_2(\sigma^{\oplus k})$ vanishes is $k=4$. The particular choice of spin structure associated with the untwisted class $\Eu(\mu_2)^2\in \mathcal{K}_{\Z/2}(\pt)$ comes from the formula~\eqref{eq:mobius} that specifies a choice of spin lift of a generator of~$\Z/2$. The other choice of spin structure corresponds to the other choice of lift.
\end{ex}

\subsection{$G$-bundles on complex analytic elliptic curves}
Let $\HH\subset \C$ be the upper half plane. The groupoid of elliptic curves is the complex analytic groupoid~$\Mell\simeq \HH\sq \SL_2(\Z)$ for the $\SL_2(\Z)$-action on $\HH$ by fractional linear transformations,
\beq
&&\gamma\cdot \tau= \frac{a\tau+b}{c\tau+d},\qquad \tau\in \HH, \ \gamma=\left[\begin{array}{cc} a & b \\ c & d\end{array}\right]\in \SL_2(\Z). \label{eq:fraction}
\eeq
A point $\tau\in \HH$ determines a complex analytic elliptic curve $E_\tau=\C/(\tau\Z\oplus\Z)$. An isomorphism of complex analytic elliptic curves~$E_\tau\simeq E_{\tau'}$ is specified by $\gamma\in \SL_2(\Z)$ with $\tau'=\gamma \cdot \tau$. As before, let $G^{(2)}\subset G^{\times 2}$
denote the set of pairs of commuting elements in $G$. Define the complex analytic groupoid
\beq
&&\Bun_G(\EE):= (\HH\times G^{(2)})\sq (G\times \SL_2(\Z))\label{eq:BunGE}
\eeq
for the previous $\SL_2(\Z)$-action on $\HH$ and the $G\times \SL_2(\Z)$-action on $G^{(2)}$ given by~\eqref{eq:G2act}. 
 A point $(\tau,g_1,g_2)\in \HH\times G^{(2)}$ determines a $G$-bundle on $E_\tau$ via the quotient
$$
P_{\tau,g_1,g_2}:=(G\times \C)/((g_1,\tau)\Z\oplus(g_2,1)\Z)\to \C/(\tau\Z\oplus \Z)=E_\tau,
$$
where $(g_1,g_2)$ generate a $\Z^2$-action on $G$ by left multiplication. An isomorphism of $G$-bundles $P_{\tau,g_1,g_2}\simeq P_{\tau',g_1',g_2'}$ covering an isomorphism $E_\tau\simeq E_{\tau'}$ of elliptic curves is specified by $(h,\gamma)\in G\times \SL_2(\Z)$ with $(g_1',g_2')=(h,\gamma)\cdot (g_1,g_2)$ and $\tau'=\gamma \cdot \tau$. A homomorphism $G\to H$ determines a map of complex analytic groupoids
\beq
\Bun_G(\EE)\to \Bun_H(\EE),\label{eq:ellnaturality}
\eeq
and we have isomorphisms of complex analytic groupoids
\beq
\Bun_{G\times H}(\EE)\simeq \Bun_G(\EE)\times_{\Mell}\Bun_H(\EE). \label{eq:ellproduct}
\eeq

\begin{rmk}
The map $\Bun_G(\EE)\to \Mell$ is an \emph{isofibration}, which implies that the strict fibered product and homotopy fibered product of groupoids in~\eqref{eq:ellproduct} coincide; e.g., see~\cite[Lemma~8.2]{CoyneNoohi}. 
\end{rmk}

Since $\HH$ is contractible and $G^{(2)}$ is discrete, any line bundle on $\Bun_G(\EE)$ can be specified in terms of a $G\times \SL_2(\Z)$-equivariant structure on the trivial line bundle on $\HH\times G^{(2)}$. We are interested in two classes of such line bundles below. The first is generated by the \emph{Hodge bundle} $\omega$ that pulls back along the projection $\Bun_G(\EE)\to \Mell$. Sections of $\omega^{\otimes k}$ over $\Mell$ are defined by the transformation property
$$
f\left(\frac{a\tau+b}{c\tau+d}\right)=(c\tau+d)^kf(\tau),\qquad f\in \mathcal{O}(\HH). 
$$
The sheaf of sections $\bigoplus_k \omega^{\otimes k}$ on $\Mell$ has a multiplication, using the evident maps $\omega^{\otimes k}\otimes \omega^{\otimes l}\stackrel{\sim}{\to} \omega^{\otimes (k+l)}$. Modular forms of weight $k$ are global sections of $\omega^{\otimes k}$; see Remark~\ref{rmk:holomorphicatinfty}. 

The second type of line bundle on $\Bun_G(\EE)$ pulls back the line bundles $\mathcal{T}^\twist$ from Example~\ref{2transgr} along the projection
\beq
\Bun_G(\EE)=(\HH\times G^{(2)})\sq (G\times \SL_2(\Z))\stackrel{pr}{\to} G^{(2)}\sq (G\times \SL_2(\Z)). \label{eq:topological}
\eeq
In the $G$-bundle description of the target afforded by Example~\ref{ex:tori},~\eqref{eq:topological} forgets the complex structure on the elliptic curve and remembers the $G$-bundle on its underlying smooth torus. We use the same notation $\mathcal{T}^\twist$ to denote the resulting line bundle on $\Bun_G(\EE)$ depending on a 3-cocycle~$\twist$ on~$G$.


\begin{rmk}\label{rmk:holomorphicatinfty}
When considering global sections of sheaves over $\Bun_G(\EE)$ and $\Mell$, one can ask for prescribed behavior of sections as $\tau\to i\infty$: often one imposes holomorphy or meromorphicity at the cusp. We choose not to impose this restriction on the sheaf $\Ell_G^\bullet(M)$ below. As a result, for example, sections of $\Ell_{\{e\}}^\bullet(\pt)$ can be identified with sections of~$\omega^{\otimes -k/2}$. These are \emph{weak} modular forms. \emph{Weakly holomorphic} modular forms are a subring of the ring of weak modular forms that are meromorphic at the cusp. 
\end{rmk}

\subsection{Complex analytic equivariant elliptic cohomology}

The following is a finite group simplification of the definitions in~\cite[\S3]{BET0} and \cite[\S6]{MicheleEll}. 

\begin{defn}\label{defn:Ell}
For a $G$-manifold $M$, \emph{complex analytic equivariant elliptic cohomology} of~$M$ is a sheaf of commutative differential graded algebras on~$\HH\times G^{(2)}$ 
that to an open subset $U \subset \HH\times G^{(2)}$ assigns the complex of $\mathcal{O}(U)$-modules, 
$$
\Ell^\bullet_G(M)(U) := \Big(\prod_{(g_1,g_2)\in p(U)}\mathcal{O}(U_{g_1,g_2};\Omega^\bullet(M^{g_1,g_2})[\beta,\beta^{-1}]),d\Big), \qquad |\beta|=-2
$$
where $d$ is the de~Rham differential, $p\colon \HH\times G^{(2)}\to G^{(2)}$ is the projection, and $U_{g_1,g_2}$ is the fiber of $U$ over $(g_1,g_2)\in G^{(2)}$. We promote this to a $G\times \SL_2(\Z)$-equivariant sheaf that to a map $f_{\gamma,h}\colon U\to U'$ assigns the isomorphism of complexes of sheaves
$$
f_{\gamma,h}^*\colon \Ell^\bullet_G(U')\to \Ell^\bullet_G(U),\quad \gamma\in \SL_2(\Z), \ h\in G
$$
by first pulling back functions using the isomorphism from left multiplication by $h\in G$,
$$
\Omega^\bullet(M^{g_1',g_2'})[\beta,\beta^{-1}]\stackrel{h^*}{\to} \Omega^\bullet(M^{g_1,g_2})[\beta,\beta^{-1}] \qquad (g_1',g_2')=(h,\gamma)\cdot (g_1,g_2)
$$
and post-composing with the algebra automorphism determined by $\beta\mapsto \beta/(c\tau+d)$. The action on the indexing set $(g_1,g_2)\in G^{(2)}$ is given by~\eqref{eq:G2act}. This determines a sheaf on the stack $[\HH\times G^{(2)}\sq (G\times \SL_2(\Z))]$. For a line bundle $\mathcal{T}$ on $\Bun_G(\EE)$, the \emph{$\mathcal{T}$-twisted complex analytic equivariant elliptic cohomology} of $M$ is the tensor product of sheaves 
$$
\Ell^{\bullet+\mathcal{T}}_G(M):=\Ell^\bullet_G(M)\otimes\mathcal{T}.
$$ 
We use the notation $\Ell^{\bullet+\alpha}_G(M)$ when $\mathcal{T}=\mathcal{T}^\alpha$ from Example~\ref{2transgr}. 
\end{defn} 

\begin{lem} There is an isomorphism of complexes of sheaves
\beq
\Ell_G^\bullet(M)\otimes \omega\stackrel{\sim}{\to} \Ell_G^{\bullet-2}(M).\label{eq:Bott}
\eeq
\end{lem}
\bp
We observe that $\beta$ determines a nonvanishing section of $\omega$ over $\HH\times G^{(2)}$. Hence the isomorphism~\eqref{eq:Bott} is given by
$$
f\otimes\sigma\mapsto f\beta\sigma, \qquad f\in \Gamma(U,\Ell_G^\bullet(M)), \ \sigma\in \Gamma(U;\omega). 
$$
\ep

\begin{rmk}
The property~\eqref{eq:Bott} is a complex analytic version of what is called \emph{weakly periodic} in~\cite[Remark~1.4]{Lurie}. It is a form of Bott periodicity that exists locally on~$\Bun_G(\EE)$. 
\end{rmk}

Let $\Gamma_{g_1,g_2}<G\times \SL_2(\Z)$ be the subgroup stabilizing $\langle g_1,g_2\rangle\colon \Z^2\to G$, and $G^{[2]}=G^{(2)}/(G\times \SL_2(\Z))$ denote the orbits. 

\begin{prop}\label{prop:dgs}
Viewing $\Bun_G(\EE)$ as a complex analytic stack, the derived global sections of $\Ell^{\bullet+\twist}_G(M)$ compute the graded vector space
$$
\prod_{[g_1,g_2]\in G^{[2]}} \Big(\H^\bullet\big(M^{g_1,g_2};\mathcal{O}(\HH)[\beta,\beta^{-1}])\otimes\chi^\twist_{g_1,g_2} \Big)^{\Gamma_{g_1,g_2}}
$$
where $\Gamma_{g_1,g_2}$ acts by pulling back differential forms along the diffeomorphisms $M^{g_1,g_2}\simeq M^{g_1',g_2'}=M^{(g_1,g_2)}$ for $(g_1',g_2')=(h,\gamma)\cdot (g_1,g_2)$ composed with the action on coefficients given by~\eqref{eq:fraction} and $\beta\mapsto \beta/(c\tau+d)$, and $\chi^\twist_{g_1,g_2}$ is the 1-dimensional representation of determined by the restriction of the 1-cocycle $\mathcal{T}^\twist$ to $\Gamma_{g_1,g_2}<G\times \SL_2(\Z)$. 
\end{prop}
\bp
The proof is similar to the proof of Proposition~\ref{prop:dgs2}, with the modification that $\Bun_G(\EE)$ is a complex analytic groupoid and we use the Dolbeault resolution to compute derived sections. First we restrict along the equivalence
\beq
\coprod_{[g_1,g_2]} \HH\sq \Gamma_{g_1,g_2}\hookrightarrow (\HH\times G^{(2)})\sq G\times \SL_2(\Z).\label{eq:anothercoprodofshit}
\eeq
On a component of~\eqref{eq:anothercoprodofshit}, the groupoid cohomology is the group hypercohomology of $\Gamma_{g_1,g_2}$ with values in the double complex $\Omega^\bullet(M^{g_1,g_2};\Omega^{0,*}(\HH)[\beta,\beta^{-1}])$ whose differentials are the de~Rham differential on $M^{g_1,g_2}$ and the Dolbeault differential on $\HH$. 
As before, the group cohomology computation collapses because the stabilizers for the action of $\Gamma_{g_1,g_2}$ on $\HH$ are finite groups whose order is invertible in the module. The Dolbeault complex is also acyclic because $\HH$ is Stein. 
%
Finally, by definition the action of $\Gamma_{g_1,g_2}$ on $\mathcal{T}^\twist_{g_1,g_2}$ is through the representation $\chi_{g_1,g_2}^\alpha$ in the statement of the proposition. So in total we find that derived global sections are computed by the product of de~Rham complexes
$$
\R\Gamma(\Bun_G(\EE),\Ell^{\bullet+\twist}_G(M))\simeq \H(\prod_{[g_1,g_2]} (\Omega^\bullet(M^{g_1,g_2};\mathcal{O}(\HH)[\beta,\beta^{-1}]\otimes \chi_{g_1,g_2}^\alpha)^{\Gamma_{g_1,g_2}}),d),
$$
proving the proposition. 
\ep

\begin{notation} Denote the derived global sections by
\beq
&&\TMF^{\bullet+\twist}_G(M)\otimes \C:=\R\Gamma(\Bun_G(\EE),\Ell^{\bullet+\twist}_G(M)).\label{eq:notation1}
\eeq
\end{notation}

\begin{rmk}
We justify the notation~\eqref{eq:notation1}. There is a sheaf of elliptic spectra on the moduli stack of elliptic curves (defined over a general commutative ring $R$) whose derived global sections are the spectrum of topological modular forms, denoted~TMF; e.g., see~\cite{HopkinsICM2002,Lurie,Goerss}. Aspects of the equivariant picture have also been worked out~\cite{LurieIII,LenartDavid}, where roughly speaking the relevant sheaf of spectra is on a moduli space of $G$-bundles over (derived, oriented) elliptic curves. Complexification is expected to recover the sheaf of cdgas from Definition~\ref{defn:Ell}, and hence we use the notation~\eqref{eq:notation1} for the derived global sections of this sheaf. 
\end{rmk}

\begin{rmk}
Devoto constructs a version of equivariant elliptic cohomology for the level~2 subgroup $\Gamma_0(2)<\SL_2(\Z)$~\cite{DevotoI,DevotoII}. The derived global sections computed in Proposition~\ref{prop:dgs} are an $\SL_2(\Z)$-version of the complexification of Devoto's theory, with the caveat that we do not impose meromorphicity at infinity (see Remark~\ref{rmk:holomorphicatinfty}). 
\end{rmk}

\begin{ex}\label{ex:ellAS} Borel equivariant complex analytic elliptic cohomology is $G$-equivariant de~Rham cohomology with coefficients in~$(\mathcal{O}(\HH)[\beta,\beta^{-1}])^{\SL_2(\Z)}$, the ring of weak modular forms. Hence, evaluating this Borel equivariant theory on a $G$-manifold~$M$ yields
$$
\H^\bullet_G(M;(\mathcal{O}(\HH)[\beta,\beta^{-1}])^{\SL_2(\Z)})\simeq \H^\bullet(M;\mathcal{O}(\HH)[\beta,\beta^{-1}])^{G\times \SL_2(\Z)}.
$$ 
The sections of $\Ell_G^\bullet(M)$ over the subgroupoid $\HH \sq (G\times \SL_2(\Z))\hookrightarrow \Bun_G(\EE)$ associated with the trivial $G$-bundle is therefore a cocycle model for Borel equivariant complex analytic elliptic cohomology. By Proposition~\ref{prop:dgs}, the restriction of a derived global section to this subgroupoid of trivial $G$-bundles determines a map
\beq
\Ell_G^\bullet(M)\to \H^\bullet(M;\mathcal{O}(\HH)[\beta,\beta^{-1}])^{G\times \SL_2(\Z)}.\label{eq:ellAScompletion}
\eeq
In view of~\eqref{eq:AScompletion}, we view this as an \emph{elliptic} Atiyah--Segal completion map. 
\end{ex}

\begin{ex}\label{ex:Ellinduction} The formula~\eqref{eq:HKR} determines a map of sheaves $\Ell_H^\bullet(\pt)\to \Ell_G^\bullet(\pt)$ that we identify with the height~2 Hopkins--Kuhn--Ravenel induction. We return to this structure in greater detail in~\S\ref{sec:discretetorsion}. 
\end{ex}

\subsection{Twisted equivariant elliptic Euler classes}\label{sec:Euler}

The equivariant K-theory classes in Examples~\ref{eq:Kex1} and~\ref{eq:Kex2} came from characters of (projective) representations. In the elliptic case, an analogous source of classes comes from \emph{theta functions}, which are related to characters of positive energy representations of loop groups. In terms of the geometry of elliptic curves, theta functions are sections of line bundles over the universal dual elliptic curve. When these line bundles are restricted to $n$-torsion points of the curve they determine (twisted) $\Z/n$-equivariant classes, which are classically known as \emph{Klein forms}. We start by explaining the relationship between $\Z/n$-equivariant elliptic cohomology and $n$-torsion points on elliptic curves. 


\begin{ex}[$\Ell_{\Z/n}^0$ as functions on $n$-torsion points]\label{ex:torsionpts}
Consider the universal dual elliptic curve $\EE^\vee=\widetilde{\EE}^\vee\sq \SL_2(\Z)$
$$
\widetilde{\EE}^\vee:=(\HH\times \C)/\Z^2\quad (m,m')\cdot (\tau,z)= (\tau,z+m-\tau m') \quad (m,m')\in \Z^2,
$$
for the $\SL_2(\Z)$-action on $\widetilde{\EE}^\vee$ determined by $(\tau,z)\mapsto ((a\tau+b)/(c\tau+d),z/(c\tau+d))$. We observe that $\EE^\vee$ is a complex analytic group over $\Mell$, where the multiplication on $\EE^\vee$ is inherited from addition in $\C$. There is an $\SL_2(\Z)$-equivariant diffeomorphism
\beq
\HH\times (\R^2/\Z^2)=\HH\times S^1\times S^1 \simeq \widetilde{\EE}^\vee,\qquad (\tau,u,v)\mapsto (\tau,u-\tau v)\label{eq:squaretorus}
\eeq
for the $\SL_2(\Z)$-action on $\HH\times \R^2$ given by $(x,y)\mapsto (du-bv,-cu+av)$. Any homomorphism $G\to S^1$ determines a map $G^{(2)}\to S^1\times S^1$. Using~\eqref{eq:squaretorus}, this yields a a functor between complex analytic groupoids,
\beq
\Bun_G(\EE)\to \EE^\vee. \label{eq:GBundual}
\eeq
In the special case that $G=\Z/n$ and $\mu_n\colon \Z/n\hookrightarrow U(1)\simeq S^1$ includes as the $n$th roots of unity,~\eqref{eq:GBundual} is injective with image the $n$-torsion points on the universal dual elliptic curve. 

We promote the functor $\Bun_{\Z/n}(\EE)\to \EE^\vee$ to a homomorphism of complex analytic groups over $\Mell$ as follows. Naturality~\eqref{eq:ellnaturality} applied to the projection and multiplication homomorphisms $\Z/n\times \Z/n\to \Z/n$ give maps
$$
m,p_1,p_2\colon \Bun_{\Z/n}(\EE)\times_{\Mell} \Bun_{\Z/n}(\EE)\simeq \Bun_{\Z/n\times \Z/n}(\EE) \to \Bun_{\Z/n}(\EE),
$$
using the equivalence~\eqref{eq:ellproduct}. This endows $\Bun_{\Z/n}(\EE)$ with the structure of a complex analytic group over $\Mell$. Since the homomorphism $\mu_n\boxtimes \mu_n \colon \Z/n\times \Z/n \to S^1\times S^1$ is compatible with this structure,~\eqref{eq:GBundual}
determines an injective homomorphism 
\beq
\Bun_{\Z/n}(\EE)\hookrightarrow \EE^\vee\label{eq:ntorsinc}
\eeq
in the category of complex analytic groups over $\Mell$. Furthermore, the sheaf $\mathcal{O}_{\EE^\vee}$ restricts to $\Ell_{\Z/n}^0(\pt)$. This identifies the coefficient sheaf $\Ell_{\Z/n}^0(\pt)$ for $\Z/n$-equivariant elliptic cohomology with functions on the $n$-torsion points of $\EE^\vee$. 
\end{ex}

\begin{rmk}
Allowing the cyclic group to vary, the homomorphisms~\eqref{eq:ntorsinc} are compatible for injective homomorphisms $\Z/n\hookrightarrow \Z/m$. In particular, we get a system of complex analytic groups for the $\Z/p^k$ torsion in $\EE^\vee$ for a fixed prime $p$ and varying $k\in \N$. This is a complex analytic version of the $p$-divisible group associated to the group $\EE^\vee$. 
\end{rmk}

In view of the formula for the multiplicative group law~\eqref{eq:multFGL}, one might attempt to extract the elliptic group law on the $n$-torsion points. This formal group law depends on a choice of coordinate on~$\EE^\vee$, i.e., a function that vanishes to first order at the identity section, $\Mell \hookrightarrow \EE^\vee$. However, since the fibers of $\EE^\vee\to \Mell$ are compact complex manifolds, no global (holomorphic) coordinate exists. One solution to this problem is to consider a section of a holomorphic line bundle that vanishes to first order at the identity section. This leads one directly to the theory of theta functions and twisted Euler classes, as the next example describes.

\begin{ex}[Twisted Euler classes and Klein forms]\label{ex:Klein}
Consider the function on $\HH\times \C$
\beq
Z(\tau,u,v)&=&e^{-\pi i(uv-u)}q^{(v^2-v)/2}\prod_{n=1}^\infty \frac{(1-q^{n-v}e^{2\pi i u})(1-q^{n+v-1}e^{-2\pi i u})}{(1-q^n)^2}\nonumber\\
&=&e^{-\pi iuv}q^{v^2/2} (e^{\pi i z}-e^{-\pi i z})\prod_{n>0} (1-q^n e^{2\pi i z})(1-q^n e^{-2\pi i z})\label{eq:determinant}
\eeq
where $(u,v)\in \R^2$ are real coordinates with $(\tau,z)=(\tau,u-\tau v)\in \HH\times \C\simeq \HH\times \R^2$ as in~\eqref{eq:squaretorus}. The function~\eqref{eq:determinant} vanishes to first order at lattice points $\langle \tau,1\rangle\colon \HH\times \Z^2\hookrightarrow \HH\times \C$. Furthermore, $Z(\tau,u,v)$ transforms under the action of $\Z^2\rtimes \SL_2(\Z)$ by a nonvanishing function on $\HH\times \C$, and so determines a section of a line bundle on the quotient $[\HH\times \C\sq \Z^2\rtimes \SL_2(\Z)]\simeq \EE^\vee$. Indeed, this is Quillen's determinant line bundle for the $\bar\partial$-operator twisted by a degree zero holomorphic line bundle~\cite{Quillen_det}; see Freed~\cite[Equation~4.11]{Freed_det}. 
The restriction of this line bundle and section along~\eqref{eq:ntorsinc} is (essentially by definition, see~\cite{AndoOrientation}) the twisted equivariant elliptic Euler class of the representation $\mu_n\colon \Z/n\hookrightarrow U(1)$
\beq
\Eu(\mu_n):=\beta^{-1}Z(\tau,u,v)\in \Gamma(\Bun_{\Z/n}(\EE);\Ell^{2+\twist}_{\Z/n}(\pt)).\label{eq:twistedEuler}
\eeq
The twisting $\twist$ and the factor of $\beta^{-1}$ in~\eqref{eq:twistedEuler} is determined by the $\Z^2\rtimes \SL_2(\Z)$-transformation properties of~$Z(\tau,u,v)$. To spell these out, we recognize the restriction of~\eqref{eq:determinant} to the torsion points as a Klein form whose transformation properties are classical. Indeed, we have
\beq
Z(\tau,u+m,v+m')&=&(-1)^{mm'+m+m'}e^{\pi i(mv-m'u)}Z(\tau,u,v)\label{eq:Klein1}\\
 Z(\gamma\tau,u,v)&=&(c\tau+d)^{-1} Z(\tau,\gamma^{-1}(u,v)),\label{eq:Klein2}
\eeq
e.g., see Lang~\cite[Chapter XV]{Lang}. Euler classes are by definition multiplicative, so~\eqref{eq:twistedEuler} determines the Euler class for sums of representations,
\beq
\Eu(\mu_n^{\oplus k}):=\beta^{-k} Z(\tau,u,v)^k\in \Gamma(\Bun_{\Z/n}(\EE);\Ell^{2k+k\twist}_{\Z/n}(\pt))\label{eq:sumofreps}
\eeq
where $\mu_n^{\oplus k}$ is the representation associated to the composition $\Z/n\stackrel{\mu_n}{\to} U(1)\stackrel{\Delta}{\hookrightarrow } U(1)^{\times k}\subset \U(k)\subset \SO(2k)$ where $\Delta$ is the diagonal map. We also have the formula for the exterior product of representations,
\beq
\Eu(\mu_n^{\boxtimes k}):=\prod_{i=1}^k\beta^{-1} Z(\tau,u_i,v_i)\in \Gamma(\Bun_{(\Z/n)^{\times k}}(\EE);\Ell^{2k+\boxtimes \twist}_{(\Z/n)^{\times k}}(\pt))\label{eq:boxsumofreps}
\eeq
where $\mu_n^{\boxtimes k}$ is the representation associated with the composition $(\Z/n)^{\times k}\to U(1)^{\times k}\subset \U(K)\subset \SO(2k)$, and we use the fibered product descripton from~\eqref{eq:ellproduct},
$$
\Bun_{(\Z/n)^{\times k}}(\EE)\simeq \Bun_{\Z/n}(\EE)\times_{\Mell}\cdots\times_{\Mell} \Bun_{\Z/n}(\EE). 
$$
\end{ex}

To get a better understanding of the twisting $\alpha$ in the above examples, we recall that there are universal twisted equivariant elliptic Euler classes~\cite[\S6]{BET0} 
\beq
\Eu_{\SO(n)} \in \Gamma(\Bun_{\SO(n)}(\EE);\mathcal{L}),\qquad \Eu_{\U(n)}\in \Gamma(\Bun_{\U(n)}(\EE);\mathcal{L})\label{eq:univEuler}
\eeq
essentially determined by a product of the functions~\eqref{eq:determinant}. These Euler classes can be pulled back along maps 
$$
\Bun_G(\EE)\to \Bun_{\SO(n)}(\EE),\qquad \Bun_G(\EE)\to \Bun_{\U(n)}(\EE)
$$
determined by a real or complex representation $\rho\colon G\to \SO(n)$ or $\rho\colon G\to \U(n)$, respectively. The pull back of~$\mathcal{L}$ is trivializable precisely when the representation $\rho$ has a \emph{string structure}, defined as follows. For a representation $\rho$ with a chosen spin structure as in~\eqref{diagam:spinlift}, the fractional Pontryagin class $[\frac{p_1}{2}(\rho)]\in \H^3(BG;U(1))$ is the pullback of the class~$[\frac{p_1}{2}]\in \H^3(B\Spin(k);U(1))\simeq \H^4(B\Spin(k);\Z)$ along $B\tilde\rho$. If $[\frac{p_1}{2}(\rho)]=0$, then $\rho$ has a \emph{string structure}. Isomorphism classes of string structures are parameterized by~$\H^2(BG;U(1))$. If $G$ is connected, then $\rho$ has a string structure if and only if it is the trivial representation. Hence, the existence and moduli of string structures on a representation is only an interesting question for finite groups. We give a couple of examples of string obstructions and string structures for representations of finite groups.


\begin{ex}[A twisted Euler class for Conway's group]\label{ex:ConwayEuler}
Consider Conway's group ${\rm Co}_0$. The cohomology groups $\H^1(B{\rm Co}_0;\Z/2)\simeq \H^2(B{\rm Co}_0;\Z/2)\simeq 0$ vanish, so any representation of ${\rm Co}_0$ has a unique spin structure. Let $\Lambda\subset \R^{24}$ denote the Leech lattice and $\Lambda\otimes \R$ be the 24-dimensional representation of ${\rm Co}_0$ inherited from its action on the Leech lattice by automorphisms. It was shown in~\cite{TheoDavid} that $\frac{p_1}{2}(\Lambda\otimes \R^{24})$ is the generator of $\H^3(B{\rm Co}_0;\U(1))\simeq \Z/24$. Hence, the line bundle associated with $\Eu(\Lambda\otimes \R)\in \Gamma(\Bun_{{\rm Co}_0}(\EE);\Ell_{{\rm Co}_0}^{24+\frac{p_1}{2}}(\pt))$ is $\omega^{24/2}$ tensored with the transgressed line bundle considered in Example~\ref{ex:Conway}.
\end{ex} 

\begin{rmk}
Let $\M$ be the monster group. In~\cite{TheoAnomaly} it was shown the order~24 element $\omega^{\natural}\in \H^3(B\M;U(1))$ associated with the moonshine anomaly does not arise as $\frac{p_1}{2}(\rho)$ for any representation $\rho$ of $\M$. Consequently, there is no twisted $\M$-equivariant elliptic Euler class with twist classified by the moonshine anomaly. 
\end{rmk} 

\begin{ex}[An untwisted Euler class for $\Z/2$]\label{eq:Z2}
Using the characteristic class techniques described in Example~\ref{ex:Kthyz2}, sums of the sign representation $\sigma^{\oplus k}\simeq \R^k\sq \Z/2$ are string if $k$ is divisible by~$8$. Furthermore, after one fixes the spin structure on the representation, the isomorphism class of string structure is unique since~$\H^2(B\Z/2;U(1))\simeq 0$. This gives untwisted Euler classes
$$
\Eu(\sigma^{\oplus 8k})=\Eu(\mu_2^{\oplus 4}) \in \Gamma(\Bun_{\Z/2}(\EE);\Ell^{8k}_{\Z/2}(\pt))
$$
where, as before, we identify the representations $\sigma\oplus\sigma\simeq \mu_2$ using $\R^2\simeq \C$. Hence, $\sigma^{\oplus 8k}$ is the composition
$$
\sigma^{\oplus 8k}\colon \Z/2\hookrightarrow (S^1)^{\times 4k}\simeq (\SO(2))^{\times 4k}\hookrightarrow {\rm O}(8k)
$$
for $k\in \N$. We will verify explicitly that $\Eu(\sigma^{\oplus 8k})$ is indeed an untwisted class; it will suffice to show for $k=1$. From the definition of the Euler class from~\eqref{eq:twistedEuler}, $\Eu(\mu_2)=\Eu(\sigma\oplus\sigma)$ is determined by the values of $Z(\tau,u,v)$ at 2-torsion points $\Bun_{\Z/2}(\EE)\hookrightarrow \EE^\vee$. We identify these values with the Jacobi theta functions as follows
\beq
Z(\tau,0,0)&=&0,\nonumber\\
Z(\tau,1/2,0)&=&2e^{\pi i/2}\prod_{n>0} \frac{(1+q^n)(1+q^{n})}{(1-q^n)^2}=e^{\pi i/2} \frac{\vartheta_2(\tau)}{\eta(\tau)^3},\nonumber\\
Z(\tau,0,1/2)&=&q^{-1/8}\prod\frac{(1-q^{n-1/2})^2}{(1-q^n)^2}=\frac{\vartheta_4(\tau)}{\eta(\tau)^3},\nonumber\\
Z(\tau,1/2,1/2)&=&e^{\pi i/4}q^{-1/8}\prod\frac{(1+q^{n-1/2})^2}{(1-q^n)^2}=e^{\pi i/4} \frac{\vartheta_3(\tau)}{\eta(\tau)^3}.\nonumber
\eeq
The transformation properties of $Z(\tau,u,v)$ from~\eqref{eq:Klein1} and~\eqref{eq:Klein2} (or equivalently from the Jacobi theta functions) show that the 4th powers
\beq
\Eu(\sigma^{\oplus 8})(\tau,1/2,1/2)= -(\vartheta_3/\eta^3)^4,& & \Eu(\sigma^{\oplus 8})(\tau,0,1/2)= (\vartheta_4/\eta^3)^4,\nonumber\\
 \Eu(\sigma^{\oplus 8})(\tau,1/2,0)= (\vartheta_2/\eta^3)^4,&& \Eu(\sigma^{\oplus 8})(\tau,0,0)= 0\nonumber
\eeq
give an $\SL_2(\Z)$-invariant function on 2-torsion points. Therefore for $k\in \N$ we have canonically (untwisted) Euler classes 
$$
\Eu(\sigma^{\oplus 8k})\in \Gamma(\Bun_{\Z/2}(\EE);\Ell^{8k}_{\Z/2}(\pt)). 
$$

\end{ex}




\section{Supermanifolds and super Lie groupoids}\label{sec:super}

\subsection{Supermanifolds}
To fix conventions and introduce some relevant examples, we give a brief survey of supermanifolds. Comprehensive introductions include~\cite{DM,Leites}. 
\begin{defn}
A \emph{$k|l$-dimensional supermanifold} is a Hausdorff, second countable, locally ringed space that is locally isomorphic to $(U,C^\infty_U\otimes \Lambda^\bullet \C^l)$ where $U\subset \R^k$ is an open submanifold, and $C^\infty_U\otimes \Lambda^\bullet \C^l$ is regarded as a sheaf of commutative super algebras. Let ${\sf SMfld}$ denote the category whose objects are supermanifolds and morphisms are maps of locally ringed spaces. 
\end{defn}

Supermanifolds as defined above (with structure sheaves over~$\C$) are sometimes called \emph{$cs$-manifolds}~\cite[\S4.8]{DM}. 

\begin{ex}\label{example1}The sheaf of $\C$-valued smooth functions allows one to view an ordinary manifold~$M$ as a supermanifold. More generally, for a complex vector bundle $E\to M$, let $\Pi E$ denote the supermanifold whose underlying topological space is $M$ with sheaf of rings the sheaf of sections of $\Lambda^\bullet E^\vee$. For $\underline{\C}^k\to \R^l$ the trivial bundle, define $\R^{k|l}:=\Pi \underline{\C}^k$. 
\end{ex}

For a supermanifold $N$, let $C^\infty(N)$ denote the global sections its structure sheaf and call these the \emph{functions} on $N$. In the notation of Example~\ref{example1}, 
\beq\label{Eq:PiE}
C^\infty(\Pi E)\simeq \Gamma(M,\Lambda^\bullet E^\vee).
\eeq
By a partition of unity argument, supermanifolds are \emph{affine}: a map $N\to N'$ of supermanifolds is determined by a map $C^\infty(N')\to C^\infty(N)$ of superalgebras. 

The geometry of supermanifolds can be subtle because $C^\infty(N)$ typically has nilpotent elements. In parallel to the study of schemes, define the \emph{reduced manifold} $N_{\rm red}$ of a supermanifold~$N$ as the locally ringed space with the same underlying space, and sheaf of rings given by the quotient by the nilpotent ideal, e.g., global sections are $C^\infty(N_{\rm red})\simeq C^\infty(N)/{\rm nil}$. This gives $N_{\rm red}$ the structure of an ordinary smooth manifold. Any supermanifold~$N$ is isomorphic to one of the form $\Pi E$ for $E\to N_{\rm red}$ a complex vector bundle over an ordinary manifold~\cite{batchelor}. The category of supermanifolds is interesting precisely because there are more morphisms $\Pi E\to \Pi E'$ than there are vector bundle maps $E\to E'$.

The functor of points is another tool in supermanifolds imported from scheme theory. For a supermanifold $N$ and a (test) supermanifold $S$, the \emph{$S$-points of $N$} are defined as the set $N(S):={\sf SMfld}(S,N)$. Varying $S$, these sets assemble into a presheaf $N\colon {\sf SMfld}^{\op}\to {\sf Set}$. By the Yoneda lemma the original supermanifold $N$ is completely determined by this presheaf. 

\begin{ex}
The $S$-points of $\R^{n|m}$ are given by 
$$
\R^{n|m}(S)\simeq \{x_1,\dots x_n \in C^\infty(S)^\ev,\ \theta_1,\dots, \theta_m \in C^\infty(S)^\odd\mid (x_i)_{\rm red}=\overline{(x_i)}_{\rm red}\}\nonumber\\
$$
where $(x_i)_{\rm red}$ denotes the image of $x_i$ along $C^\infty(S)\to C^\infty(S)/{\rm nil}=C^\infty(S_{\rm red})$; being a smooth manifold, $\C$-valued functions on $S_{\rm red}$ have a well-defined complex conjugate $\overline{(x_i)}_{\rm red}$. 
\end{ex}

\begin{rmk} A function on a supermanifold does not usually have a (canonical) complex conjugate. In the description from~\eqref{Eq:PiE} such a conjugation depends on a choice of real structure on the complex vector bundle~$E$, i.e., an isomorphism $E\simeq E_\R\otimes \C$ for some real vector bundle $E_\R\to N_{\rm red}$ (which need not exist). \end{rmk} 

\begin{ex}\label{ex:R21}
We require two descriptions of $S$-points of~$\R^{2|1}$
\beq
\R^{2|1}(S)&\simeq& \{x,y \in C^\infty(S)^\ev,\ \theta \in C^\infty(S)^\odd\mid (x)_{\rm red}=\overline{(x)}_{\rm red}, (y)_{\rm red}=\overline{(y)}_{\rm red}\}\label{eq:r211}\\
&\simeq& \{z,w \in C^\infty(S)^{\ev}, \theta\in C^\infty(S)^\odd\mid (z)_{\rm red}=\overline{(w)}_{\rm red}\}.\label{eq:r212}
\eeq
The isomorphism between~\eqref{eq:r211} and~\eqref{eq:r212} is $(x,y)\mapsto (x+iy,x-iy)=(z,w)$. When working with~\eqref{eq:r212} we shall adopt the standard (though potentially misleading) notation $\overline{z}:=w$. We emphasize that $\overline{z}$ is only the conjugate of $z$ after restriction to~$S_{\rm red}$. We use similar notation for $S$-points of complex manifolds regarded as supermanifolds, e.g., $(u,\bar u)\in \C^\times (S)$. 
\end{ex}

Sometimes geometric objects naturally determine non-representable presheaves on the site of supermanifolds. To emphasize these geometric origins, we often call a presheaf on supermanifolds a \emph{generalized supermanifold}. 


\begin{defn}\label{defn:map}
For supermanifolds $N$ and $M$, define the generalized supermanifold $\Map(N,M)$ whose values are
\beq
\Map(N,M)(S):={\sf SMfld}(S\times N,M).\label{ex:mapssmfld}
\eeq
These sets pull back along base changes $b\colon S'\to S$ by precomposing with $b\times \id\colon S'\times N\to S\times N$. 
\end{defn}

The generalized supermanifold $\Map(N,M)$ fails to be representable for essentially the same reasons that maps between ordinary manifolds generally do not form a (finite-dimensional) manifold. There is an exception when~$N$ is zero dimensional, and similarly the presheaf $\Map(N,M)$ on supermanifolds is representable when $N$ is $0|l$-dimensional. The following example will be crucial in our constructions below. 

\begin{ex}\label{ex:PiT}
For an ordinary smooth manifold $M$ regarded as a supermanifold, consider the mapping object $\Map(\R^{0|1},M)$. This generalized supermanifold is representable, $\Map(\R^{0|1},M)\simeq \Pi TM$, e.g., see~\cite[Proposition~3.1]{HKST}. In particular, functions 
$$
C^\infty(\Map(\R^{0|1},M))\simeq C^\infty(\Pi TM)\simeq \Gamma(M,\Lambda^\bullet TM^\vee)\simeq \Omega^\bullet(M)
$$
are ($\C$-valued) differential forms as a $\Z/2$-graded algebra. 
\end{ex}

\begin{defn}\label{defn:vb} A \emph{vector bundle} on a supermanifold is a locally free sheaf of $\Z/2$-graded modules over the structure sheaf. The \emph{parity reversal} of a vector bundle $V\to N$ is denoted by $\Pi V\to N$, and is defined as the same sheaf of modules with the opposite $\Z/2$-grading. A \emph{line bundle} is a rank one locally free sheaf.  \end{defn}

By a partition of unity argument, a vector bundle on a supermanifold $N$ is determined by a module over the global sections~$C^\infty(N)$. 

\begin{ex}\label{ex:lineonsuper}
For any supermanifold $N$, the trivial line bundle $\underline{\C}\to N$ is given by~$C^\infty(N)$ as a module over itself. The trivial \emph{odd} line bundle $\Pi \underline{\C}\to N$ is $C^\infty(N)$ regarded as a module with the opposite grading. We observe that automorphisms of the trivial (odd) line bundle are given by $(C^\infty(N)^\times)^{\ev}$, the invertible even functions on $M$. 
\end{ex}

\begin{ex} A \emph{vector field} $X$ on a supermanifold~$N$ is a derivation $X\colon C^\infty(N)\to C^\infty(N)$. Such derivations can have even or odd parity, and we call the vector field even or odd accordingly. Working locally on~$N$, derivations form a sheaf of modules over the structure sheaf. This sheaf of modules is the \emph{tangent bundle}~$TN$ of the supermanifold~$N$, which is an example of a vector bundle on~$N$. Just as for ordinary manifolds, a map $f\colon N\to N'$ of supermanifolds determines a map $df\colon TN\to TN'$ on tangent bundles. This allows one to import notions like immersions and submersions to the category of supermanifolds, defined in terms of properties of the map of vector bundles~$df$. 
\end{ex}

\subsection{Super Lie groups}
A \emph{super Lie group} is a group object in supermanifolds. Explicitly, a super Lie group is a supermanifold~$G$ where $G(S)$ is endowed with the structure of a group for each $S$, and the maps of sets $G(S')\to G(S)$ associated to maps of supermanifolds~$S\to S'$ are group homomorphisms. 

\begin{ex}\label{ex:confgroups} Define the super Lie groups $\E^{0|1},\E^{1|1}$ and $\E^{2|1}$ whose underlying supermanifolds are $\R^{0|1}, \R^{1|1}$ and $\R^{2|1}$, respectively, and whose multiplications $m\colon \E^{d|1}\times \E^{d|1}\to \E^{d|1}$ are determined by the $S$-point formulas
\beq
m(\theta, \eta)&=&(\theta+\eta)\qquad \theta,\eta\in \E^{0|1}(S),\nonumber\\
m((t,\theta),(s,\eta))&=&(t+s+\theta\eta,\theta+\eta)\qquad (t,\theta),(s,\eta)\in \E^{1|1}(S),\nonumber\\
m((z,\bar z,\theta),(w,\bar w,\eta))&=&(z+w,\bar z+\bar w+\theta\eta,\theta+\theta)\qquad (z,\bar z,\theta),(w,\bar w,\eta)\in \E^{2|1}(S), \nonumber
\eeq
where the final formula uses the notation from~\eqref{eq:r212}. Define the semidirect products $\E^{0|1}\rtimes \R^\times$, $\E^{0|1}\rtimes \C^\times$, $\E^{1|1}\rtimes \R^\times$ and $\E^{2|1}\rtimes \C^\times$ for the actions
\beq
\begin{array}{rclrclll}
\dil \cdot \theta&=&\dil\theta, & \dil\cdot(t,\theta)&=&(\dil^2t,\dil\theta), & \dil\in \R^\times(S),\\
(\dil,\bar\dil)\cdot \theta&=&\bar\dil\theta,& (\dil,\bar\dil)\cdot (z,\bar z,\theta)&=&(\dil^2z,\bar\dil^2\bar z,\bar\dil\theta),& (\dil,\bar\dil)\in \C^\times(S).
\end{array}\nonumber
\eeq
There are (left) action maps 
\beq
&&\E^{2|1}\rtimes \C^\times\times \R^{2|1}\to \R^{2|1},\quad \E^{1|1}\rtimes \R^\times\times \R^{1|1}\to \R^{1|1},\quad \E^{0|1}\rtimes \C^\times \times \R^{0|1}\to \R^{0|1}.\label{eq:actionmapsuper}
\eeq
We also observe that the projections define homomorphisms
\beq
\E^{2|1}\rtimes \C^\times \to \E^{0|1}\rtimes \C^\times,\qquad \E^{1|1}\rtimes \R^\times\to \E^{0|1}\rtimes \R^\times\label{eq:projhomo}. 
\eeq
\end{ex}

\begin{ex} Let $\Aut(\R^{0|1})$ denote the presheaf whose $S$-points are isomorphisms $S\times \R^{0|1}\to S\times \R^{0|1}$ over $S$. Composition of isomorphisms endows $\Aut(\R^{0|1})$ with the structure of a group object in presheaves. There are homomorphisms of group objects in presheaves
\beq
\E^{0|1}\rtimes \C^\times\to \Aut(\R^{0|1})\qquad \E^{0|1}\rtimes \R^\times\to \Aut(\R^{0|1})\label{eq:autE01}
\eeq
using the action on $\R^{0|1}$ from~\eqref{eq:projhomo}. 
\end{ex}

\begin{ex}\label{eq:deRham}
The presheaf $\Map(\R^{0|1},M)$ from Example~\ref{ex:PiT} has a left $\Aut(\R^{0|1})$-action by precomposing a map with an automorphism. The homomorphisms~\eqref{eq:autE01} give left actions of $\E^{0|1}\rtimes \C^\times$ and $\E^{0|1}\rtimes \R^\times$ on $\Map(\R^{0|1},M)$. The derivative at zero of each of these actions gives a pair of derivations on $C^\infty(\Map(\R^{0|1},M))\simeq \Omega^\bullet(M)$: the $\E^{0|1}$-action is generated by the odd derivation~$-d$, and the $\R^\times$- or $\C^\times$-action is generated by the even derivation $-\deg$, e.g., see~\cite[Proposition~3.9]{HKST}. Here $d$ is the de~Rham differential and $\deg(\omega)=k\omega$ for $\omega\in \Omega^k(M)$. The minus signs on these derivations come from turning the precomposition action (usually a right action) into a left action. Explicitly, the action maps are given by
$$
C^\infty(\Map(\R^{0|1},M))\to C^\infty(\Map(\R^{0|1},M)\times \E^{0|1}\rtimes \R^\times),\quad f\mapsto u^{-\deg} (f-\eta df)
$$
$$ 
C^\infty(\Map(\R^{0|1},M))\to C^\infty(\Map(\R^{0|1},M)\times \E^{0|1}\rtimes \C^\times),\quad f\mapsto \bar u^{-\deg} (f-\eta df)
$$
where in the first line, $(\eta,u)$ are coordinates on $\E^{0|1}\rtimes \R^\times$, and in the second line $(\eta,u,\bar u)$ are coordinates on $\E^{0|1}\rtimes \C^\times$. 

\end{ex}

We will often construct actions of a super Lie group $G$ on a supermanifold $N$ in terms group actions (in sets) $G(S)\times N(S)\to N(S)$ that are natural in the supermanifold~$S$. 

\begin{ex}
In terms of the functor of points, the left action in Example~\ref{eq:deRham} is determined by maps of sets $\Aut(\R^{0|1})(S)\times \Map(\R^{0|1},M)(S)\to \Map(\R^{0|1},M)(S)$ defined as $(g,\phi)\mapsto (\phi\circ g^{-1})$ for the composition
$$
S\times \R^{0|1}\stackrel{g^{-1}}{\to} S\times \R^{0|1}\stackrel{\phi}{\to} M,\qquad \phi\in \Map(\R^{0|1},M)(S), \ g\in \Aut(\R^{0|1})(S).
$$
\end{ex}

\subsection{Super Lie groupoids}\label{sec:superLie}
\begin{defn}
A \emph{super Lie groupoid} $\mathcal{G}=\{\mathcal{G}_1\rightrightarrows \mathcal{G}_0\}$ consists of data~\eqref{eq:grpd} where~$\mathcal{G}_1$ and~$\mathcal{G}_0$ are supermanifolds and the structure maps $s,t,u,c$ are maps between supermanifolds. We further require that $s$ be a submersion so that $\mathcal{G}_1\times_{\mathcal{G}_0}\mathcal{G}_1$ exists in the category of supermanifolds. As usual, these data are required to satisfy the axioms for a groupoid. 
\end{defn}
\begin{defn}
 A \emph{functor} between super Lie groupoids $F\colon \mathcal{G}\to \mathcal{H}$ is the data of maps of supermanifolds $F_1\colon \mathcal{G}_1\to \mathcal{H}_1$, $F_0\colon \mathcal{G}_0\to \mathcal{H}_0$ satisfying the axioms for a functor. A \emph{natural transformation} between functors $F,F'\colon \mathcal{G}\rightrightarrows \mathcal{H}$ is a map of supermanifolds $\mathcal{G}_0\to \mathcal{H}_1$ satisfying the axioms of a natural transformation $F\Rightarrow F'$. 
\end{defn}

Super Lie groupoids, functors, and natural transformations form a strict 2-category. 

\begin{ex} Any supermanifold $N$ determines a (discrete) super Lie groupoid $\{N\rightrightarrows N\}$ with only identity morphisms. Typically we use the same notation~$N$ for the discrete super Lie groupoid associated to a supermanifold~$N$. \end{ex}

\begin{ex}
For a super Lie group $G$ acting on a supermanifold $N$, the \emph{action groupoid} $N\sq G$ has objects $N$ and morphisms $G\times N$. The source map is projection and target map is the action. The unit and composition are determined by the unit and multiplication in~$G$. 
\end{ex}

\begin{defn}
A \emph{generalized super Lie groupoid} is a groupoid object in the category of presheaves on supermanifolds, i.e., $\mathcal{G}_1,\mathcal{G}_0$ are presheaves and the structure maps $s,t,u,c$ as in~\eqref{eq:grpd} are maps of presheaves. 
\end{defn}
Natural constructions involving super Lie groupoids often land in generalized super Lie groupoids, as the following example shows. 

\begin{ex}
For super Lie groupoids $\mathcal{G}$ and $\mathcal{H}$, let $\Fun(\mathcal{G},\mathcal{H})$ denote the generalized super Lie groupoid whose object presheaf assigns to a supermanifold $S$ the set of functors $S\times\mathcal{G}\to \mathcal{H}$, and whose morphism presheaf assigns to $S$ the set of natural transformations between functors. Pulling back functors and natural transformations along base changes $S'\to S$ promotes these values at each~$S$ to presheaves. These presheaves are usually not representable, just as the mapping objects~\eqref{ex:mapssmfld} between supermanifolds are usually not representable. We will use the notation
$$
\Map(\mathcal{G},\mathcal{H}):={\rm Obj}(\Fun(\mathcal{G},\mathcal{H}))
$$
to denote the object presheaf of $\Fun(\mathcal{G},\mathcal{H})$. 
\end{ex}

\begin{defn}\label{defn:equivalence} A functor $\mathcal{G}\to \mathcal{H}$ between (generalized) super Lie groupoids is an \emph{equivalence} if for each supermanifold~$S$, the functor $\mathcal{G}(S)\to \mathcal{H}(S)$ between groupoids internal to sets is fully faithful and essentially surjective. \end{defn}

An equivalence of super Lie groupoids need not admit an inverse functor, as the following example shows. 

\begin{ex}
Suppose a super Lie group $G$ acts freely on a supermanifold $M$. Then the quotient $M/G$ is a supermanifold. The canonical functor $q\colon M\sq G\to M/G$ is an equivalence of super Lie groupoids, viewing the target as a discrete super Lie groupoid. The functor $q$ admits an inverse if and only if $M\to M/G$ is a trivial principal $G$-bundle, i.e., there exists a smooth section $M/G \to M$.
\end{ex}

\begin{rmk} Localization of the 2-category of super Lie groupoids, functors and natural transformations with respect to equivalences leads to the bicategory of geometric super stacks. Just as in the case of smooth stacks, this localization has an explicit description in terms of super Lie groupoids, bibundles and maps between bibundles, e.g., see~\cite[\S2.2]{SchommerPries}.
\end{rmk} 

We define functions, vector bundles, and sheaves on super Lie groupoids identically to the Definitions~\ref{defn:grpdfun},~\ref{defn:grpdvb} and~\ref{defn:grpdsheaf}. We give a feel for these definitions in the following examples. 

\begin{ex} A function $f\in C^\infty(\mathcal{G})$ on a super Lie groupoid $\mathcal{G}$ is determined by the data of a function $f\in C^\infty(\mathcal{G}_0)$ on its supermanifold of objects with the property $s^*f=t^*f\in C^\infty(\mathcal{G}_1)$ on the supermanifold of morphisms. If $\mathcal{G}=N\sq G$ is an action groupoid, then $C^\infty(N\sq G)=C^\infty(N)^G$ is the set of $G$-invariant functions. A function on a super Lie groupoid is \emph{even} or \emph{odd} according to the parity of the function on~$\mathcal{G}_0$. 
\end{ex}

\begin{ex} \label{ex:canonicalline} A line bundle on a super Lie groupoid $\mathcal{G}$ is a line bundle $\mathcal{L}\to \mathcal{G}_0$ on the supermanifold of objects and an isomorphism of line bundles $s^*\mathcal{L}\simeq t^*\mathcal{L}$ on the supermanifold of morphisms $\mathcal{G}_1$ satisfying a cocycle condition over $\mathcal{G}_1\times_{\mathcal{G}_0}\mathcal{G}_1$. A line bundle on $N\sq G$ is a $G$-equivariant line bundle on~$N$. The super Lie groupoid $\pt\sq \C^\times$ comes equipped with a canonical line bundle specified by the trivial line bundle on $\pt$ (which is just the 1-dimensional vector space $\C$ regarded as a module over $C^\infty(\pt)\simeq \C$), and the automorphism of the trivial line bundle on $\C^\times$ given by the (invertible, even) coordinate function $z\in C^\infty(\C^\times)$ (see Example~\ref{ex:lineonsuper}).
\end{ex}

\begin{ex}
By the previous example, a functor $F\colon \mathcal{G} \to \pt\sq \C^\times$ of super Lie groupoids determines a line bundle on $\mathcal{G}$ by pull back of the canonical line $\mathcal{L}\to \pt\sq \C^\times$. The pullback $F^*\mathcal{L}$ is trivial over the objects $\mathcal{G}_0$ and has the data of an automorphism of the trivial line bundle on $\mathcal{G}_1$. The global sections of such a line bundle are given by
\beq
\Gamma(\mathcal{G};F^*\mathcal{L})=\{f\in C^\infty(\mathcal{G}_0)^\ev \mid t^*f=(F_1^*z)s^*f\in C^\infty(\mathcal{G}_1)\}\label{eq:grpodsections}
\eeq
where $F_1\colon \mathcal{G}_1\to \C^\times$ is the functor $F$ on morphisms and $z\in C^\infty(\C^\times)$ is the coordinate function. Sections of the pullback of the parity reversal $\Pi\mathcal{L}$ of the canonical line along $F\colon \mathcal{G}\to \pt\sq \C^\times$ are given by~\eqref{eq:grpodsections}, but where $f\in C^\infty(\mathcal{G}_0)^\odd$. 
\end{ex}

\begin{ex}\label{Ex:differentialdescent}
Sheaves on super Lie groupoids are defined completely analogously to line bundles (compare Example~\ref{defn:grpdsheaf}). There is an important relationship between $\E^{0|1}$-equivariant sheaves and chain complexes, as we now explain. Suppose we are given a sheaf $\F$ of $C^\infty(X)$-modules on a supermanifold $X$. Equip $X$ with the trivial $\E^{0|1}$-action. Then descent data for $\F$ along the map $X\to X\sq \E^{0|1}$ is an isomorphism $p^*\F\simeq \mu^*\F$ over $\E^{0|1}\times X$ where $p$ is the projection and $\mu$ is the action. Since the action is trivial, $p=\mu$. Hence, the descent data is equivalent to an $\E^{0|1}$-action on the sections of $\F$ that is linear over $C^\infty(X)$. Differentiating at $0\in \E^{0|1}$, such an action is the data of an odd, square zero, $C^\infty(X)$-linear endomorphism of~$\F$. This identifies $\F$ with a $\Z/2$-graded complex of sheaves in~$C^\infty(X)$-modules. 
\end{ex}

\begin{ex}\label{ex:deRhamsheaf}
The $\E^{0|1}\rtimes \C^\times$-action on $\Map(\R^{0|1},M)$ determines the action super Lie groupoid $\Map(\R^{0|1},M)\sq (\E^{0|1}\rtimes \C^\times)$. From the description of this action in Example~\ref{eq:deRham}, functions are
$$
C^\infty(\Map(\R^{0|1},M)\sq (\E^{0|1}\rtimes \C^\times))\simeq C^\infty(\Map(\R^{0|1},M))^{\E^{0|1}\rtimes \C^\times}\simeq (\Omega^\bullet_\cl(M))^{\C^\times}\simeq \Omega^0_\cl(M),
$$
i.e., the closed $0$-forms on $M$. We can also consider the sheaf of functions on the super Lie groupoid $\Map(\R^{0|1},M)\sq (\E^{0|1}\rtimes \C^\times)$. This is given by the sheaf of functions on $\Map(\R^{0|1},M)$ (which we identify with the sheaf of differential forms on~$M$) together with the isomorphism between the pullback along source and target maps given by the action of $\E^{0|1}\rtimes \C^\times$ on differential forms. This action is generated by the operators $-d$ and $-\deg$. Hence, the sheaf of functions on $\Map(\R^{0|1},M)\sq (\E^{0|1}\rtimes \C^\times)$ encodes the entire de~Rham complex of~$M$. To pick out closed forms in a fixed degree (not necessarily zero), observe that there is a functor
$$
\Map(\R^{0|1},M)\sq (\E^{0|1}\rtimes \C^\times)\to \pt\sq \C^\times,
$$
determined by the canonical map $\Map(\R^{0|1},M)\to \pt$ and the projection $\E^{0|1}\rtimes \C^\times\to \C^\times$. We define a line bundle $\mathcal{L}\to \Map(\R^{0|1},M)\sq (\E^{0|1}\rtimes \C^\times)$ as the pullback of the parity reversal of the line from Example~\ref{ex:canonicalline}, i.e., the pullback of the canonical \emph{odd} line bundle on $\pt\sq \C^\times$. Sections of tensor powers of this line bundle are computed in~\cite[\S6]{HKST}, using the description~\eqref{eq:grpodsections}. The result is the vector space of closed differential forms of degree $k$, 
$$
\Gamma(\Map(\R^{0|1},M)\sq (\E^{0|1}\rtimes \C^\times);\mathcal{L}^{\otimes k})\simeq \Omega^k_\cl(M).
$$
\end{ex}

\section{Principal $G$-bundles over super circles and super tori}\label{sec:circlesandtori}

\subsection{Rigid conformal geometries}
In the next two definitions, we briefly review the (fibered) category whose objects are families of supermanifolds with a rigid geometry from~\cite[\S2.5 and~\S4.1]{ST11}; we refer to this reference for details. 
\begin{defn} 
A \emph{rigid geometry} is the data of a supermanifold $\M$ and a super Lie group $\Iso(\M)$ acting on $\M$. We call $\M$ the \emph{model space} and $\Iso(\M)$ the \emph{group of isometries} for the rigid geometry.
\end{defn} 

\begin{defn}[\cite{ST11}, Definitions~2.33 and~4.4] An $S$-family of supermanifolds with $(\M,\Iso(\M))$-structure is the data of 
\begin{enumerate}
\item a smooth submersion $p\colon T\to S$;
\item a maximal atlas $(U_i)$ of $T$ with charts equipped with isomorphisms over $S$, $\varphi_i \colon U_i \xrightarrow{\sim} V_i \subset S \times \M$ where $V_i$ is open;
\item transition data $g_{ij}\colon p(U_i\bigcap U_j)\to \Iso(\M)$. 
\end{enumerate}
The isomorphisms $\varphi_i$ are required to be compatible with transition data $g_{ij}$ and the transition data must satisfy a further cocycle condition. An \emph{isometry} between $S$-families of supermanifolds with $(\M,\Iso(\M))$-structure is a map $T\to T'$ over $S$ that on an open cover is determined by the action of~$\Iso(\M)$ on the open sub supermanifolds $V_i\subset \M$, satisfying a compatibility condition with the $(\M,\Iso(\M))$ structures on $T\to S$ and $T'\to S$. 
\end{defn}

\begin{rmk}
Fixing $(\M,\Iso(\M))$, there is a fibered category (in fact, a stack) whose objects are supermanifolds with $(\M,\Iso(\M))$-structure and morphisms are isometries between $S$-families~\cite[\S2]{ST11}. 
\end{rmk} 

\begin{defn}
Define $1|1$-dimensional and $2|1$-dimensional \emph{rigid conformal geometries} by
$$
(\M,\Iso(\M))=(\R^{1|1},\E^{1|1}\rtimes \R^\times),\qquad (\M,\Iso(\M))=(\R^{1|1},\E^{2|1}\rtimes \C^\times)
$$
respectively, using the left action of $\E^{1|1}\rtimes \R^\times$ on $\R^{1|1}$, and the left action of $\E^{2|1}\rtimes \C^\times$ on~$\R^{2|1}$ as defined in Example~\ref{ex:confgroups}. 
\end{defn} 

We adopt the following terminology: an $(\M,\Iso(\M))$-structure as above will be referred to as a \emph{$d|1$-dimensional rigid conformal structure} and a family with such a structure is a family of \emph{$d|1$-dimensional rigid conformal supermanifolds}. Finally, a morphism between families with $(\M,\Iso(\M))$-structure as above is a family of \emph{rigid conformal maps}. 


\begin{defn} \label{defn:supercircle}
Let $\R_{>0}\subset \R$ be the open submanifold of positive real numbers regarded as a supermanifold. Define the quotient, 
\beq
&&\mathcal{T}^{1|1}:=(\R_{>0}\times \R^{1|1})/\Z, \quad n\cdot (t,\theta)=(t+n\ell ,\theta)\quad n\in \Z(S),\ (t,\theta)\in \R^{1|1}(S)\label{eq:aZaction}
\eeq
for the $\Z$-action defined by the functor of points formula. 
\end{defn}

\begin{lem}\label{lem:rigidconfstru}The family $\mathcal{T}^{1|1}\to\R_{>0}$ has a canonical $1|1$-dimensional rigid conformal structure, defining an $\R_{>0}$-family of $1|1$-dimensional rigid conformal supermanifolds. \end{lem}
\bp
Giving $\R_{>0}\times \R^{1|1}=\R_{>0}\times \M$ the canonical structure of an $\R_{>0}$-family of $1|1$-dimensional rigid conformal manifolds, the $\Z$-action~\eqref{eq:aZaction} is through rigid conformal isometries. Furthermore, the quotient map by the $\Z$-action is a local isomorphism. Hence, the quotient $\mathcal{T}^{1|1}$ has a canonical rigid conformal structure. 
\ep

Let $\Lat\subset \C^\times \times \C^\times$ denote the subset $(\ell_1,\ell_2)\in \C^\times\times \C^\times$ such that $
\ell_2/\ell_1\in \HH$ where $\HH\subset \C$ is the upper half plane. The isomorphism $\Lat\simeq \C^\times\times \HH$ given by $(\ell_1,\ell_2)\mapsto (\ell_1,\ell_2/\ell_1)$ shows that $\Lat\subset \C^\times\times \C^\times$ is an open submanifold. When viewing $\Lat$ as a supermanifold, an $S$-point will be denoted by~$\ell\in (\ell_1,\bar\ell_1,\ell_2,\bar\ell_2)\in \Lat(S)$, following the complex conjugate notation discussed in Example~\ref{ex:R21}. 

\begin{defn} \label{defn:supertorus}
Define the quotient supermanifold,
$$
\mathcal{T}^{2|1}:=(\Lat \times \R^{2|1})/\Z^2
 $$
 for $\Z^2$-action on $\Lat\times \R^{2|1}$ specified by the functor of points formula
\beq
(n,m)\cdot(z,\bar z,\theta)=(z+n\ell_1+m\ell_2,\bar z+n\bar\ell_1+m\bar\ell_2,\theta),\\
(\ell_1,\bar\ell_1,\ell_2,\bar\ell_2)\in \Lat(S), \ (z,\bar z,\theta)\in \R^{2|1}(S), \ (n,m)\in \Z^2(S).\nonumber \label{eq:Z2action}
\eeq
\end{defn}

\begin{lem}The family $\mathcal{T}^{2|1}\to\Lat$ has a canonical $2|1$-dimensional rigid conformal structure, defining a $\Lat$-family of $2|1$-dimensional rigid conformal supermanifolds. \end{lem}

The proof of the above lemma is identical to that of Lemma~\ref{lem:rigidconfstru}. 

\begin{notation}
For a map $\ell\colon S\to \R_{>0}$ or $\ell\colon S\to \Lat$, we adopt the notation
$$
T^{1|1}_\ell:=\ell^*\mathcal{T}^{1|1},\qquad T^{2|1}_\ell:=\ell^*\mathcal{T}^{2|1}
$$
for the pullback of the families of supermanifolds in Definitions~\ref{defn:supercircle} and~\ref{defn:supertorus}, respectively. We call $T^{1|1}_\ell$ an \emph{$S$-family of super circles} and $T^{2|1}_\ell$ and \emph{$S$-family of super tori}. By construction, these families have rigid conformal structures. 
\end{notation}

\begin{rmk} \label{rmk:supercircle} There is a more general notion of super circle and super torus, where $\R_{>0}$ and $\Lat$ are replaced by spaces of based lattices in $\E^{d|1}$ rather than $\E^d\subset \E^{d|1}$, for $d=1,2$ respectively. The resulting geometry is both richer and more complicated. For example, families of super circles need not carry an action of a group object $S\times \E^{d|1}/\Z^d$ over~$S$: generic $S$-families of subgroups~$\Z^d\subset \E^{d|1}$ are not be normal. One upshot when $d=2$ is that this more complicated moduli space has a natural kind of holomorphic structure; we refer to~\cite{DBEChern} for details. 
\end{rmk}

\subsection{Super Lie groupoids of rigid conformal super circles and tori}

\begin{defn} 
Let $\mathcal{M}^{d|1}$ denote the generalized super Lie groupoid whose objects are the representable sheaf $\R_{>0}$ when $d=1$ and $\Lat$ when $d=2$, and whose morphisms are the set of rigid conformal maps $T^{d|1}_\ell\to T^{d|1}_{\ell'}$ over $S$.

\end{defn} 

\begin{lem} \label{lem:grpdT}
For $d=1,2$, the generalized super Lie groupoid $\mathcal{M}^{d|1}$ is a super Lie groupoid represented by
\beq
\Mst^{1|1}\simeq\{ (\R_{>0}\times \E^{1|1}\rtimes \R^\times)/\Z\rightrightarrows \R_{>0}\}\label{eq:11grpod}
\eeq
\beq
\mathcal{M}^{2|1}\simeq \{(\E^{2|1}\rtimes \C^\times\times \SL_2(\Z)\times \Lat)/\Z^2\rightrightarrows \Lat\}.\label{eq:Mtorigrpd}
\eeq

\end{lem}

\bp
The claim amounts to showing that the morphism presheaves are representable, the target, unit, and composition are maps of supermanifolds, and the source map is a submersion. By definition, a rigid conformal isometry $f\colon T^{d|1}_\ell\to T^{d|1}_{\ell'}$ over $S$ lifts to an isometry $\tilde{f}$ on the open cover defining the rigid conformal structure on $T^{d|1}_\ell$,
\beq
\begin{tikzpicture}[baseline=(basepoint)];
\node (A) at (0,0) {$S\times \R^{d|1}$};
\node (B) at (4,0) {$S\times \R^{d|1}$};
\node (C) at (0,-1.25) {$T^{d|1}_\ell$};
\node (D) at (4,-1.25) {$T^{d|1}_{\ell'}$}; 
\draw[->,dashed] (A) to node [above] {$\tilde{f}$} (B);
\draw[->>] (A) to  (C);
\draw[->] (C) to node [above] {$f$} (D);
\draw[->>] (B) to (D);
\path (0,-.75) coordinate (basepoint);
\end{tikzpicture}\qquad \tilde{f}\in \left\{\begin{array}{ll} (\E^{1|1}\rtimes \R^\times)(S) & d=1 \\ (\E^{2|1}\rtimes \C^\times)(S) & d=2 \end{array}\right.\label{eq:tildef}
\eeq
where we identify the map $\tilde{f}$ with an $S$-point of the indicated isometry group. When $d=1$, such an $S$-point $\tilde{f}$ yields the identity map $f\colon T^{1|1}_\ell\to T^{1|1}_{\ell}$ when it is in the image of the fiberwise homomorphism
$$
S\times \Z\hookrightarrow \E\hookrightarrow  \E^{1|1}\hookrightarrow \E^{1|1}\rtimes \R^\times
$$
where the first inclusion is determined by the image of the generator $S\times \{1\}\simeq S\stackrel{\ell}\to \R_{>0}\subset \E$ determines by the map $\ell$. Hence, the supermanifold of morphisms with source $\ell\in \R_{>0}$ is the quotient of $(\E^{1|1}\times \R^\times)(S)$ by the $\Z$-action 
\beq
n\cdot (\ell,s,\eta,\dil)=(\ell,s+n\ell,\eta,\dil)\qquad \ell\in \R_{>0}(S), \ (s,\eta)\in \E^{1|1}(S), \ \dil\in \R^\times(S). \label{eq:Lambdaact1}
\eeq
In the universal case $S=\R_{>0}$ this recovers the claimed supermanifold of morphisms, using that the $\Z$-action is free and hence the quotient is indeed represented by a supermanifold. Examining the actions in the square~\eqref{eq:tildef}, the source map in~\eqref{eq:11grpod} is the projection and the target map is the composition
\beq
(\R_{>0}\times \E^{1|1}\rtimes \R^\times)/\Z\stackrel{pr}{\to} \R_{>0}\times \R^\times \stackrel{{\rm act}}{\to} \R_{>0}\label{eq:lataction1}
\eeq
where the first arrow comes from the obvious projection and the second is $(\ell,\dil)\mapsto \dil^2\ell.$ The unit and composition are determined by the unit and group structure on $\E^{1|1}\rtimes \R^\times$. 

When $d=2$, a rigid conformal map $f$ in~\eqref{eq:tildef} is more data than just $\tilde{f}$. Indeed, to give a well-defined map on the quotient, $\tilde{f}$ must be $\Z^2$-equivariant relative to a family of homomorphisms $S\times \Z^2\to S\times \Z^2$ given by $\gamma\in \SL_2(\Z)(S)$, and this homomorphism is additional data. The restriction to $\SL_2(\Z)<\GL_2(\Z)$ is so that the image of the lattice $\gamma\cdot (\ell_1, \bar\ell_1,\ell_2,\bar\ell_2)\in \Lat(S)$ remains oriented, i.e., $(\ell_1/\ell_2,\bar\ell_1/\bar\ell_2)\in \HH(S)$. Next observe that an $S$-point $(\tilde{f},\gamma)\in (\E^{2|1}\rtimes \C^\times\times \SL_2(\Z))(S)$ determines the identity map $f=\id\colon T^{2|1}_\ell\to T^{2|1}_\ell$ when it is in the image of the fiberwise homomorphism
$$
S\times \Z^2\hookrightarrow \E^2\hookrightarrow  \E^{2|1}\hookrightarrow \E^{2|1}\rtimes \C^\times\times \SL_2(\Z)
$$
where the first inclusion is determined by the image of the generators $(\ell_1,\bar\ell_1) \colon S\times \{1,0\}\hookrightarrow \C\simeq \E^2$ and $(\ell_2,\bar\ell_2) \colon S\times \{0,1\}\hookrightarrow  \C\simeq  \E^2$. Hence, the supermanifold of morphisms with source $\ell\in \Lat(S)$ is the quotient of $(\E^{2|1}\rtimes \C^\times\times \SL_2(\Z))(S)$ by the $\Z^2$-action 
\beq
&&(m,n)\cdot (w,\bar w,\eta,\dil,\bar\dil,\gamma,\ell)\mapsto (w+n\ell_1+m\ell_2,\bar w+n\bar\ell_1+m\bar\ell_2,\eta,\dil,\bar\dil,\gamma,\ell), \label{eq:Lambdaact}
\eeq
$$
(w,\bar w,\eta)\in \E^{2|1}(S), \ \gamma \in \SL_2(\Z)(S), \ (\dil,\bar\dil) \in \C^\times(S), 
$$$$ 
\ell=(\ell_1,\bar\ell_1,\ell_2,\bar\ell_2)\in \Lat(S), \ (m,n)\in \Z^2(S).
$$
In the universal case $S=\Lat$, this recovers the claimed supermanifold of morphisms, using that the $\Z^2$-action is free and hence the quotient is represented by a supermanifold. Examining the actions in the square~\eqref{eq:tildef}, the source map in~\eqref{eq:Mtorigrpd} is the projection and the target map projects along 
$(\E^{2|1}\rtimes \C^\times\times \SL_2(\Z)\times \Lat)/\Z^2\to \C^\times\times \SL_2(\Z)\times \Lat$, and applies the $\C^\times\times \SL_2(\Z)$-action on $\Lat$,
\beq
(\dil,\bar\dil,\left[\begin{array}{cc} a& b \\ c & d\end{array}\right])\cdot (\ell_1,\bar\ell_1,\ell_2,\bar\ell_2)&=&(\dil^{2}(a\ell_1+b\ell_2),\bar\dil^{2}(a\bar\ell_1+b\bar\ell_2),\label{eq:lataction}\\ 
&& \dil^{2}(c\ell_1+d\ell_2),\bar\dil^{2}(c\bar\ell_1+d\bar\ell_2))\in \Lat(S)\nonumber
\eeq
$$\quad (\dil,\bar\dil) \in \C^\times(S), \ \left[\begin{array}{cc} a& b \\ c & d\end{array}\right]\in \SL_2(\Z)(S),\quad (\ell_1,\bar\ell_1,\ell_2,\bar\ell_2)\in \Lat(S)\subset (\C\times \C)(S).
$$ 
The unit and composition are determined by the unit and group structure on $\E^{2|1}\rtimes \C^\times\times \SL_2(\Z)$. 
\ep

Let $T^{1|1}=\R^{1|1}/\Z$ and $T^{2|1}=\R^{2|1}/\Z^2$ denote the quotients for the standard subgroups $\Z\subset \E\subset \E^{1|1}$ and $\Z^2\subset \E^2\subset \E^{2|1}$, i.e., the super circle of circumference $\ell=1$, and the super torus associated to the square lattice of unit volume. There is an $\SL_2(\Z)$-action on $T^{2|1}$ inherited from the $\SL_2(\Z)$-action on $\R^{2|1}$ given by 
$$
\gamma\cdot (x,y,\theta)=(ax+by,cx+dy,\theta)\qquad (x,y,\theta)\in \R^{2|1}(S), \ \gamma=\left[\begin{array}{cc} a & b \\ c & d\end{array}\right]\in \SL_2(\Z). 
$$
We call this the action of the \emph{mapping class group} on $T^{2|1}$, and observe that it restricts the standard action of the mapping class group on $T^2\subset T^{2|1}$. The following allows one to view~$\R_{>0}$ as the moduli space of rigid conformal structures on $T^{1|1}$ and $\Lat$ as the moduli space of rigid conformal structures on $T^{2|1}$. 

\begin{lem}  There exists an isomorphism of supermanifolds over~$\R_{>0}$ and $\Lat$,
\beq
&&\R_{>0} \times T^{1|1}\stackrel{\sim}{\to} \mathcal{T}^{1|1},\qquad \Lat \times T^{2|1}\stackrel{\sim}{\to} \mathcal{T}^{2|1}, \label{eq:isotostdS}
\eeq
from the constant family with fiber $T^{d|1}$ to the family $\mathcal{T}^{d|1}$. The latter isomorphism is $\SL_2(\Z)$-equivariant for the action on $\mathcal{T}^{2|1}$ from Lemma~\ref{lem:grpdT} and the mapping class group action on~$T^{2|1}$. 
\end{lem}
\bp
Define the maps
\beq
\R_{>0}\times \R^{1|1}&\to& \R_{>0}\times \R^{1|1},\qquad \ell \in \R_{>0}(S),\ (t,\theta)\in \R^{1|1}(S)\label{eq:Zequiv}\\
(\ell,t,\theta)&\mapsto& (\ell,\ell t,\theta)\nonumber
\eeq
\beq
\Lat \times \R^{2|1}&\to& \Lat \times \R^{2|1},\qquad (\ell_1,\bar\ell_1,\ell_2,\bar\ell_2) \in \Lat (S),\  (x,y,\theta)\in \R^{2|1}(S)\label{eq:Z2equiv}\\
(\ell_1,\bar\ell_1,\ell_2,\bar\ell_2,x,y,\theta)&\mapsto& (\ell_1,\bar\ell_1,\ell_2,\bar\ell_2,\ell_1 x+\ell_2 y,x\bar \ell_1+y\bar \ell_2,\theta),\nonumber
\eeq
where in~\eqref{eq:Z2equiv} an $S$-point of the source is specified using~\eqref{eq:r211} whereas $S$-points of the target are given in the description~\eqref{eq:r212}. Observe that~\eqref{eq:Zequiv} is $\Z$-equivariant for the actions on source and target,
\beq
n\cdot (\ell,t,\theta)=(\ell,t+n,\theta), && n\cdot (\ell,t,\theta)=(\ell,t+ n \ell,\theta)\nonumber
\eeq 
and~\eqref{eq:Z2equiv} is $\Z^2$-equivariant for the action on source and target,
\beq
(n,m)\cdot (\ell_1,\bar\ell_1,\ell_2,\bar\ell_2,x,y,\theta)&=&(\ell_1,\bar\ell_1,\ell_2,\bar\ell_2,x+n,y+m,\theta) \nonumber\\
(n,m)\cdot (\ell_1,\bar\ell_1,\ell_2,\bar\ell_2,z,\bar z,\theta)&=&(\ell_1,\bar\ell_1,\ell_2,\bar\ell_2,z+n\ell_1+m\ell_2,\bar z+n\bar\ell_1+m\bar\ell_2,\theta).\nonumber
\eeq
Hence~\eqref{eq:Zequiv} and~\eqref{eq:Z2equiv} determine maps between the respective quotients, defining maps~\eqref{eq:isotostdS}. These are easily seen to be isomorphisms. Checking that the isomorphism $\Lat\times T^{2|1}\to \mathcal{T}^{2|1}$ is $\SL_2(\Z)$-equivariant for the previously defined actions is an easy computation. 
\ep

\begin{cor} \label{cor:isowstd}
For any $S$-family of super circles or super tori, there is a canonical isomorphism $T^{1|1}_\ell\simeq T^{1|1}\times S$ and $T^{2|1}_\ell\simeq T^{2|1}\times S$ with the standard family. 
\end{cor} 
\bp
This follows immediately from pulling back the isomorphisms~\eqref{eq:isotostdS} along a map $\ell\colon S\to \R_{>0}$ or $\ell\colon S\to \Lat$. 
\ep

\subsection{$G$-bundles on super circles and super tori}

Using the canonical open cover $\R^{d|1}\to T^{d|1}$, a $G$-bundle on $T^{d|1}$ can be identified with a functor $\R^{d|1}\sq \Z^d\to \pt\sq G$, which in turn is given by a homomorphism $\Z^d\to G$ since the functor on objects $\R^{d|1}\to \pt$ is uniquely determined. Similarly, a $G$-bundle over a family of super conformal circles or tori is classified by an $S$-point of 
\beq
\R_{>0}\times \Map(\R^{1|1}\sq \Z,\pt\sq G)&\simeq& \R_{>0}\times \Map(\pt\sq \Z,\pt\sq G)\simeq \R_{>0}\times G,\label{eq:Gbundleobj1}\\
 \Lat\times \Map(\R^{2|1}\sq \Z^2,\pt\sq G)&\simeq& \Lat\times \Map(\pt\sq \Z^2,\pt\sq G)\simeq \Lat\times G^{(2)}. \label{eq:Gbundleobj2}
\eeq
where we deduce representability from the descriptions on the right. Indeed, an $S$-point $(\ell\in \R_{>0}(S),g\in G(S))$ of~\eqref{eq:Gbundleobj1} or $(\ell\in \Lat(S),(g_1,g_2)\in G^{(2)}(S))$ of~\eqref{eq:Gbundleobj2} determines 
\beq
p\colon P_g\to T^{1|1}\times S\stackrel{\sim}{\to} T^{1|1}_\ell,\qquad p\colon P_{g_1,g_2} \to T^{2|1}\times S\stackrel{\sim}{\to} T^{2|1}_\ell,\label{eq:Gbundlesrct}
\eeq
respectively, where the isomorphisms in~\eqref{eq:Gbundlesrct} come from Corollary~\ref{cor:isowstd}, and 
\beq
P_g:=(G\times \R^{1|1}\times S)/\Z\qquad P_{g_1,g_2}:=(G\times\R^{2|1}\times S)/\Z^2\label{eq:Gbundledefnagain}
\eeq
where the left $\Z$-action defining $P_g$ is the standard action on $\R^{1|1}$ and the action on $G$ generated by $g$, and similarly the left $\Z^2$-action defining $P_{g_1,g_2}$ is the standard action on~$\R^{2|1}$ and the action on $G$ generated by $g_1,g_2\in G$. 

\begin{defn}\label{defn:Gbundlesct}
Define generalized super Lie groupoids $\Bun_G(\mathcal{T}^{1|1})$ and $\Bun_G(\mathcal{T}^{2|1})$ whose presheaves of objects are the (representable) presheaves~\eqref{eq:Gbundleobj1} and~\eqref{eq:Gbundleobj2}, respectively, and whose presheaves of morphisms have as $S$-points commuting squares
\beq\begin{tikzpicture}[baseline=(basepoint)];
\node (A) at (0,0) {$P_{g}$};
\node (B) at (4,0) {$P_{g'}$};
\node (C) at (0,-1.25) {$T^{1|1}_\ell$};
\node (D) at (4,-1.25) {$T^{1|1}_{\ell'}$}; 
\draw[->] (A) to (B);
\draw[->] (A) to  (C);
\draw[->] (C) to (D);
\draw[->] (B) to (D);
\path (0,-.75) coordinate (basepoint);
\end{tikzpicture}
\quad {\rm and} \quad 
\begin{tikzpicture}[baseline=(basepoint)];
\node (A) at (0,0) {$P_{g_1,g_2}$};
\node (B) at (4,0) {$P_{g_1',g_2'}$};
\node (C) at (0,-1.25) {$T^{2|1}_\ell$};
\node (D) at (4,-1.25) {$T^{2|1}_{\ell'}$}; 
\draw[->] (A) to (B);
\draw[->] (A) to  (C);
\draw[->] (C) to (D);
\draw[->] (B) to (D);
\path (0,-.75) coordinate (basepoint);
\end{tikzpicture}\label{eq:Gbundleisos}
\eeq
respectively, given by an isomorphism of $G$-bundles covering a rigid conformal map. The source and target are respectively given by taking the left and right columns of the above squares. The unit map is the square whose horizontal arrows are identities. Composition comes from composing squares horizontally. 
\end{defn}

\begin{lem} The generalized super Lie groupoids $\Bun_G(\mathcal{T}^{1|1})$ and $\Bun_G(\mathcal{T}^{2|1})$ are super Lie groupoids, i.e., their presheaves of objects and morphisms are representable, and the source map is a submersion.
\end{lem}
\bp
The presheaves of objects~\eqref{eq:Gbundleobj1} and~\eqref{eq:Gbundleobj2} are representable by definition. To show that morphisms are representable, we pullback along $\R^{d|1}\times S\to T^{d|1}_\ell$ so that the $G$-bundles trivialize
\beq
\begin{tikzpicture}[baseline=(basepoint)];
\node (A) at (0,0) {$G\times\R^{1|1}\times S$};
\node (B) at (3,0) {$P_{g}$};
\node (C) at (0,-1.25) {$\R^{1|1}\times S$};
\node (D) at (3,-1.25) {$T^{1|1}_{\ell'}$}; 
\draw[->>] (A) to  (B);
\draw[->] (A) to node [left] {$\tilde{p}$} (C);
\draw[->>] (C) to  (D);
\draw[->] (B) to node [right] {$p$} (D);
\path (0,-.75) coordinate (basepoint);
\end{tikzpicture}
\qquad\qquad \begin{tikzpicture}[baseline=(basepoint)];
\node (A) at (0,0) {$G\times\R^{2|1}\times S$};
\node (B) at (3,0) {$P_{g_1,g_2}$};
\node (C) at (0,-1.25) {$\R^{2|1}\times S$};
\node (D) at (3,-1.25) {$T^{2|1}_{\ell'}.$}; 
\draw[->] (A) to  (B);
\draw[->] (A) to node [left] {$\tilde{p}$} (C);
\draw[->>] (C) to (D);
\draw[->] (B) to node [right] {$p$}  (D);
\path (0,-.75) coordinate (basepoint);
\end{tikzpicture}
\eeq
Then isomorphisms~\eqref{eq:Gbundleisos} can be described analogously to~\eqref{eq:tildef}, namely as isomorphisms of trivial $G$-bundles 
$$
\begin{tikzpicture}[baseline=(basepoint)];
\node (A) at (0,0) {$G\times\R^{1|1}\times S$};
\node (B) at (4,0) {$G\times\R^{1|1}\times S$};
\node (C) at (0,-1.25) {$P_{g}$};
\node (D) at (4,-1.25) {$P_{g'}$}; 
\draw[->] (A) to node [above] {$(h,\tilde{f})$} (B);
\draw[->>] (A) to  (C);
\draw[->] (C) to (D);
\draw[->>] (B) to (D);
\path (0,-.75) coordinate (basepoint);
\end{tikzpicture}\qquad \begin{tikzpicture}[baseline=(basepoint)];
\node (A) at (0,0) {$G\times\R^{2|1}\times S$};
\node (B) at (4,0) {$G\times\R^{2|1}\times S$};
\node (C) at (0,-1.25) {$P_{g_1,g_2}$};
\node (D) at (4,-1.25) {$P_{g_1',g_2'}$}; 
\draw[->] (A) to node [above] {$(h,\tilde{f})$} (B);
\draw[->>] (A) to  (C);
\draw[->] (C) to (D);
\draw[->>] (B) to (D);
\path (0,-.75) coordinate (basepoint);
\end{tikzpicture}
$$
determined by data $h\in G(S)$ and $\tilde{f}\in (\E^{1|1}\rtimes \R^\times)(S)$ when $d=1$ and $\tilde{f}\in (\E^{2|1}\rtimes \C^\times)(S)$ when $d=2$. The upper horizontal arrow in the diagrams above is then given by left multiplication by $h$ on $G\times S$ and the left action of $\tilde{f}$ on $\R^{d|1}\times S$ by a rigid conformal isometry. When $d=2$, we have the additional data of $\gamma\in \SL_2(\Z)(S)$ acting on $\ell$ and $(g_1,g_2)$. Hence, $g'=hgh^{-1}$ when $d=1$, and $(g_1,'g_2')$ is given by~\eqref{eq:G2act}. The data $(h,\tilde{f})$ determine the identity map in~\eqref{eq:Gbundleisos} if they are in the image of the subgroup 
\beq
S\times \Z&\stackrel{\langle g,\ell\rangle}{\hookrightarrow}& G\times \E\hookrightarrow G\times \E^{1|1}\rtimes \R^\times\label{eq:Zsubgroup}\\
S\times \Z^2&\stackrel{\langle (g_1,g_2),\ell\rangle}{\hookrightarrow}& G\times \E^2\hookrightarrow G\times \E^{2|1}\rtimes \C^\times. \label{eq:Z2subgroup}
\eeq
Therefore, the presheaves of isomorphisms whose $S$-points are the squares~\eqref{eq:Gbundleisos} are 
$$
(G\times \E^{1|1}\rtimes \R^\times\times \R_{>0}\times G)/\Z \quad {\rm and}\quad (G\times \E^{2|1}\rtimes \C^\times\times\SL_2(\Z)\times \Lat\times G^{(2)})/\Z^2
$$
respectively, where the quotient is by the $S=\R_{>0}\times G$ family of subgroups~\eqref{eq:Zsubgroup} and the $S=\Lat\times G^{(2)}$ family of subgroups~\eqref{eq:Z2subgroup}. These are quotients of a supermanifold by a free action, and hence are representable presheaves. The source map in both cases is determined by the projection map, which is a submersion. The proposition is proved. 
\ep

\section{Fields and inertia fields as groupoids}\label{sec:cocycles}

\subsection{Motivation from physics: Twist fields}\label{sec:twistfields}

Recall Definition~\ref{defn:map} for mapping objects between supermanifolds. Let $M$ be a Riemannian manifold. Classical mechanics with $N=1$ supersymmetry studies the super path space $\Map(\R^{1|1},M)$. The classical sigma model with $\mathcal{N}=(0,1)$ supersymmetry studies the mapping space $\Map(\R^{2|1},M)$. To motivate the definitions below, we briefly review some standard constructions in the context of \emph{compactification} of these classical theories. We focus on the case $d=2$; $d=1$ is similar. 

Compactification first restricts to a subspace of fields with periodic boundary conditions
\beq
\phi(z+n\ell_1+m\ell_2,\bar z+n\bar\ell_1+m\bar\ell_2,\theta)=\phi(z,\bar z,\theta),\quad \phi\colon S\times \R^{2|1}\to M,\label{eq:periodic}
\eeq
for a based lattice $(\ell_1,\bar\ell_1,\ell_2,\bar\ell_2)\in \Lat(S)$ determining a $\Z^2$-action on $S\times \R^{2|1}$. The terminology comes from the fact that maps satisfying~\eqref{eq:periodic} are equivalent to maps $T^{2|1}_\ell\to M$ where the source is \emph{compact} (fiberwise over $S$). Indeed, fields satisfying~\eqref{eq:periodic} can be identified with the generalized supermanifold~$\Lat\times \Map(T^{2|1},M)$ where $\ell\in \Lat(S)$ and $\phi\in \Map(T^{2|1},M)(S)$ determines 
$$
T^{2|1}_\ell \stackrel{\sim}{\leftarrow} S\times T^{2|1}\stackrel{\phi}{\to} M
$$ 
with the isomorphism on the left coming from from Corollary~\ref{cor:isowstd}. 
Next, \emph{1-loop quantization} looks to integrate along the fibers $\Lat\times \Map(T^{2|1},M)\to \Lat$. These fibers are infinite-dimensional, so a priori integration is an ill-defined operation. However, \emph{supersymmetric localization} formally identifies the integral with one over a finite-dimensional subspace of maps~\eqref{eq:periodic} invariant under the translation action of $\E^2\simeq \C$,
$$
(w,\bar w)\cdot \phi(z,\bar z,\theta)=\phi(z+w,\bar z+\bar w,\theta),\qquad (w,\bar w)\in \C(S)\simeq \E^2(S).
$$
This $\E^2$-fixed subspace is $\Lat\times \Map(\R^{0|1},M)\subset \Lat\times \Map(T^{2|1},M)$. The isomorphism $C^\infty(\Lat\times \Map(\R^{0|1},M))\simeq \Omega^\bullet(M;C^\infty(\Lat))$ from Example~\ref{ex:PiT} reduces the path integral to an integral of a differential form on a finite-dimensional manifold $M$ with coefficients in~$C^\infty(\Lat)$. We refer to~\cite{BElocalization} for a construction of this path integral via localization. 

When a finite group $G$ acts isometrically on $M$, compactification can be implemented relative to the \emph{twisted} boundary conditions
\beq
\phi(z+n\ell_1+m\ell_2,\bar z+n\bar\ell_1+m\bar\ell_2,\theta)=g_1^ng_2^m\phi(z,\bar z,\theta),\label{eq:twistedperiodic}
\eeq
where $g_1,g_2\in G$ are commuting elements; e.g., see~\cite[Equations~2 and~3]{Vafatorsion}. The set of all maps satisfying~\eqref{eq:twistedperiodic} for a fixed $(g_1,g_2)$ is called a \emph{twisted sector}. The \emph{twist fields} are then the union of twisted sectors, 
\beq
\coprod_{(g_1,g_2)\in G^{(2)}} \Lat\times \Map(\R^{2|1},M)^{\Z^2}\simeq \Lat\times \Map(\R^{2|1}\sq \Z^2,M\sq G)\label{eq:twistfields}
\eeq
which can be re-expressed as the generalized supermanifold of functors on the right hand side above, where we identify a functor with a map satisfying~\eqref{eq:twistedperiodic}. 
The 1-loop quantization of twist fields again considers a localized path integral, but in this case the $\E^2$-fixed subspace onto which the path integral localizes is
\beq
\Lat\times \Map(\R^{0|1}\sq \Z^2,M\sq G)\subset \Lat\times \Map(\R^{2|1}\sq \Z^2,M\sq G), \label{eq:inertia}
\eeq
using that the $\E^2$-action on the groupoid $\R^{2|1}\sq \Z^2$ is free with quotient $\R^{0|1}\sq \Z^2$. We recognize this subspace as a generalization of the double inertia groupoid of $M\sq G$,
$$
\Lat\times \Map(\R^{0|1}\sq \Z^2,M\sq G)\simeq \Lat\times \Map(\R^{0|1},\Map(\pt\sq \Z^2,M\sq G)).
$$
Twist fields enhance the double inertia groupoid by adding a conformal structure to the torus and a nilpotent direction from $\R^{0|1}$. By Example~\ref{ex:PiT}, functions on twist fields are differential forms on the double inertia groupoid with coefficients in $C^\infty(\Lat)$, 
\beq
C^\infty(\Lat\times \Map(\R^{0|1},\Map(\pt\sq \Z^2,M\sq G)))&\simeq& \Omega^\bullet(\Map(\pt\sq \Z^2,M\sq G);C^\infty(\Lat))\nonumber \\
&\simeq& \prod_{(g_1,g_2)\in G^{(2)}} \Omega^\bullet(M^{g_1,g_2};C^\infty(\Lat)),\nonumber
\eeq
as we shall explain in greater detail below. By Proposition~\ref{prop:dgs}, this algebra of functions is closely related to the $G$-equivariant elliptic cohomology of $M$. However, it is missing actions by~$G\times \SL_2(\Z)$ and the de~Rham differential. From Example~\ref{ex:deRhamsheaf}, we might hope to recover these by an action of a super Lie group. It turns out that the symmetries of the supersymmetric sigma model lead to precisely such an action. 

To describe these symmetries it is useful to change points of view. Given a functor $\Phi\colon S\times \R^{2|1}\sq \Z^2\to M\sq G$, consider 
\beq
&&T^{2|1}\times S\simeq \R^{2|1}/\Z^2\times S \leftarrow P_{g_1,g_2} :=(G\times \R^{2|1}\times S)/\Z^2 \stackrel{\id_G\times \Phi_0}{\to} G\times M\stackrel{{\rm act}}{\to} M. \label{eq:bundleize}
\eeq
where $(g_1,g_2)\in G^{(2)}(S)$ comes from the value of $\Phi$ on morphisms. If we also specify $\ell\in \Lat(S)$,~\eqref{eq:bundleize} determines a span $T^{2|1}_\ell\leftarrow P_{g_1,g_2}\stackrel{\phi}{\to} M$
where $P_{g_1,g_2}\to T^{2|1}$ is a principal $G$-bundle $\phi$ is a $G$-equivariant map. In this description, symmetries (i.e., isomorphisms between twist fields) are specified by a rigid conformal map $T^{2|1}_\ell\to T^{2|1}_{\ell}$ covered by an isomorphism of $G$-bundles $P_{g_1,g_2}\to P_{g_1',g_2'}$ that is compatible with the maps~$\phi$ and~$\phi'$ to~$M$. This leads to our definition of fields below. 

\subsection{Definitions of fields and inertia fields}

Below we use the boldface notation ${\bf g}\in G^{(d)}$ to denote an element $g\in G$ when $d=1$, and a pair of commuting elements $g_1,g_2\in G$ when $d=2$. 

\begin{defn}\label{defn:fields} Define a generalized super Lie groupoid $\mathcal{L}^{d|1}(M\sq G)$ of \emph{fields} whose objects are the presheaves with $S$-points given by triples $(\ell,{\bf g},\phi)$ determining
\beq
T^{d|1}_\ell\stackrel{p}{\leftarrow} P_{{\bf g}}\stackrel{\phi}{\to} M 
\eeq
where $P_{\bf g}$ and $p$ are from~\eqref{eq:Gbundlesrct} and $\phi$ is a $G$-equivariant map.
The presheaf of morphisms has $S$-points given by pairs $(h,f)$ that sit in a commutative diagram
\beq
&&\begin{tikzpicture}[baseline=(basepoint)];
\node (A) at (0,0) {$T^{d|1}_\ell$};
\node (B) at (3,0) {$P_{\bf g}$};
\node (BB) at (3,-1.5) {$P_{\bf g'}$};
\node (C) at (0,-1.5) {$T^{d|1}_{\ell'}$};
\node (D) at (6,-.75) {$M$}; 
\draw[->] (B) to (A);
\draw[->] (A) to node [left] {$f$} (C);
\draw[->] (BB) to  (C);
\draw[->,bend left=10] (B) to node [above] {$\phi$} (D);
\draw[->,bend right=10] (BB) to node [below] {$\phi'$} (D);
\draw[->] (B) to node [left] {$h$} (BB);
\path (0,-.75) coordinate (basepoint);
\end{tikzpicture}\label{eq:isooffields}
\eeq
where $h$ is an isomorphism of $G$-bundles covering a rigid conformal isometry~$f$. Composition of morphisms in $\mathcal{L}^{d|1}(M\sq G)$ is determined by $(h,f)\circ (h',f')=(h\circ h',f\circ f')$, and the unit map is determined by the diagram~\eqref{eq:isooffields} with $(f,h)=(\id,\id)$. 
\end{defn}

\begin{prop} \label{prop:natural}
A functor $N\sq H\to M\sq G$ between Lie groupoids determines a functor $\mathcal{L}^{d|1}(N\sq H)\to \mathcal{L}^{d|1}(M\sq G)$ of generalized super Lie groupoids. 
\end{prop}
\bp
A functor between action Lie groupoids is equivalent to an equivariant map $N\to M$ relative to a homomorphism $\zeta\colon H\to G$. The composition $P_{\bf g} \stackrel{\phi}{\to} N \to M$ is $H$-equivariant, and so determines a $G$-equivariant map~$\phi_\zeta$ from the associated $G$-bundle,
$$
\phi_\zeta\colon  P_{\zeta({\bf g})}\simeq P_{\bf g}\times_\zeta G\to M
 $$ 
Define a functor $\mathcal{L}^{d|1}(N\sq H)\to \mathcal{L}^{d|1}(M\sq G)$ whose map on objects is the assignment $(\ell,{\bf g},\phi)\mapsto (\ell,\zeta({\bf g}),\phi_\zeta)$. On morphisms the functor is defined as follows. Given an isomorphism $(f,h)$ as in~\eqref{eq:isooffields}, we observe that naturality of the associated bundle construction gives a commutative diagram
\beq
&&\begin{tikzpicture}[baseline=(basepoint)];
\node (A) at (0,0) {$T^{d|1}_\ell$};
\node (B) at (3,0) {$P_{\zeta({\bf g})}\simeq P_{\bf g}\times_\zeta G$};
\node (BB) at (3,-1.5) {$P_{\zeta({\bf g'})}\simeq P_{\bf g'}\times_\zeta G$};
\node (C) at (0,-1.5) {$T^{d|1}_{\ell'}$};
\node (D) at (7,-.75) {$M$}; 
\draw[->] (B) to (A);
\draw[->] (A) to node [left] {$f$} (C);
\draw[->] (BB) to  (C);
\draw[->,bend left=10] (B) to node [above] {$\phi_\zeta$} (D);
\draw[->,bend right=10] (BB) to node [below] {$\phi_\zeta'$} (D);
\draw[->] (B) to node [left] {$\zeta(h)$} (BB);
\path (0,-.75) coordinate (basepoint);
\end{tikzpicture}\label{diag:assocmap}
\eeq
for the map of $G$-bundles $\zeta(h)$ associated to $h$ via $\zeta$. Hence, the assignment $(f,h)\mapsto (f,\zeta(h))$ defines the functor $\mathcal{L}^{d|1}(N\sq H)\to \mathcal{L}^{d|1}(M\sq G)$ on morphisms. It is easy to check that the required compatibility conditions are satisfied, so that these assignments do indeed determine a functor.
\ep

\begin{rmk}
We observe that there is an isomorphism of generalized super Lie groupoids $\Bun_G(\mathcal{T}^{d|1})\simeq \mathcal{L}^{d|1}(\pt\sq G)$, and hence by the previous lemma there are canonical functors
$$
\mathcal{L}^{d|1}(M\sq G)\to \Bun_G(\mathcal{T}^{d|1})
$$
induced by the $G$-equivariant map $M\to \pt$. 
\end{rmk}
The canonical covers $\R^{d|1}\times S\to T^{d|1}$ and $G\times \R^{d|1}\times S\to P_{\bf g}$ yield the commutative diagram 
\beq
&&\begin{tikzpicture}[baseline=(basepoint)];
\node (A) at (0,0) {$\R^{d|1}\times S$};
\node (B) at (3,0) {$G\times \R^{d|1}\times S$};
\node (BB) at (3,-1.5) {$P_{\bf g}$};
\node (C) at (0,-1.5) {$T^{d|1}_{\ell}$};
\node (D) at (6,-1.5) {$M.$}; 
\draw[->] (B) to (A);
\draw[->>] (A) to (C);
\draw[->] (BB) to  (C);
\draw[->] (B) to node [above] {$\tilde{\phi}$} (D);
\draw[->] (BB) to node [below] {$\phi$} (D);
\draw[->>] (B) to (BB);
\path (0,-.75) coordinate (basepoint);
\end{tikzpicture}\label{eq:Eactiondiag}
\eeq
There is an $\E^d(S)$-action on triples $(\ell,{\bf g},\tilde{\phi})$ that is trivial on $\ell$ and ${\bf g}$, and precomposes $\tilde{\phi}$ with the $\E^d(S)$-action on $\R^{d|1}\times S$ and the trivial action on $G$. 

\begin{lem}
The $\E^d(S)$-action on $(\ell,{\bf g},\tilde{\phi})$ determines an $\E^d$-action on the objects of~$\mathcal{L}^{d|1}(M\sq G)$. 
\end{lem}
\bp
Since the $\E^d(S)$-action commutes with the $\Z^d(S)$-action determined by $\ell$, it yields a well-defined action on the $\Z^d$-quotient, $(\ell,{\bf g},\phi)$, in~\eqref{eq:Eactiondiag}. It is straightforward to verify that this action is natural with respect to base changes $S'\to S$, yielding the claimed action. 
\ep

\begin{defn}
For $d=1,2$, define generalized super Lie groupoids $\mathcal{L}_0^{d|1}(M\sq G)$ of \emph{inertia fields} as the full subgroupoids of $\mathcal{L}^{d|1}(M\sq G)$ with objects the sub-presheaves
\beq
&&{\rm Obj}(\mathcal{L}_0^{1|1}(M\sq G))\simeq {\rm Obj}(\mathcal{L}^{1|1}(M\sq G))^\E,\qquad 
{\rm Obj}(\mathcal{L}_0^{2|1}(M\sq G))\simeq {\rm Obj}(\mathcal{L}^{2|1}(M\sq G))^{\E^2}\nonumber
\eeq
fixed by the $\E^d$-action. 
\end{defn}


\begin{cor}\label{cor:natural}
A functor $N\sq H\to M\sq G$ between Lie groupoids determines a functor $\mathcal{L}_0^{d|1}(N\sq H)\to \mathcal{L}_0^{d|1}(M\sq G)$ of generalized Lie groupoids.
\end{cor}
\bp
One applies the same construction as in Lemma~\ref{prop:natural}, observing that if $P_{\bf g}\to N$ is $\E^2$-invariant, then so is the map out of the associated bundle $P_{\zeta({\bf g})}\to M$.
\ep

\begin{rmk}
Proposition~\ref{prop:natural} and Corollary~\ref{cor:natural} can also be deduced by composing functors $S\times \R^{d|1}\sq \Z^2\to N\sq H\to M\sq G$. We have chosen the proof involving associated bundles because it will make for more explicit indexing sets in the finite path integrals in~\S\ref{sec:discretetorsion}. 
\end{rmk}

\subsection{Inertia fields are a super Lie groupoid}
\begin{prop} \label{prop:inertiapresent}
For $d=1,2$, the generalized super Lie groupoids $\mathcal{L}_0^{d|1}(M\sq G)$ are super Lie groupoids: the object and morphism presheaves are represented by supermanifolds, the structure maps are maps of supermanifolds, and the source map is a submersion. 
\end{prop}

\bp 
We observe a map $\phi\colon P_{\bf g}\to M$ being $\E^d$-invariant is equivalent to the existence of the dashed arrow
\beq
&&\begin{tikzpicture}[baseline=(basepoint)];
\node (A) at (0,0) {$G\times \R^{d|1}\times S$};
\node (B) at (4,0) {$G\times \R^{0|1}\times S$};
\node (C) at (0,-1.5) {$P_{\bf g}$};
\node (D) at (4,-1.5) {$M$};
\draw[->] (A) to node [above] {$q$} (B);
\draw[->] (A) to (C);
\draw[->,dashed] (B) to (D);
\draw[->] (C) to node [above] {$\phi$} (D);
\path (0,-.75) coordinate (basepoint);
\end{tikzpicture}\nonumber
\eeq
where $q$ is determined by the quotient map $\R^{d|1}\to \R^{d|1}/\E^d\simeq \R^{0|1}$ for the (free) $\E^d$-action. Restricting the dashed arrow to $\{e\}\times \R^{0|1}\times S$ together with the data ${\bf g}$ determines a functor
$$
(\R^{0|1}\sq \Z^d)\times S\to M\sq G,
$$
for the trivial $\Z^d$-action on $\R^{0|1}$. Conversely such a functor determines a diagram~\eqref{eq:Eactiondiag} and hence an $S$-point of the object presheaf. This provides isomorphisms of presheaves, 
\beq
{\rm Obj}(\mathcal{L}_0^{1|1}(M\sq G))&\simeq&\R_{>0}\times \Map(\R^{0|1}\sq \Z,M\sq G)\label{eq:inertia11}\\
&\simeq& \R_{>0}\times \Map(\R^{0|1},\Map(\pt\sq \Z,M\sq G))\nonumber\\
&\simeq& \R_{>0}\times \coprod_{g\in G} \Map(\R^{0|1},M^g)\simeq  \coprod_{g\in G}\R_{>0}\times \Map(\R^{0|1},M^g)\nonumber\\
{\rm Obj}(\mathcal{L}_0^{2|1}(M\sq G))&\simeq&\Lat\times \Map(\R^{0|1}\sq \Z^2,M\sq G)\label{eq:inertia21}\\
&\simeq&\Lat\times \Map(\R^{0|1},\Map(\pt\sq \Z^2,M\sq G))\nonumber \\
&\simeq& \Lat \times \coprod_{(g_1,g_2)\in G^{(2)}} \Map(\R^{0|1},M^{g_1,g_2})\nonumber\\
&\simeq& \coprod_{(g_1,g_2)\in G^{(2)}}  \Lat \times \Map(\R^{0|1},M^{g_1,g_2}).\nonumber
\eeq
demonstrating that the object sheaves are indeed representable. 

Next suppose we are given an $S$-point of the morphism presheaf~$(f,h)$ with source $(\ell,{\bf g},\phi)\in {\rm Obj}(\mathcal{L}_0^{d|1}(M\sq G))(S)$. The pull back of $P_{\bf g}$ to the cover $S\times \R^{d|1}$ trivializes, and the pullback of $h$ is determined by maps $\tilde{f}\colon S\times \R^{d|1}\to S\times \R^{d|1}$ and $\tilde{h}\colon S\times G\to S\times G$
\beq
\begin{tikzpicture}[baseline=(basepoint)];
\node (A) at (0,0) {$G\times \R^{d|1}\times S$};
\node (B) at (4,0) {$G\times \R^{d|1}\times S$};
\node (C) at (0,-1.25) {$P_{\bf g}$};
\node (D) at (4,-1.25) {$P_{\bf g'}$}; 
\draw[->] (A) to node [above] {$\tilde{h}\times \tilde{f}$} (B);
\draw[->>] (A) to  (C);
\draw[->] (C) to node [above] {$h$} (D);
\draw[->>] (B) to (D);
\path (0,-.75) coordinate (basepoint);
\end{tikzpicture}\label{eq:supertorusmap}
\eeq
where we identify the map $\tilde{h}$ with left multiplication by $\tilde{h}\in G(S)$ and $\tilde{f}$ with the left action by $\tilde{f}\in (\E^{1|1}\rtimes \R^\times)(S)$ when $d=1$ and $\tilde{f}\in (\E^{2|1}\rtimes \C^\times)(S)$ when $d=2$. The data $\tilde{h}\times \tilde{f}$ determines the identity map $(h,f)=(\id,\id)$ on $(\ell,{\bf g},\phi)$ if $(\tilde{h},\tilde{f})$ are in the image of
$$
S\times \Z\stackrel{\langle \ell,g\rangle }{\hookrightarrow} (\E^{1|1}\times \R^\times \times G),\qquad S\times \Z^2\stackrel{\langle \ell,(g_1,g_2)\rangle}{\hookrightarrow} (\SL_2(\Z) \times \E^{2|1}\times \C^\times \times G).
$$
This shows that the presheaf of morphisms for $\mathcal{L}_0^{d|1}(M\sq G)$ is given by the $\Z^d$-quotients,
\beq
{\rm Mor}(\mathcal{L}_0^{1|1}(M\sq G))&=&\coprod_{g\in G} (\E^{1|1}\rtimes \R^\times\times G \times \R_{>0}\times \Map(\R^{0|1},M^g))/\Z\nonumber \\
{\rm Mor}(\mathcal{L}_0^{2|1}(M\sq G))&=&\coprod_{(g_1,g_2)\in G^{(2)}} (\SL_2(\Z) \times \E^{2|1}\rtimes \C^\times\times G \times \Lat \times \Map(\R^{0|1},M^{g_1,g_2}))/\Z^2.\nonumber
\eeq
These are free $\Z^d$-actions on representable presheaves, so the quotients are also representable. The structure maps of the generalized Lie groupoid $\mathcal{L}_0^{d|1}(M\sq G)$ are therefore morphisms between representable presheaves, and hence determine maps of supermanifolds. The source maps come from the obvious projections, which are submersions. This completes the proof.
\ep

\begin{rmk}\label{rmk:fail}
Proposition~\ref{prop:inertiapresent} fails spectacularly for $G$ compact Lie but not finite. The reason is that ${\bf g}\in G^{(d)}$ can vary in continuous families while the fixed point sets~$M^{\bf g}$ do not, e.g., the dimension of a fixed point set can change under a continuous deformation. This results in the presheaf of objects being singular, and hence non-representable. 
\end{rmk}

For computations below, we will need a more explicit description of the target map in the super Lie groupoids $\mathcal{L}_0^{d|1}(M\sq G)$. 

\begin{lem}\label{lem:grpdtarget} The target map in the super Lie groupoid $\mathcal{L}_0^{d|1}(M\sq G)$ is a composition, 
\beq\resizebox{1.1\textwidth}{!}{$\begin{array}{rll}
\coprod (\E^{1|1}\rtimes \R^\times\times G \times \R_{>0}\times \Map(\R^{0|1},M^g))/\Z &\stackrel{p}{\to}& \coprod \E^{0|1}\rtimes \R^\times\times G \times \R_{>0}\times \Map(\R^{0|1},M^g)\nonumber\\
&\stackrel{\mu}{\to}& \coprod \R_{>0}\times \Map(\R^{0|1},M^g)\nonumber\\
\coprod (\SL_2(\Z) \times \E^{2|1}\rtimes \C^\times\times G \times \Lat \times \Map(\R^{0|1},M^{g_1,g_2}))/\Z^2&\stackrel{p}{\to}& \coprod \SL_2(\Z) \times \E^{0|1}\rtimes \C^\times\times G \times \Lat \times \Map(\R^{0|1},M^{g_1,g_2})\nonumber\\
&\stackrel{\mu}{\to}&  \coprod  \Lat \times \Map(\R^{0|1},M^{g_1,g_2}).\end{array}$}\nonumber
\eeq
when $d=1$ and $2$, respectively, where in both cases $p$ is the projection determined by $\E^{d|1}\to \E^{d|1}/\E^d\simeq \E^{0|1}$, and $\mu$ is an action map. The action is diagonal for the actions on $\R_{>0}$ and $\Lat$ determined by~\eqref{eq:lataction1} and~\eqref{eq:lataction}, respectively. The $\E^{0|1}\rtimes \R^\times\times G$-action on $\coprod \Map(\R^{0|1},M^g)$ is determined by the left action $h\colon  M^g\to M^{hgh^{-1}}$ and the $\E^{0|1}\times \R^\times$-action on $\Map(\R^{0|1},M^g)$ by precomposition. Similarly the $\E^{0|1}\rtimes \C^\times \times G$-action on $\coprod \Map(\R^{0|1},M^{g_1,g_2})$ is determined by the left action $h\colon  M^{g_1,g_2}\to M^{hg_1h^{-1},hg_2h^{-1}}$ and the $\E^{0|1}\times \C^\times$-action on $\Map(\R^{0|1},M^{g_1,g_2})$ by precomposition. 
\end{lem}
\bp
In terms of the description of objects from~\eqref{eq:inertia11} and~\eqref{eq:inertia21}, suppose we are given $S$-points of objects determined by functors $\Phi,\Phi'\colon S\times \R^{0|1}\sq \Z^d\to M\sq G$ and $\ell\in \R_{>0}(S)$ or $\ell\in \Lat(S)$ determining lattices. Then an isomorphism $(\ell,\Phi)\to (\ell,\Phi')$ between $S$-points of the super Lie groupoid $\mathcal{L}_0^{d|1}(M\sq G)$ gives a diagram of groupoids
\beq
&&\begin{tikzpicture}[baseline=(basepoint)];
\node (A) at (-2,0) {$T^{d|1}_\ell$};
\node (AAA) at (1,0) {$S\times \R^{d|1}\sq \Z^d$};
\node (BBB) at (1,-1.5) {$S\times \R^{d|1}\sq \Z^d$};
\node (B) at (4,0) {$S\times \R^{0|1}\sq \Z^d$};
\node (BB) at (4,-1.5) {$S\times \R^{0|1}\sq \Z^d$};
\node (C) at (-2,-1.5) {$T^{d|1}_{\ell'}$};
\node (D) at (8,-.75) {$M\sq G$}; 
\node (E) at (6.2,-.75) {$\tilde{h} \ \twocommute$};
\draw[->] (AAA) to node [above] {$\sim$} (A);
\draw[->] (AAA) to (B);
\draw[->] (A) to node [left] {$f$} (C);
\draw[->] (BBB) to node [below] {$\sim$} (C);
\draw[->,dashed] (AAA) to node [left] {$(\tilde{f},\gamma)$} (BBB);
\draw[->>] (BBB) to (BB);
\draw[->,bend left=10] (B) to node [above] {$\Phi$} (D);
\draw[->,bend right=10] (BB) to node [below] {$\Phi'$} (D);
\draw[->,dashed ] (B) to node [left] {$(\tilde{f}_0,\gamma)$} (BB);
\path (0,-.75) coordinate (basepoint);
\end{tikzpicture}\label{diag:2pullback}
\eeq
where the squares strictly commute and the triangle commutes up to 2-isomorphism. The datum~$\tilde{h}\in G(S)$ requires a choice of lift as in~\eqref{eq:supertorusmap}. Similarly, the map~$f$ has a (non-unique) lift to a functor determined by the map~$\tilde{f}$ from~\eqref{eq:supertorusmap} together with $\gamma\in \SL_2(\Z)(S)$ when $d=2$. The lift $\tilde{f}$ is only unique up to translations in the lattice specified by~$\ell$. But then the map $\tilde{f}_0$ on the $\E^2$-quotient is uniquely determined by $f$ and $\gamma$, since the ambiguity in the lift $\tilde{f}$ is removed by passing to the quotient by translations (and in particular translations in the lattice). Furthermore, $\tilde{f}_0$ acts on $S\times \R^{0|1}$ by the image of $\tilde{f}$ under the projections
$$
\E^{1|1}\rtimes \R^\times\to \E^{0|1}\rtimes \R^\times, \qquad \E^{2|1}\rtimes \C^\times\to \E^{0|1}\rtimes \C^\times. 
$$
This shows that the target maps factor through the claimed projections. The remaining claims of the lemma follow from the fact that $\Phi$ and $\Phi'$ are related by the action of an $S$-point of $\E^{0|1}\rtimes \R^\times\times G$ when $d=1$ and $\SL_2(\Z)\times \E^{0|1}\rtimes \C^\times\times G$ when $d=2$, as in the 2-commuting triangle in~\eqref{diag:2pullback}. 
\ep

\section{Cocycles for cohomology theories}\label{sec:coccyles}

In this section we construct line bundles $\omega^{\bullet/2+\twist}$ over the super Lie groupoids $\mathcal{L}^{d|1}_0(M\sq G)$ and show that sections of these line bundles give cocycles for $\twist$-twisted equivariant K-theory with complex coefficients when $d=1$ and $\twist$-twisted complex analytic equivariant elliptic cohomology when $d=2$. The statement when $d=2$ requires a restriction to holomorphic sections; see Definition~\ref{defn:holomorphic}. This holomorphic structure is an expected physical consequence of $\mathcal{N}=(0,1)$ supersymmetry (see Remark~\ref{rmk:holophys}) and a mathematical consequence of a more technical treatment of the moduli stack of super tori (see Remark~\ref{rmk:supercircle}). We also provide a sheafy refinement of the cocycle model by way of an isomorphism with the sheaf $\mathcal{K}_G^{\bullet+\twist}(M)$ when $d=1$ and $\Ell_G^{\bullet+\twist}(M)$ when $d=2$. 

\subsection{Sections of line bundles over inertia fields}
The description of the target map in Lemma~\ref{lem:grpdtarget} gives functors between super Lie groupids,
\beq
\mathcal{L}_0^{1|1}(M\sq G)\to \pt\sq \R^\times,\qquad \mathcal{L}_0^{2|1}(M\sq G)\to \pt\sq \C^\times,\label{eq:Hodgemap}
\eeq
determined by the projection homomorphisms, 
$$
\E^{1|1}\rtimes \R^\times\times G\to \R^\times, \qquad \SL_2(\Z)\times \E^{2|1}\rtimes \C^\times\times G\to \C^\times. 
$$
We recall the canonical odd line bundle on $\pt\sq \C^\times$ from Example~\ref{ex:canonicalline}, and we pull this odd line bundle back to $\pt\sq \R^\times$ using the inclusion homomorphism $\R^\times\hookrightarrow \C^\times$. 
\begin{defn} \label{defn:Hodge} Let $\omega^{-k/2}$ denote the pullback of the $k$th tensor power of the canonical odd line bundles along~\eqref{eq:Hodgemap}. 
\end{defn}

\begin{lem} \label{naturalline1}
For a functor of Lie groupoids $N\sq H\to M\sq G$, the pullback of the line bundle $\omega^{\bullet/2}$ on $\mathcal{L}^{d|1}_0(M\sq G)$ along the functor $\mathcal{L}^{d|1}_0(N\sq H)\to \mathcal{L}^{d|1}_0(M\sq G)$ is canonically isomorphic to the line bundle $\omega^{\bullet/2}$ on $\mathcal{L}^{d|1}_0(N\sq H)$. 
\end{lem}

\bp
This follows from the fact that the diagrams of super Lie groupoids
\beq
&&
\begin{tikzpicture}[baseline=(basepoint)];
\node (A) at (0,0) {$\mathcal{L}^{1|1}_0(N\sq H)$};
\node (B) at (3,0) {$\mathcal{L}^{1|1}_0(M\sq G)$};
\node (C) at (1.5,-1) {$\pt\sq \R^\times$};
\draw[->] (A) to (B);
\draw[->] (A) to (C);
\draw[->] (B) to  (C);
\path (0,-.75) coordinate (basepoint);
\end{tikzpicture}
\qquad 
\begin{tikzpicture}[baseline=(basepoint)];
\node (A) at (0,0) {$\mathcal{L}^{2|1}_0(N\sq H)$};
\node (B) at (3,0) {$\mathcal{L}^{2|1}_0(M\sq G)$};
\node (C) at (1.5,-1) {$\pt\sq \C^\times$};
\draw[->] (A) to (B);
\draw[->] (A) to (C);
\draw[->] (B) to  (C);
\path (0,-.75) coordinate (basepoint);
\end{tikzpicture}\nonumber
\eeq
strictly commute. 
\ep

For $d=1,2$ we observe there are functors between super Lie groupoids
\beq
\mathcal{L}_0^{2|1}(M\sq G)\to G^{(2)}\sq (G\times \SL_2(\Z)),\qquad \mathcal{L}_0^{1|1}(M\sq G)\to G\sq G
\eeq
that on objects send the component indexed by $G^{(d)}$ to the associated point of $G^{(d)}$, 
\beq
\coprod_{g\in G}\R_{>0}\times \Map(\R^{0|1},M^g)&\to& \coprod_{g\in G}\pt\simeq G\label{eq:component1}\\
\coprod_{(g_1,g_2)\in G^{(2)}}  \Lat \times \Map(\R^{0|1},M^{g_1,g_2}) &\to& \coprod_{(g_1,g_2)\in G^{(2)}}\pt \simeq G^{(2)}\label{eq:component2}
\eeq
and on morphisms are determined by the projection homomorphisms,
$$
\E^{1|1}\rtimes \R^\times\times G\to G,\qquad \SL_2(\Z)\times \E^{2|1}\rtimes \C^\times\times G\to \SL_2(\Z)\times G. 
$$

\begin{defn} \label{eq:defntwist} Given a $(d+1)$-cocycle $\twist$ on $G$ with values in $U(1)$, let $\mathcal{T}^\twist$ denote the line bundle on $\mathcal{L}_0^{d|1}(M\sq G)$ gotten by pulling back the line bundle also denoted $\mathcal{T}^\alpha$ on $G\sq G$ or $G^{(d)}\sq(G\times \SL_2(\Z))$ constructed in Examples~\ref{1transgr} and~\ref{2transgr}. 
\end{defn}

\begin{lem}  \label{naturalline2}
Suppose we are given a homomorphism of groups $\zeta\colon H\to G$, a $(d+1)$-cocycle $\twist\in Z^{d+1}(G;U(1))$, and a $G$-manifold~$M$. Then the line bundle $\mathcal{T}^\twist$ on $\mathcal{L}^{d|1}_0(M\sq G)$ pulls back to the line bundle $\mathcal{T}^{\zeta^*\alpha}$ on $\mathcal{L}^{d|1}_0(M\sq H)$ for the $H$-action on $M$ through $\zeta\colon H\to G$. 
\end{lem}
\bp
This follows from the fact that the diagrams of super Lie groupoids
\beq
&&
\begin{tikzpicture}[baseline=(basepoint)];
\node (A) at (0,0) {$\mathcal{L}^{1|1}_0(M\sq H)$};
\node (B) at (3,0) {$\mathcal{L}^{1|1}_0(M\sq G)$};
\node (C) at (0,-1) {$H\sq H$};
\node (D) at (3,-1) {$G\sq G$};
\draw[->] (A) to (B);
\draw[->] (A) to (C);
\draw[->] (C) to (D);
\draw[->] (B) to  (D);
\path (0,-.75) coordinate (basepoint);
\end{tikzpicture}
\qquad 
\begin{tikzpicture}[baseline=(basepoint)];
\node (A) at (0,0) {$\mathcal{L}^{2|1}_0(M\sq H)$};
\node (B) at (3.5,0) {$\mathcal{L}^{2|1}_0(M\sq G)$};
\node (C) at (0,-1) {$H^{(2)}\sq (H\times \SL_2(\Z))$};
\node (D) at (3.5,-1) {$G^{(2)}\sq (G\times \SL_2(\Z))$};
\draw[->] (A) to (B);
\draw[->] (A) to (C);
\draw[->] (C) to (D);
\draw[->] (B) to  (D);
\path (0,-.75) coordinate (basepoint);
\end{tikzpicture}\nonumber
\eeq
strictly commute. 
\ep

\subsection{Computing global sections}

We observe that the line bundles $\omega^{k/2+\alpha}:=\omega^{k/2}\otimes \mathcal{T}^\twist$ on $\mathcal{L}_0^{d|1}(M\sq G)$ are trivial line bundles over their supermanifolds of objects, being the tensor product of line bundles that pullback from~$\pt\sq \C^\times$. This implies that global sections of~$\omega^{k/2+\alpha}$ are a subspace of sections of the trivial line over the objects of the super Lie groupoid, i.e., a subspace of functions
\beq
&&\resizebox{.95\textwidth}{!}{$
\begin{array}{l} \displaystyle
\Gamma(\mathcal{L}_0^{1|1}(M\sq G);\omega^{k/2+\twist})\hookrightarrow \prod_{g} C^\infty(\R_{>0}\times \Map(\R^{0|1},M^g))\simeq \prod_{g} \Omega^\bullet(M^g;C^\infty(\R_{>0}))\\\displaystyle
\Gamma(\mathcal{L}_0^{2|1}(M\sq G);\omega^{k/2+\twist})\hookrightarrow \prod_{(g_1,g_2)} C^\infty(\Lat\times \Map(\R^{0|1},M^{g_1,g_2}))\simeq \prod_{(g_1,g_2)} \Omega^\bullet(M^{g_1,g_2};C^\infty(\Lat)).
\end{array}$}
\label{eq:sectonatlas}
\eeq
The subspace of sections is characterized by the equivariance condition~\eqref{eq:grpodsections}, so it remains to compute the action of morphisms on the above functions. For the $\C^\times$-action with $d=2$, it will be useful to consider 
$$
C^\infty_k(\Lat):=\{ f\in C^\infty(\Lat)\mid f(\dil^2\ell_1,\bar\dil^2\bar\ell_1,\dil^2\ell_2,\bar\dil^2\bar\ell_2)=\dil^{-k} f(\ell_1,\bar\ell_1,\ell_2,\bar\ell_2)\},
$$
as the subspace of smooth functions that transform with weight $k/2\in \frac{1}{2}\Z$. We view $\bigoplus_{k\in \Z} C^\infty_k(\Lat)$ as a graded ring with $C^\infty_k(\Lat)$ in degree $-k$. We observe that the multiplication is graded commutative because $C^\infty_k(\Lat)=\{0\}$ for $k$ odd. 

\begin{lem}\label{lem:wtspace}
There is an isomorphism of graded rings
\beq
\bigoplus_{k\in \Z} C^\infty_k(\Lat)\simeq C^\infty(\HH)[\beta,\beta^{-1}],\qquad |\beta|=-2.\label{eq:LattoH}
\eeq
\end{lem}
\bp
Define maps
\beq
C^\infty_{2k}(\Lat)&\to& C^\infty(\HH)\cdot \beta^{k}\qquad f(\ell_1,\bar\ell_1,\ell_2,\bar\ell_2)\mapsto \beta^kf(\ell_1/\ell_2,\bar\ell_1/\bar\ell_2,1,1)\nonumber\\
C^\infty(\HH)\cdot \beta^k&\to&C^\infty_{2k}(\Lat)\qquad F(\tau,\bar\tau)\mapsto \ell_2^{-k}F(\ell_1/\ell_2,\bar\ell_1/\bar\ell_2). 
\eeq
Using that $\ell_2\mapsto \dil^2\ell_2$ under the $\C^\times$-action, the map is well-defined. It is easy to check that it determines isomorphisms of vector spaces in each fixed degree that are compatible with multiplication. The result follows. 
\ep

Next let $\deg\colon \Omega^\bullet(M)\to \Omega^\bullet(M)$ be the derivation on differential forms with $\deg(f)=kf$ for $f\in \Omega^k(M)$. Let $\ell\in C^\infty(\R_{>0})$ denote the standard coordinate function and define
$$
\vol:=(\ell_1\bar\ell_2-\ell_2\bar\ell_1)/2i \in C^\infty(\Lat)
$$ 
as the function that reads off the volume of a torus associated to a given lattice, relative to the flat metric. Then define algebra automorphisms 
$$
\ell^{\deg/2}\colon \Omega^\bullet(M;C^\infty(\R_{>0}))\to \Omega^\bullet(M;C^\infty(\R_{>0})),\quad \vol^{\deg/2}\colon \Omega^\bullet(M;C^\infty(\Lat))\to \Omega^\bullet(M;C^\infty(\Lat))
$$ 
by $\ell^{\deg/2}f=\ell^{k/2}f$ for $f\in \Omega^k(M;C^\infty(\R_{>0}))$ and $\vol^{\deg/2}f=\vol^{k/2}f$ for $f\in \Omega^k(M;C^\infty(\Lat))$. We use that $\ell$ and $\vol$ take values in the positive real numbers to take their (positive) square roots. We adopt the notation
\beq
&&f_g\in \Omega^\bullet(M^g), \qquad f_{g_1,g_2}\in  \Omega^\bullet_\cl(M^{g_1,g_2};C^\infty(\HH)[\beta,\beta^{-1}]),\label{eq:formonsubalg}
\eeq
to denote a differential form on the associated fixed point set. There are natural inclusions
\beq
\Omega^\bullet(M^g)\hookrightarrow \Omega^\bullet(M^g;C^\infty(\R_{>0})),\qquad 
 \Omega^\bullet(M^{g_1,g_2};C^\infty(\HH)[\beta,\beta^{-1}])\hookrightarrow \Omega^\bullet(M^{g_1,g_2};C^\infty(\Lat)).\nonumber
\eeq
The first inclusion is induced by the map on coefficients $\C\hookrightarrow C^\infty(\R_{>0})$ that include constant functions into all functions. The second inclusion is induced by the map on coefficients $\oplus C^\infty_k(\Lat)\hookrightarrow C^\infty(\Lat)$ including functions that transform with weight into all functions together with Lemma~\ref{eq:LattoH}. Finally, let $\Omega^\bullet_\cl(M)\subset \Omega^\bullet(M)$ denote the subalgebra of closed differential forms. 

\begin{prop} \label{prop:Cinfty}
The image of sections of $\omega^{k/2+\alpha}$ under~\eqref{eq:sectonatlas} are the subspaces given by
\beq
\resizebox{\textwidth}{!}{$\begin{array}{l} \displaystyle
\prod_g \bigoplus_{j+l=k}\Omega^{l}_\cl(M^g;\C[\beta,\beta^{-1}]^j)\hookrightarrow  \prod_g \Omega^{\bullet}(M^g;C^\infty(\R_{>0})),  \ \{f_g\}\mapsto \{\ell^{(\deg+k)/2}f_g\},\nonumber\\
\displaystyle \prod_{(g_1,g_2)}\bigoplus_{j+l=k}\Omega^l_\cl(M^{g_1,g_2};C^\infty(\HH)[\beta,\beta^{-1}]^j)\hookrightarrow \prod_{(g_1,g_2)} \Omega^\bullet_\cl(M^{g_1,g_2};C^\infty(\Lat)), \ \{f_{g_1,g_2}\}\mapsto \{\vol^{\deg/2}f_{g_1,g_2}\} \end{array}$}\nonumber
\eeq
where the sets $\{f_g\}_{g\in G}$ and $\{f_{g_1,g_2}\}_{(g_1,g_2)\in G^{(2)}}$ satisfy the further equivariance property
\beq
&&h^*f_g=\mathcal{T}^\twist(h)f_{hgh^{-1}},\quad (h,\gamma)^*f_{g_1,g_2}= \mathcal{T}^\twist(h,\gamma) f_{(h,\gamma)\cdot (g_1,g_2)},\qquad h\in G, \gamma\in \SL_2(\Z)\label{eq:equivforG}
\eeq
where $\mathcal{T}^\twist(h)\in U(1)$ and $\mathcal{T}^\twist(h,\gamma)\in U(1)$ are the values of the 1-cocycles~\eqref{eq:11twist} and~\eqref{eq:21twist} determining the equivariant line bundles $\mathcal{T}^\twist$ on $G\sq G$ and $G^{(2)}\sq (G\times \SL_2(\Z))$, respectively. Finally, $(h,\gamma)^*f_{g_1,g_2}$ is the pullback of $f_{g_1,g_2}$ along the action by $(h,\gamma)$ followed by the action on coefficients $\beta\mapsto \beta/(c\tau+d)$. 
\end{prop}
\begin{proof}
As stated previously, we compute sections using~\eqref{eq:grpodsections}. Using Lemma~\ref{lem:grpdtarget}, this boils down to an equivariance property for the actions of $G$, $\E^{0|1}$, and $\R^\times$ when $d=1$ or $\C^\times\times \SL_2(\Z)$ when $d=2$. 
 
By definition, sections of~$\omega^{\bullet/2+\alpha}$ are invariant under the action of $\E^{0|1}$. As reviewed in Example~\ref{eq:deRham}, this action is generated by (minus) the de~Rham differential, viewed as an odd derivation on the super algebras
$$
C^\infty(\Map(\R^{0|1},M^g))\simeq \Omega^\bullet(M^g),\qquad C^\infty(\Map(\R^{0|1},M^{g_1,g_2}))\simeq \Omega^\bullet(M^{g_1,g_2}).
$$
Hence, the image of~\eqref{eq:sectonatlas} is contained in the subspace of the closed differential forms.  
  
From Example~\ref{eq:deRham} and Definition~\ref{defn:Hodge}, the action of $\dil\in \R^\times$ on $f\in \Omega^k(M^g)$ is by $f\mapsto \dil^{-k}f$. By~\eqref{eq:lataction1}, the $\R^\times$-action on $C^\infty(\R_{>0})$ sends the standard coordinate function $\ell$ to~$\dil^2\ell$. Hence, for $d=1$ the $\R^\times$-equivariance restricts the image of~\eqref{eq:sectonatlas} to the span of elements of the form
$$
\ell^{(\deg+k)/2}f_g\in \bigoplus_{k=l \ {\rm mod} \ 2} \Omega^l(M^g;C^\infty(\R_{>0})), \qquad f\in \bigoplus_{j+l=k} \Omega^l(M^g;(\C[\beta,\beta^{-1}])^{j}), 
$$
where the mod 2 restriction comes from equivariance for the action of $\{\pm 1\}\in \R^\times$, which acts by the parity involution on $\Omega^k(M^g)$ and trivially on $C^\infty(\R_{>0})$. Similarly, the action of $(\dil,\bar\dil)\in \C^\times$ on $f\in \Omega^k(M^g)$ is by $f\mapsto \bar\dil^{-k}f$ and by~\eqref{eq:lataction} the action on $C^\infty(\Lat)$ is $(\ell_1,\bar\ell_1,\ell_2,\bar\ell_2)\mapsto (\dil^2\ell_1,\bar\dil^2\bar\ell_1,\dil^2\ell_2,\bar\dil^2\bar\ell_2)$. Using that $\vol\mapsto u^2\bar u^2\vol$, the $k$th $\C^\times$-weight space is therefore spanned by
$$
\vol^{\deg/2}f_{g_1,g_2}\in \Omega^\bullet(M^{g_1,g_2};C^\infty(\Lat)),\qquad f_{g_1,g_2}\in \bigoplus_{-j+l=k} \Omega^l(M^{g_1,g_2};C^\infty_j(\Lat)).
$$
Using Lemma~\ref{lem:wtspace}, the element $f_{g_1,g_2}$ is equivalent to $f_{g_1,g_2}\in \Omega^\bullet(M^{g_1,g_2};C^\infty(\HH)[\beta,\beta^{-1}])$ of total degree $k$.

It remains to check equivariance for the $G$-action and (for $d=2$) the $\SL_2(\Z)$-action. That the equivariance property is as claimed in~\eqref{eq:equivforG} comes directly from Definition~\ref{eq:defntwist}: the cocycle pulls back from $G\sq G$ when $d=1$ and $G^{(2)}\sq (G\times \SL_2(\Z))$ when $d=2$. The claimed $\SL_2(\Z)$-action on coefficients $C^\infty(\HH)[\beta,\beta^{-1}]$ follows from the definition of the isomorphism in Lemma~\ref{lem:wtspace} and the computation
\beq
\beta=\ell_2^{-1}\mapsto (c\ell_1+d\ell_2)^{-1}=\ell_2^{-1}(c\tau+d)^{-1}=\frac{\beta}{c\tau+d}.\label{eq:coeffaction}
\eeq
Finally, we note that $\ell\in C^\infty(\R_{>0})$ is $G$-invariant (for the case $d=1$), and $\vol\in C^\infty(\Lat)$ is $G\times \SL_2(\Z)$-invariant (for the case $d=2$), so these normalizations to not change the transformation properties of~\eqref{eq:equivforG}. 
\ep

%

\begin{cor}\label{cor:main}
There are isomorphisms of graded vector spaces of smooth sections
\beq
&&\Gamma(\mathcal{L}^{1|1}_0(M\sq G);\omega^{\bullet/2+\twist})\simeq (\prod_{g\in G} \bigoplus_{j+l=\bullet} \Omega^l_\cl(M^g;\C[\beta,\beta^{-1}]^j))^G\nonumber \\
&&\Gamma(\mathcal{L}^{2|1}_0(M\sq G);\omega^{\bullet/2+\twist})\simeq (\prod_{(g_1,g_2)\in G^{(2)}} \bigoplus_{j+l=\bullet} \Omega^l_\cl(M^{g_1,g_2};C^\infty(\HH)[\beta,\beta^{-1}]^j))^{G\times \SL_2(\Z)}\label{cor:smococycle21}
\eeq
where $G$- and $G\times\SL_2(\Z)$-invariants are taken with respect to the actions in the statement of Proposition~\ref{prop:Cinfty}. 
\end{cor}

\bp
This follows immediately from Proposition~\ref{prop:Cinfty}.
\ep

\begin{defn}\label{defn:holomorphic} A section $f\in \Gamma(\mathcal{L}^{d|1}_0(M\sq G);\omega^{\bullet/2+\alpha})$ is \emph{holomorphic} if under the isomorphism~\eqref{cor:smococycle21} it is in the image of the subspace 
$$
(\prod_{(g_1,g_2)\in G^{(2)}} \Omega^\bullet_\cl(M^{g_1,g_2};\mathcal{O}(\HH)[\beta,\beta^{-1}]))^{G\times \SL_2(\Z)}\hookrightarrow \Gamma(\mathcal{L}^{2|1}_0(M\sq G);\omega^{\bullet/2+\twist})
$$
for the standard inclusion $\mathcal{O}(\HH)[\beta,\beta^{-1}]\subset C^\infty(\HH)[\beta,\beta^{-1}]$ of holomorphic functions. 
\end{defn}

\begin{notation}\label{not:holomorphic}
Hereafter, we use the notation $\Gamma(\mathcal{L}^{2|1}_0(M\sq G);\omega^{\bullet/2+\twist})$ for the space of holomorphic sections and $\Gamma_{C^\infty}(\mathcal{L}^{2|1}_0(M\sq G);\omega^{\bullet/2+\twist})$ for the smooth sections. 
\end{notation}

\begin{rmk} \label{rmk:holophys}
Holomorphic dependence of a partition function on the modular parameter~$\tau$ is a consequence of $\mathcal{N}=(0,1)$ supersymmetry. For example, see~\cite[\S4.4]{DualityMock} for the physics argument that the partition function of the $\mathcal{N}=(0,1)$ supersymmetric sigma model with compact target is holomorphic. Hence, we expect field theories to determine functions with this additional property. Alternatively, restricting to holomorphic functions gives the correct coefficients for complex analytic elliptic cohomology. 
\end{rmk}

\begin{proof}[Proof of Theorem~\ref{mainthm}] By sending a closed differential form to its underlying de~Rham cohomology class, the statement is immediate from Theorem~\ref{thm:FHTcomplex}, Proposition~\ref{prop:dgs}, Corollary~\ref{cor:main} and Definition~\ref{defn:holomorphic}. 
\end{proof}

\subsection{A refinement to complexes of sheaves}

Consider the projection homomorphisms
\beq
&&\E^{1|1}\rtimes \R^\times\times G\to \E^{0|1}\times G,\qquad \E^{2|1}\rtimes \C^\times\times G\times \SL_2(\Z)\to \E^{0|1}\times G\times \SL_2(\Z)\label{eq:homomorphisms}
\eeq
and the maps 
\beq
&&\R_{>0}\times \coprod_{h\in G} \Map(\R^{0|1},M^g) \to G,\qquad \Lat \times \coprod_{(g_1,g_2)\in G^{(2)}} \Map(\R^{0|1},M^{g_1,g_2}) \to \HH\times G^{(2)}\label{eq:atlasproj}
\eeq
induced by $\Map(\R^{0|1},M^g)\to \pt$, $\Map(\R^{0|1},M^{g_1,g_2})\to \pt$, and the (coarse) quotients $\R_{>0}/\R^\times\simeq \pt$ and $\Lat\to \Lat/\C^\times\simeq \HH$. The maps~\eqref{eq:homomorphisms} and~\eqref{eq:atlasproj} determine functors between super Lie groupoids, 
\beq
\mathcal{L}^{1|1}_0(M\sq G)&\stackrel{\pi}{\to}& G\sq (G\times\E^{0|1}),\nonumber\\
\mathcal{L}^{2|1}_0(M\sq G)&\stackrel{\pi}{\to}& (\HH\times G^{(2)})\sq (G\times \SL_2(\Z)\times \E^{0|1}).\nonumber
\eeq
This gives the correspondences
\beq\label{eq:correspondences}
\begin{array}{c}
\mathcal{L}^{1|1}_0(M\sq G)\stackrel{\pi}{\to} G\sq (G\times \E^{0|1}) \stackrel{q}{\leftarrow} G\sq G,\\
\mathcal{L}^{2|1}_0(M\sq G)\stackrel{\pi}{\to}(\HH\times G^{(2)})\sq (G\times \SL_2(\Z)\times \E^{0|1}) \stackrel{q}{\leftarrow} \Bun_G(\EE),\end{array}
\eeq
where the $\E^{0|1}$-action in the middle is trivial. 
We recall from Example~\ref{Ex:differentialdescent} that the descent data of a sheaf on $X$ to a sheaf on $X\to X\sq \E^{0|1}$ is the data of an odd, square zero endomorphism of the sheaf. 

\begin{prop}\label{mainpropsheaf}
When $d=1$, the direct image sheaf $\bigoplus_{k\in \Z} \pi_*\omega^{\bullet/2+\alpha}$ on $G$ is
\beq\label{eq:Kthysheafqft}
U\mapsto \prod_{h\in U} \Omega^\bullet(M^g;\C[\beta,\beta^{-1}]), \qquad U\subset G.
\eeq
When $d=2$, the direct image sheaf $\bigoplus_{k\in \Z} \pi_*\omega^{\bullet/2+\alpha}$ on $\HH\times G^{(2)}$ is given by
\beq\label{eq:ellsheafqft}
U\mapsto \prod_{(g_1,g_2)\in p(U)}C^\infty (U_{g_1,g_2};\Omega^\bullet(M^{g_1,g_2})[\beta,\beta^{-1}])\qquad U\subset\HH\times G^{(2)}
\eeq
where (in the notation of Definition~\ref{defn:Ell}) $p(U)\subset G^{(2)}$, and $U=\coprod_{p(U)} U_{g_1,g_2}$. The $G$-equivariant structure on~\eqref{eq:Kthysheafqft} and $G\times \SL_2(\Z)$-equivariant structure on~\eqref{eq:ellsheafqft} obey the same formulas as in Definitions~\ref{defn:Ksheaf} and~\ref{defn:Ell} for $d=1$ and $d=2$, respectively. The $\E^{0|1}$-equivariant structure on is generated by~$-d$, minus the de~Rham differential. 
\end{prop}
\bp
The direct image sheaves along the maps~\eqref{eq:atlasproj} takes $\R^\times$-invariant sections over the coarse $\R^\times$-quotient when $d=1$ and $\C^\times$-invariant sections over the coarse $\C^\times$-quotient when $d=2$. By the calculations in the proof of Proposition~\ref{prop:Cinfty}, this gives the claimed sheaves on the objects. It remains to match the equivariant structures. For the $G$-actions, this is evident since the equivariant structure comes from $\mathcal{T}^\twist$. For the $\SL_2(\Z)$-action, this follows from Lemma~\ref{lem:wtspace} and~\eqref{eq:coeffaction}. Finally, as previously calculated, the $\E^{0|1}$-action is through minus the de~Rham differential, as claimed.
\ep

\begin{defn} \label{defn:sheafholo}
A section of the sheaf $\bigoplus_{k\in \Z} \pi_*\omega^{\bullet/2+\alpha}$ on $(\HH\times G^{(2)})\sq (G\times \SL_2(\Z)\times \E^{0|1})$ is \emph{holomorphic} if it takes values in the subspace, 
$$
\prod_{(g_1,g_2)\in p(U)}\mathcal{O} (U_{g_1,g_2};\Omega^\bullet(M^{g_1,g_2})[\beta,\beta^{-1}])\subset \prod_{(g_1,g_2)\in p(U)}C^\infty (U_{g_1,g_2};\Omega^\bullet(M^{g_1,g_2})[\beta,\beta^{-1}])
$$
using the standard holomorphic structure on $U_{g_1,g_2}\subset \HH$.
\end{defn}

\begin{thm}\label{mainthmsheaf}
Let $\omega^{\bullet/2+\alpha}$ denote the sheaf of smooth sections of the associated line bundle on~$\mathcal{L}^{1|1}_0(M\sq G)$ and the sheaf of holomorphic sections on $\mathcal{L}^{2|1}_0(M\sq G)$. There are isomorphisms of complexes of sheaves
\beq
\pi_*\omega^{\bullet/2+\alpha}\simeq \left\{\begin{array}{ll} \mathcal{K}_G^{\bullet+\twist}(M) & {\rm on} \ G\sq G \ {\rm when} \ d=1 \\ \Ell_G^{\bullet+\twist}(M) & {\rm on} \ \Bun_G(\EE) \ {\rm when} \ d=2\end{array}\right.\label{eq:sheafiso}
\eeq
where sheaf of graded vector spaces is $q^*\pi_*\omega^{\bullet/2+\alpha}$, and descent along $q$ in~\eqref{eq:correspondences} endows this sheaf of graded vector spaces with its differential. When $\twist$ is trivial, the isomorphism~\eqref{eq:sheafiso} is of sheaves of commutative differential graded algebras. 
\end{thm}
\bp
This follows from Definitions~\ref{defn:Ksheaf},~\ref{defn:Ell},~\ref{defn:sheafholo} and Proposition~\ref{mainpropsheaf}. 
\ep

\subsection{Free fermion partition functions as Euler classes} Theorem~\ref{mainthm} allows one to translate examples from physics into cocycle representatives of cohomology classes. One class of such examples comes from the chiral free fermion theories, as we review briefly. 

Let $V$ be a vector space with inner product $\langle-,-\rangle$, and $T_\tau =\C/\tau\Z\oplus\Z$ be a genus~1 Riemann surface. We recall that the classical free fermion theory in dimension $d=1,2$ has as fields the vector spaces of functions 
$$
\mathcal{F}=\left\{\begin{array}{ll} C^\infty(S^1;V) & d=1 \\ C^\infty(T_\tau;V) & d=2 \end{array}\right.\quad 
\mathcal{S}(\phi)=\left\{\begin{array}{ll} \int_{S^1}\langle \phi,i\partial_t\phi\rangle dt & d=1 \\ \int_{T_\tau} \langle \phi,\bar\partial\phi\rangle d^2z & d=2. \end{array}\right.
$$ 
with the indicated classical action functionals.\footnote{More precisely, fields are spinors and the classical action depends on the Dirac operator; the case treated here is for the odd spin structure where the spinor bundle is trivializable.} The \emph{partition function} of the free fermion theory is the determinant of the operator $i\partial_t$ when $d=1$ or $\bar\partial$ when $d=2$. Because $\partial_t$ and $\bar\partial$ have a zero eigenspace (given by constant functions), the determinants are zero and the partition function vanishes. However, a nontrivial partition function arises if one chooses a $G$-representation $\rho\colon G\to \End(V)$ and considers the twist fields (analogous to~\S\ref{sec:twistfields}) 
$$
\mathcal{F}_g=\{f\in C^\infty(\R)\mid f(t+\ell n)=g^nf(t)\}\quad \mathcal{F}_{g_1,g_2}=\{f\in C^\infty(\C)\mid f(z+n\tau+m)=g_1^ng_2^mf(z)\}
$$
depending on an element $g\in G$ when $d=1$ and a pair of commuting elements $g_1,g_2\in G$ when $d=2$. The determinants of $i\partial_t$ acting on $\F_g$ and $\bar\partial$ acting on $F_{g_1,g_2}$ can be rigorously constructed using techniques of $\zeta$-regularization, e.g., see Quillen~\cite{Quillen_det} and Freed~\cite{Freed_det}. The result is a section of a line bundle over the moduli space of flat $G$-bundles on metrized circles when $d=1$ and flat $G$-bundles on elliptic curves when $d=2$. When $d=1$, the determinant section is built from the $\sinh$ function (as in~\eqref{eq:mobius}), while for $d=2$ the determinant is built out of the theta function~\eqref{eq:determinant}, e.g., see~\cite[Equation~18]{Vafatorsion}. This gives identifications
\beq
\det(i\partial_t\otimes V)=\Eu(V)\quad {\rm and} \quad \det(\bar\partial\otimes V)=\Eu(V)\label{eq:Eulerpartition}
\eeq
between the determinant section of the relevant $V$-twisted operator (i.e., the partition function of $V$-valued free fermions with background $G$-symmetry) and the equivariant Euler class of the representation~$V$ valued in complexified equivariant K-theory when $d=1$ and complex analytic equivariant elliptic cohomology when $d=2$. The equality~\eqref{eq:Eulerpartition} is as sections of determinant line bundles on $\R_{>0}\times G\sq G$ when $d=1$ and $\Bun_G(\EE)$ when~$d=2$. The nontriviality of the determinant line bundle is an \emph{anomaly} of the field theory, leading to the interpretation of the twist in equivariant elliptic cohomology as an anomaly. As discussed in~\S\ref{sec:Euler}, we find that the (gauge) anomaly vanishes if and only if the representation has a string structure, as is expected from the physics. We refer to~\cite[\S7]{BET1} for a detailed discussion of these determinant lines and their relationship to free fermion partition functions. 

Many examples in physics arise from manipulating free theories to produce new theories. One such manipulation for the chiral free fermions is called \emph{gauging}. In the next section we show that gauging determines a wrong-way map in equivariant elliptic cohomology. Applying this wrong-way map to equivariant elliptic Euler classes is then identified with gauging the free fermion theories. Following several well-studied examples in the physics literature then leads to interesting examples of pushforwards in equivariant elliptic cohomology. 

\section{Gauging and discrete torsion} \label{sec:discretetorsion}

%


\subsection{Motivation from path integral quantization of finite gauge theories}
Quantizing the classical theory reviewed in~\S\ref{sec:twistfields} happens in two steps, namely integration along the fibers of the two maps
\beq
\mathcal{L}^{d|1}(M\sq G)\to \mathcal{L}^{d|1}(\pt\sq G)\to \mathcal{L}^{d|1}(\pt), \label{eq:pushforwardforquant}
\eeq
where the first arrow is induced by the $G$-equivariant map $M\to \pt$, and the second arrow is from the homomorphism $G\to *$. Integration along the first map uses the techniques of supersymmetric localization briefly described in~\S\ref{sec:twistfields}. Integration along the second map is quantization of a gauge theory. This can be implemented rigorously via the path integral: since $G$ is finite, the path integral is a finite sum weighted by automorphism groups~\cite{DW,FreedQuinn}. Indeed, this integral follows the paradigm of integration on finite groupoids as reviewed in~\S\ref{sec:integration}. We show that these integrals construct the expected transfer maps in equivariant elliptic cohomology. 

A cocycle $B\in Z^2(G;U(1))$ allows one to modify these transfer maps. Physically this modification comes from a choice of $B$-field~\cite{sharpetorsion,sharpetorsionhet}. From a computational perspective, the effect is to modify the contribution to the integral of each twisted sector by a phase depending on~$B$. Vafa calls this \emph{discrete torsion}~\cite{Vafatorsion}. When applied to equivariant elliptic Euler classes, this results in a family of pushforwards depending on the choice of string structure on the representation. A similar effect was studied in~\cite{AndoFrench}, albeit in the different context of Borel equivariant elliptic cohomology at a super singular elliptic curve. 


\subsection{Wrong-way maps and discrete torsion}

We construct a wrong-way map in a more general context than integration along the second arrow in~\eqref{eq:pushforwardforquant}. Namely, given a $G$-manifold $M$ and a homomorphism $\zeta\colon  H\to G$, the induced functor between Lie groupoids $M\sq H\to M\sq G$ determines a functor 
\beq
\mathcal{L}_0^{d|1}(M\sq H)\to \mathcal{L}_0^{d|1}(M\sq G)\label{eq:transfer}
\eeq 
using Corollary~\ref{cor:natural}. This functor sends an $H$-bundle to the associated $G$-bundle using~$\zeta\colon H\to G$ to define a left $H$-action on $G$. We further recall from Lemmas~\ref{naturalline1} and~\ref{naturalline2} that for $\alpha\in Z^{d+1}(G;U(1))$, the line bundle $\omega^{\bullet/2+\twist}$ on $\mathcal{L}_0^{d|1}(M\sq G)$ pulls back along~\eqref{eq:transfer} to the line bundle $\omega^{\bullet/2+\zeta^*\twist}$ on $\mathcal{L}_0^{d|1}(M\sq H)$. Next, recall the maps of super Lie groupoids
\beq
\mathcal{L}_0^{1|1}(M\sq H)\to H\sq H,\qquad \mathcal{L}_0^{2|1}(M\sq H)\to H^{(2)}\sq (H\times \SL_2(\Z))\label{eq:projectsomewhere}
\eeq
defined using~\eqref{eq:component1} and~\eqref{eq:component2}. Given a $d$-cocycle $B\in Z^d(H;U(1))$, transgression from Example~\ref{1transgr} and Example~\ref{2transgr} determines a function on $H\sq H$ when $d=1$ and $H^{(2)}\sq (H\times \SL_2(\Z))$ when $d=2$. Define $\epsilon_B\in C^\infty(\mathcal{L}_0^{1|1}(M\sq H))$ as the pullback of this function along~\eqref{eq:projectsomewhere}. 

\begin{defn}
Given $\zeta$ and $B$ as above and a section $f\in \Gamma (\mathcal{L}^{d|1}_0(M\sq H);\omega^{\bullet/2+\zeta^*\twist})$, define a function $\zeta_!^Bf$ on the objects of $\mathcal{L}^{d|1}_0(M\sq G)$ whose value at an $S$-point $(\ell,{\bf g},\phi)$ (using the notation of Definition~\ref{defn:fields} and~\eqref{diag:assocmap}) is 
\beq
&&(\zeta_!^Bf)(\ell,{\bf g},\phi)=\sum_{\begin{smallmatrix} [P_{\bf g}\simeq P_{\bf h}\times_\zeta G]\\ \phi=\phi_\zeta'\end{smallmatrix}} \frac{(\epsilon_B\cdot f)(\ell,{\bf h},\phi')}{|\Aut(P_{\bf h})|/|\Aut(P_{\bf g})|}\in C^\infty({\rm Obj}(\mathcal{L}_0^{d|1}(M\sq G)) \label{eq:induction2}
\eeq
where the sum is indexed by isomorphism classes of principal $H$-bundles $P_{\bf h}$ such that there exists an isomorphism of $G$-bundles $P_{\bf g}\simeq P_{\bf h}\times_\zeta G$ compatible with the maps to $M$.
%
\end{defn}

\begin{rmk}
We observe that the denominator in~\eqref{eq:induction2} is the fiberwise volume form for the map of finite groupoids~\eqref{eq:inducedcharacterd} pulled back along the map
$$
\mathcal{L}_0^{d|1}(M\sq H)\to H^{(d)}\sq H. 
$$
This fiberwise volume form is modified by the choice of discrete torsion~$\epsilon_B$ to define~\eqref{eq:induction2}. 
\end{rmk}

\begin{thm} \label{prop:pushforward}
The assignment~\eqref{eq:induction2} determines maps for $d=1,2$
\beq
\zeta_!^B\colon \Gamma(\mathcal{L}^{d|1}_0(M\sq H);\omega^{\bullet/2+\zeta^*\twist})\to \Gamma(\mathcal{L}^{d|1}_0(M\sq G);\omega^{\bullet/2+\twist})\label{eq:inductionmapsagain}
\eeq
parameterized by the choice of $[B]\in \H^d(BG;U(1))$. 
\end{thm}

\begin{proof}
The proposition amounts to showing that if $f$ is invariant under~$\E^{0|1}$ and equivariant for $\R^\times\times H$ when $d=1$ and invariant under $\E^{0|1}$ and equivariant for $\C^\times \times H\times \SL_2(\Z)$ when $d=2$, then $\zeta_!^Bf$ is similarly $\E^{0|1}$-invariant and $\R^\times\times G$- respectively, $\C^\times \times G\times \SL_2(\Z)$-equivariant. To show this, it is convenient to rewrite~\eqref{eq:induction2} using the orbit stabilizer argument in the proof of Lemma~\ref{lem:Frob},
\beq
(\zeta_!^Bf)(\ell,{\bf g},\phi)=\frac{1}{|H|}\sum_{\begin{smallmatrix} P_{\bf g}\simeq P_{\bf h}\times_\zeta G\\ \phi=\phi_\zeta'\end{smallmatrix}} (\epsilon_B\cdot f)(\ell,{\bf h},\phi')\in C^\infty({\rm Obj}(\mathcal{L}_0(M\sq G)) \label{eq:induction3}
\eeq
where the sum in~\eqref{eq:induction3} is indexed by the set of pairs given by a principal $H$-bundle $P_{\bf h}$ for ${\bf h}\in H^{(2)}$ and an isomorphism $P_{\bf g}\simeq P_{\bf h}\times_\zeta G$ from $x\in G$, where the isomorphism is required to be compatible with the $G$-equivariant maps to $M$, $\phi\colon P_{\bf g}\to M$ and $\phi'_\zeta\colon P_{\bf h}\times_\zeta G\to M$.

Using Proposition~\ref{prop:Cinfty}, we express a section of $\omega^{\bullet/2+\twist}$ in terms of a set of differential forms $\{f_g\}$ when $d=1$ and $\{f_{g_1,g_2}\}$ when $d=2$, where
$$
f_g\in \Omega^\bullet(M^g;\C[\beta,\beta^{-1}]),\qquad f_{g_1,g_2}\in \Omega^\bullet(M^{g_1,g_2};\mathcal{O}(\HH)[\beta,\beta^{-1}]).
$$
In terms of these data we can re-express~\eqref{eq:induction2} as (compare Proposition~\ref{prop:HKR})
\beq
 (\zeta_!^Bf)_g&=&\displaystyle\frac{1}{|H|}\sum_{g=x\zeta(h)x^{-1}} \epsilon_B(h)x^*f_h \label{eq:HKRagain1}\\
 (\zeta_!^Bf)_{g_1,g_2}&=&\displaystyle\frac{1}{|H|}\sum_{\begin{smallmatrix} (g_1,g_2)= \\ x\zeta(h_1,h_2)x^{-1}\end{smallmatrix}} \epsilon_B(h_1,h_2)x^*f_{h_1,h_2} \label{eq:HKRagain2}
\eeq
where the first sum is indexed by pairs $(h,x)\in H\times G$ such that $g=x\zeta(h)x^{-1}$ and the second sum is indexed by $(h_1,h_2,x)\in H^{(2)}\times G$ such that 
$$
(g_1,g_2)=(x\zeta(h_1)x^{-1},x\zeta(h_2)x^{-1})=:x\zeta(h_1,h_2)x^{-1}.
$$
By definition, $H$ acts on $M$ through $\zeta\colon H\to G$, so if $\zeta(h)=g$ or $(\zeta(h_1),\zeta(h_2))=(g_1,g_2)$ we have 
\beq
\Omega^\bullet(M^{h};\C[\beta,\beta^{-1}])&=&\Omega^\bullet(M^g;\C[\beta,\beta^{-1}]),\nonumber\\
\Omega^\bullet(M^{h_1,h_2};\mathcal{O}(\HH)[\beta,\beta^{-1}])&=&\Omega^\bullet(M^{g_1,g_2};\mathcal{O}(\HH)[\beta,\beta^{-1}]),\nonumber
\eeq
and similarly if $x\zeta(h)x^{-1}=g$ or $x(\zeta(h_1),\zeta(h_2))x^{-1}=(g_1,g_2)$, then $x^*f_h$ and $x^*f_{h_1,h_2}$ give differential forms on the $g$-fixed points and $(g_1,g_2)$-fixed points, respectively. Hence the sums in~\eqref{eq:HKRagain1} and~\eqref{eq:HKRagain2} give differential forms on the correct fixed point sets. 

In this description, we identify $\E^{0|1}$-invariant sections with closed differential forms. Then $\E^{0|1}$-invariance of the sum is clear: a sum of closed forms is closed. We also observe that the sum of exact forms is exact, so a coboundary is sent to a coboundary. The $\R^\times$ or $\C^\times$-equivariant property follows similarly: a sum of degree~$k$ differential forms is a degree~$k$ differential form. The $G\times \SL_2(\Z)$-equivariance follows from $\epsilon_B$ being a function, and hence \emph{invariant} under this action. Hence, the transformation properties of a section of $\zeta^*\alpha$ are the same as for an $\alpha$-twisted section, and the sum of differential forms with this transformation property inherits the same equivariance property. This completes the proof. 
\ep

\begin{cor}\label{cor:inductionexists}
The maps~\eqref{eq:inductionmapsagain} determine wrong-way maps 
\beq
&&\zeta_!^B\colon \K_H^{\bullet+\zeta^*\alpha}(M)\to \K_G^{\bullet+\alpha}(M)\qquad \zeta_!^B\colon \TMF_H^{\bullet+\zeta^*\alpha}(M)\otimes \C\to \TMF_G^{\bullet+\alpha}(M)\otimes \C. \label{eq:transfer2}
\eeq
\end{cor}

\bp
This is an immediate consequence of Theorem~\ref{mainthm}, using that~\eqref{eq:inductionmapsagain} sends cocycles to cocycles and coboundaries to coboundaries, as seen in the last paragraph of the proof of Theorem~\ref{prop:pushforward}. 
\ep
%

\begin{cor}\label{cor:HKR}
For an inclusion of groups $\zeta\colon H\hookrightarrow G$, the map determined by~\eqref{eq:induction2} with $M=\pt$ and $B\equiv 1$ recovers the character formula for an induced representation when~$d=1$, and the height~2 Hopkins--Kuhn--Ravenel character formula when $d=2$. 
\end{cor}

\bp
This follows immediately from specializing~\eqref{eq:HKRagain1} and~\eqref{eq:HKRagain2} to this case together with the descriptions of these character formulas as in Proposition~\ref{prop:HKR}. 
\ep

\begin{cor}\label{cor:Vafa}
When $M=\pt$, $d=2$, and $\zeta\colon G\to  \{e\}$ is the terminal homomorphism,~\eqref{eq:induction2} recovers the formula for gauging with discrete torsion $\epsilon_B$ from Vafa~\cite{Vafatorsion}. 
\end{cor}

\bp
In this case~\eqref{eq:induction2} simplifies to 
\beq
(\zeta_!^Bf)(\ell,(g_1,g_2))&=&\frac{1}{|G|}\sum_{(g_1,g_2)\in G^{(2)}} \epsilon_B(g_1,g_2)f(\ell,(g_1,g_2))\nonumber\\
&=&\frac{1}{|G|}\sum_{g_1g_2=g_2g_1}\frac{B(g_1,g_2)}{B(g_2,g_1)} f(\ell,(g_1,g_2)),\nonumber
\eeq
which recovers Vafa's discrete torsion, as can be extracted from~\cite[Equation~51]{Vafatorsion}; see also~\cite[Equation~1.7]{AndoFrench}, though beware the misprint in their formula:  one must pair the 2-cocycle on $G$ with the fundamental class of the torus to reproduce the correct invariant (see~\cite[pg.~605 footnote]{Vafatorsion}). 
\ep

\subsection{Invariants of representations from gauging and discrete torsion}\label{sec:examples}

The following example connects gauging and discrete torsion with classical results for the character ring. 
\begin{ex}\label{ex:Kthyfixed}
We compute the wrong-way map applied to a character $\chi(\rho)\in C^\infty(G\sq G)$ of a representation $\rho\colon G\to \End(V)$. Identifying the character of the representation with a function $\chi(\rho)\in C^\infty(\mathcal{L}^{1|1}_0(\pt\sq G))$, the wrong-way map applied to $\zeta\colon G\to *$ is
\beq
\zeta_!(\chi(\rho))=\sum_{g\in G} {\rm Tr}(\rho(g))={\rm dim}(V^G)\in C^\infty(\mathcal{L}^{1|1}(\pt)),\label{eq:Ginvt}
\eeq
and hence records the dimension of the $G$-fixed subspace of $V$. Indeed, the above sum is the formula for the inner product $\langle \chi,{\bf 1}\rangle$ on the character ring where ${\bf 1}$ is the trivial representation. A choice of discrete torsion is a 1-dimensional representation $B\colon G\to U(1)$. Hence, $\zeta_!^B(\chi(\rho))$ is the dimension of the subspace of $V$ on which $G$ acts by~$\bar B\colon G\to U(1)$. 
\end{ex}

\begin{rmk}
In physics, the goal of gauging a theory with $G$-symmetry is to produce a $G$-invariant field theory: the dependence on $G$ is ``integrated out." The relationship between the gauged theory and the theory with (background) $G$-symmetry is analogous to the relationship between a representation and its $G$-invariant subspace, with the sum~\eqref{eq:Ginvt} being an easy version of the path integral. 
\end{rmk}

Next we apply wrong-way maps to the $(\Z/n)^{\times k}$-equivariant Euler classes from Examples~\ref{eq:Kex2} and~\ref{ex:Klein}. From the discussion there, if a representation has a spin or string structure the Euler class is canonically untwisted, and the untwisting gives a canonical choice of spin or string structure in these respective cases. With this canonical structure fixed, we then may identify any other spin structure with a 1-cocycle~$B\in Z^1(BG;\Z/2)$, and similarly any other string structure with a 2-cocycle~$B\in Z^2(BG;U(1))$. Let $s$ denote the map sending a spin or a string structure to the associated cocycle in $Z^1(BG;\Z/2)$, respectively, $Z^2(BG;U(1)$. Applying~\eqref{eq:transfer2} to the Euler class of a representation with spin structure, respectively a representation with string structure, therefore gives maps
\beq
&&\resizebox{.95\textwidth}{!}{$
\begin{array}{l}
\{ G-{\rm representations \ with \ spin \ structure} \}\stackrel{s\times \Eu}{\to}Z^1(BG;\Z/2)\times \K_{G}^\bullet(\pt)\stackrel{\zeta_!^B}{\to} \K^\bullet(\pt)\\
\{ G-{\rm representations \ with \ string \ structure} \}\stackrel{s\times \Eu}{\to}Z^2(BG;U(1))\times\Ell_{G}^\bullet(\pt)\stackrel{\zeta_!^B}{\to} \Ell^\bullet(\pt)\end{array}$}
\eeq
where $G=(\Z/n)^{\times k}$ and $\zeta_!^B$ denotes the wrong-way map relative to spin or string structure~$B$. The resulting element of $\K^\bullet(\pt)\simeq \C[\beta,\beta^{-1}]$, respectively, $\Ell^\bullet(\pt)\simeq (\mathcal{O}(\HH)[\beta,\beta^{-1}])^{\SL_2(\Z)}$ is an invariant of the representation with its choice of spin structure, respectively, string structure. We illustrate the non-triviality of this construction in the following examples.

\begin{ex} 
As before, let $\sigma$ denote the sign representation of $\Z/2$, i.e., the action on $\R$ by $\{\pm 1\}$. From Example~\ref{ex:Kthyz2}, $\sigma^{\oplus 4}=\mu_2^{\oplus 2}$ has a spin structure and a canonically untwisted Euler class. For this Euler class we compute
$$
\zeta_!(\Eu(\sigma^{\oplus 4}))=\frac{1}{2}(e^{\pi i/2}-e^{-\pi i/2})^2\beta^{-1}=\frac{1}{2}(2i)^2\beta^{-1}=-2\beta^{-1}\in \K^2(\pt). 
$$
The choice of nontrivial class $[B]\in \H^1(\Z/2;\Z/2)\simeq \Z/2$ gives 
$$
\zeta_!^B(\Eu(\sigma\oplus\sigma))=-\frac{1}{2}(e^{\pi i/2}-e^{-\pi i/2})^2\beta^{-1}=-\frac{1}{2}(2i)^2=2\beta^{-1}\in \K^2(\pt),
$$
and hence changes the value by a sign.
\end{ex}

Next we consider two examples of classes pushed forward along
$$
\zeta_!\colon \Gamma(\Bun_{\Z/2};\Ell_{\Z/2}^\bullet(\pt))\to \Gamma(\Mell;\Ell^\bullet(\pt))\simeq (\mathcal{O}(\HH)[\beta,\beta^{-1}])^{\SL_2(\Z)}
$$
for the homomorphism $\Z/2\to *$ and $B\equiv 1$. One example shows that a nontrivial class can pushforward to~0, while the other shows that infinitely many distinct classes are in the image.

\begin{ex}\label{ex:Jacobi}
The Jacobi's identity for theta functions is
$$
\vartheta^4_3=\vartheta_2^2+\vartheta_4^4. 
$$
Hence, applying~\eqref{eq:induction2} to the Euler class $\Eu(\sigma^{\oplus 8})\in \Gamma(\Bun_{\Z/2}(\EE);\Ell^8_{\Z/2}(\pt))$ from Example~\ref{eq:Z2}, we find
$$
\zeta_! (\Eu(\sigma^{\oplus 8}))=\zeta_! (\beta^{-4} Z(\tau,u,v)^4)=\frac{\beta^{-4}}{\eta^{12}}(\vartheta_2^4-\vartheta_3^4+\vartheta_4^4)=0\in \Gamma(\Mell;\Ell^8(\pt)). 
$$
This example is well-known in physics, arising in the computation of the partition function for 8 chiral free fermions, e.g., see~\cite[pg.~105]{Ginsparg1988} or~\cite[\S2]{Theta}.
\end{ex}

\begin{ex}[The $E_8$-theta function]
We recall that the theta function of a rank~$n$ even unimodular lattice $\Lambda$ is the modular form of weight $n/2$ defined by
$$
\vartheta_\Lambda=\sum_{{\bf x}\in \Lambda} e^{i\pi \tau \|{\bf x}\|^2}. 
$$
The $E_8$ lattice $\Lambda_{8}$ has rank $8$, and its $k$th direct sum with itself $\Lambda_{8k}:=\Lambda_8^{\oplus k}$ therefore has rank~$8k$. The associated theta function can be expressed in terms of Jacobi theta functions
$$
\vartheta_{\Lambda_{8k}}=\frac{1}{2}(\vartheta_2^{8k}+\vartheta_3^{8k}+\vartheta_4^{8k}).
$$
Hence, for $\Eu(\sigma^{\oplus 16k})\in \Gamma(\Bun_{\Z/2}(\EE);\Ell^{16k}_{\Z/2}(\pt))$, we have
\beq
\zeta_!(\Eu(\sigma^{\oplus 16k}))&=&\zeta_!(\beta^{-8k}Z(\tau,u,v)^{8k})=\frac{\beta^{-8k}}{2\eta^{24k}}(\vartheta_2^{8k}+\vartheta_3^{8k}+\vartheta_4^{8k})\nonumber\\
&=&\beta^{-8k}\frac{\vartheta_{\Lambda_{8k}}}{\Delta^k} \in \Gamma(\Mell;\Ell^{16k}(\pt)) \label{eq:E8}
\eeq
where $\Delta$ is the modular discriminant. Up to normalization factors by $\eta$,~\eqref{eq:E8} is the partition function of $16k$ chiral free fermions~\cite[pg.~106]{Ginsparg1988}. With $k=1$ we get
\beq
\zeta_!((\beta^{12}\Delta)\cdot \Eu(\sigma^{\oplus 8}))=\beta^4\vartheta_{\Lambda_{8}}=E_4\in \Gamma(\Mell;\Ell^{-8}(\pt)),\label{eq:E8theta}
\eeq
where $E_4=\vartheta_{\Lambda_{8}}$ is the $4$th Eisenstein series. 
\end{ex}


\begin{ex}[Klein's $j$-function]
Consider the $(\Z/2)^{\times 3}$-equivariant class
\beq
&&\resizebox{.92\textwidth}{!}{$
\Eu(\sigma^{\oplus 16})^{\boxtimes 3}=\beta^{-24}(Z(\tau,u,v)^8\boxtimes Z(\tau,u,v)^8\boxtimes Z(\tau,u,v)^8)\in \Gamma(\Bun_{(\Z/2)^{\times 3}}(\EE);\Ell_{(\Z/2)^{\times 3}}^{48}(\pt))$}.\label{eq:preimageofj}
\eeq
We recall the identity 
$$
j=\frac{1}{8\Delta}(\vartheta_2^8+\vartheta_3^8+\vartheta_4^8)^3=\frac{1}{8\eta^{24}}(\vartheta_2^8+\vartheta_3^8+\vartheta_4^8)^3
$$
for Klein's $j$-function, and hence the pushforward along $\zeta\colon (\Z/2)^{\times 3}\to *$ is 
\beq
\zeta_!(\Eu(\sigma^{\oplus 16})^{\boxtimes 3})&=&\zeta_!(\beta^{-24}Z(\tau,u,v)^8\boxtimes Z(\tau,u,v)^8\boxtimes Z(\tau,u,v)^8)\nonumber\\
&=&\frac{\beta^{-24}}{8\eta^{72}}(\vartheta_2^8+\vartheta_3^8+\vartheta_4^8)^3=\frac{j}{\beta^{24}\Delta^2}\in \Gamma(\Mell;\Ell^{48}(\pt)).\label{eq:jdelta}
\eeq
In particular, 
\beq
\zeta_!((\beta^{24}\Delta^2)\cdot \Eu(\sigma^{\oplus 16})^{\boxtimes 3})=j\label{eq:jdelta2}
\eeq
recovering the $j$-function in the image of~\eqref{eq:preimageofj}. 
\end{ex}

\begin{rmk}\label{rmk:modularintlift}
We recall that the modular form $\Delta$ does not lift to a topological modular form; indeed $24\Delta$ and $\Delta^{24}$ are the smallest sum and product of $\Delta$ that lift~\cite[\S4.3]{HopkinsICM2002}. As such, we would expect the computations~\eqref{eq:E8} and~\eqref{eq:jdelta} to lift to integral statements (as equivariant elliptic Euler classes ought to lift to classes over $\Z$), but~\eqref{eq:E8theta} and~\eqref{eq:jdelta2} need to by modified owing to the nonexistence of the class $\Delta$. For example, one might hope that $\zeta_!((\beta^{24}(24\Delta)^2)\cdot \Eu(\sigma^{\oplus 16})^{\boxtimes 3})=576\cdot j$ has an integral lift. 
\end{rmk}

%
%
%
%
%

%


\begin{ex}[Discrete torsion for $\Z/2\times \Z/2$]\label{ex:discretetorsion}
Since $\H^2(\Z/2\times \Z/2;U(1))\simeq \Z/2$, the Klein 4-group is the easiest example in which we might explore nontrivial choices of discrete torsion. For $a,b,a',b'\in \{0,1\}\simeq \Z/2$, a choice of generator of $\H^2(\Z/2\times \Z/2;U(1))$ is given by the 2-cocycle
\beq
B((a,b),(a',b'))=(-1)^{ab'}\in \Z/2\simeq \{\pm 1\}\subset U(1)\label{eq:distorsionklein}
\eeq
for $(a,b),(a',b')\in \Z/2\times \Z/2$. So we find that 
$$
\epsilon_B((a,b),(a',b'))=\frac{B((a,b),(a',b'))}{B((a',b'),(a,b))}= (-1)^{ab'-ba'}\qquad \begin{tabular}{c|cccc}
\null &$(0,0)$ & $(0,1)$ & $(1,0)$ & $(1,1)$\\
\hline
$(0,0)$ & + & + & + & + \\
$(0,1)$ & + & $+$ & $-$ & $-$ \\ 
$(1,0)$ & + & $-$ & $+$ & $-$ \\ 
$(1,1)$ & + & $-$ & $-$ &$+$ 
\end{tabular}
$$
where the table on the right provides a way of organizing the signs for pairs of elements $(a,b),(a',b')\in \Z/2\times \Z/2$; the first row lists the possible choices for $(a,a')$ and the first column lists the choices for $(b,b')$, with the entries on the interior recording the sign associated with $((a,b),(a',b'))$.
%
Now consider the $\Z/2\times\Z/2$-equivariant class
$$
\Eu(\sigma^{\oplus 8})^{\boxtimes 2}=\beta^{-8}Z(\tau,u,v)^4\boxtimes Z(\tau,u,v)^4 \in \Gamma(\Bun_{\Z/2\times \Z/2}(\EE);\Ell^{16}_{\Z/2\times \Z/2}(\pt)). 
$$
%
%
where we have a similar table describing this equivariant class as a function on $\HH$ for any pair of commuting elements $(a,b),(a',b')$ in $\Z/2\times \Z/2$
\beq
\begin{tabular}{c|cccc}
\null &$(0,0)$ & $(0,1)$ & $(1,0)$ & $(1,1)$\\
\hline
$(0,0)$ & 0 & 0 & 0 & 0 \\
$(0,1)$ & 0 & $\vartheta_2^8/\Delta$ & $-\vartheta_2^4\vartheta_3^4/\Delta$ & $\vartheta_2^4\vartheta_4^4/\Delta$ \\ 
$(1,0)$ & 0 & $-\vartheta_2^4\vartheta_3^4/\Delta$ & $\vartheta_3^8/\Delta$ & $-\vartheta_3^4\vartheta_4^4/\Delta$ \\ 
$(1,1)$ & 0 & $\vartheta_2^4\vartheta_4^4/\Delta$ & $-\vartheta_3^4\vartheta_4^4/\Delta$ &$\vartheta_4^8/\Delta$ 
\end{tabular}\label{eq:table}
\eeq
Then for $\zeta\colon \Z/2\times \Z/2\to *$, we compute
\beq
\zeta_!(\Eu(\sigma^{\oplus 8})^{\boxtimes 2})&=&\zeta_!(\beta^{-8}Z(\tau,u,v)^4\boxtimes Z(\tau,u,v)^4))\nonumber\\
&=&\frac{\beta^{-8}}{4\Delta}\Big(\vartheta_2^8-\vartheta_2^4\vartheta_3^4+\vartheta_2^4\vartheta_4^4-\vartheta_2^4\vartheta_3^4+\vartheta_3^8-\vartheta_3^4\vartheta_4^4+\vartheta_2^4\vartheta_4^4-\vartheta_3^4\vartheta_4^4+\vartheta_4^8\Big)\nonumber\\ 
&=&\frac{1}{4\beta^8\Delta}(\vartheta_2^4-\vartheta_3^4+\vartheta_4^4)^2\nonumber\\
&=&0\in \Gamma(\Mell;\Ell^{16}(\pt)).\label{eq:Klein11}
\eeq
On the other hand, with $B$ from~\eqref{eq:distorsionklein} we calculate
\beq
\zeta_!^B(\Eu(\sigma^{\oplus 8})^{\boxtimes 2})&=&\zeta_!^B(\Eu(\sigma^{\oplus 4})^{\boxtimes 2})+\zeta_!(\Eu(\sigma^{\oplus 4})^{\boxtimes 2})=\frac{\beta^{-8}}{4\Delta}(2(\vartheta_2^8+\vartheta_3^8+\vartheta_4^8))\nonumber\\
&=&\frac{\vartheta_{E_8}}{\beta^8\Delta}\in \Gamma(\Mell;\Ell^{16}(\pt)),\label{eq:Klein21}
\eeq
where addition of $\zeta_!(\Eu(\sigma^{\oplus 8})^{\boxtimes 2})$ (which equals zero) cancels all the off-diagonal terms and doubles the diagonal terms in the description from the table~\eqref{eq:table}. Hence,~\eqref{eq:Klein11} and~\eqref{eq:Klein21} give the invariants of the $\Z/2\times \Z/2$-representation $(\sigma^{\oplus 8})^{\boxtimes 2}$ on $\R^{16}$ that depends on the two possible choices of string structure. The string structure with $B=0$ comes from the product of the canonical string structures on~$\R^8$ discussed in Example~\ref{eq:Z2}, and the $B\ne 0$ choice gives the other isomorphism class of string structure. 
\end{ex}


\appendix

\bibliographystyle{amsalpha}
\bibliography{references}

\end{document}